\input amstexl


\catcode`\@=11
\ifx\amstexloaded@\relax\else
 \errmessage{AmS-TeX must be loaded before LamS-TeX}\fi
\ifx\laxread@\undefined\else\catcode`\@=\active \fi
\def\err@#1{\errmessage{LamS-TeX error: #1}}
\def^^L{\par}
\let\+\tabalign
\def\newcount{\alloc@0\count\countdef\insc@unt}
\def\newdimen{\alloc@1\dimen\dimendef\insc@unt}
\def\newskip{\alloc@2\skip\skipdef\insc@unt}
\def\newmuskip{\alloc@3\muskip\muskipdef\@cclvi}
\def\newbox{\alloc@4\box\chardef\insc@unt}
\let\newtoks\relax
\def\newhelp#1#2{\newtoks#1#1\expandafter{\csname#2\endcsname}}
\def\newtoks{\alloc@5\toks\toksdef\@cclvi}
\def\newread{\alloc@6\read\chardef\sixt@@n}
\def\newwrite{\alloc@7\write\chardef\sixt@@n}
\def\newfam{\alloc@8\fam\chardef\sixt@@n}
\def\newlanguage{\alloc@9\language\chardef\@cclvi}
\def\newinsert#1{\global\advance\insc@unt by\m@ne
  \ch@ck0\insc@unt\count
  \ch@ck1\insc@unt\dimen
  \ch@ck2\insc@unt\skip
  \ch@ck4\insc@unt\box
  \allocationnumber=\insc@unt
  \global\chardef#1=\allocationnumber
  \wlog{\string#1=\string\insert\the\allocationnumber}}
\def\newif#1{\count@\escapechar \escapechar\m@ne
  \expandafter\expandafter\expandafter
   \edef\@if#1{true}{\let\noexpand#1=\noexpand\iftrue}%
  \expandafter\expandafter\expandafter
   \edef\@if#1{false}{\let\noexpand#1=\noexpand\iffalse}%
  \@if#1{false}\escapechar\count@}

\def\Err@#1{\errhelp\defaulthelp@\err@{#1}}
{\catcode`\@=\active
 \edef\next{\gdef\noexpand@{\futurelet\noexpand\next
  \csname at\string@\endcsname}}
 \next
}
\def\at@{\ifcat\noexpand\next a\let\next@\at@@\else
 \ifcat\noexpand\next0\let\next@\at@@\else
 \ifcat\noexpand\next\relax\let\next@\at@@\else
 \let\next@\at@@@\fi\fi\fi\next@}
\def\at@@@{\errhelp\athelp@\err@{Invalid use of @}}
\def\at@@#1{\expandafter
 \ifx\csname\string#1@at\endcsname\relax\let\next@\at@@@\else
 \DN@{\csname\string#1@at\endcsname}\fi\next@}
\def\atdef@#1{\expandafter\def\csname\string#1@at\endcsname}
\newif\iftest@
\def\tagin@#1{\tagin@false
 \DN@##1\tag##2##3\next@{\test@true\ifx\tagin@##2\test@false\fi}%
 \next@#1\tag\tagin@\next@\tagin@false\iftest@\tagin@true\fi}
\let\lkerns@\relax
\def\nolinebreak{\RIfM@\mathmodeerr@\nolinebreak\else
 \ifhmode\saveskip@\lastskip\unskip
 \nobreak\ifdim\saveskip@>\z@\hskip\saveskip@\fi\lkerns@
 \else\vmodeerr@\nolinebreak\fi\fi}
\def\allowlinebreak{\RIfM@\mathmodeerr@\allowlinebreak\else
 \ifhmode\saveskip@\lastskip\unskip
 \allowbreak\ifdim\saveskip@>\z@\hskip\saveskip@\fi\lkerns@
 \else\vmodeerr@\allowlinebreak\fi\fi}
\def\linebreak{\RIfM@\mathmodeerr@\linebreak\else
 \ifhmode\unskip\unkern\break\lkerns@
 \else\vmodeerr@\linebreak\fi\fi}
\let\nkerns@\relax
\def\newline{\RIfM@\mathmodeerr@\newline\else
 \ifhmode\unskip\unkern\null\hfill\break\nkerns@
 \else\vmodeerr@\newline\fi\fi}%
\def\newbox@{\alloc@@4\box\chardef\insc@unt}
\def\newcount@{\alloc@@0\count\countdef\insc@unt}
\def\accentedsymbol#1#2{\expandafter\newbox@\csname\exstring@#1@box\endcsname
 \setbox\csname\exstring@#1@box\endcsname\hbox{$\m@th#2$}%
 \define#1{\copy\csname\exstring@#1@box\endcsname{}}}
\def\rightadd@#1\to#2{\toks@{\\#1}\toks@@\expandafter{#2}\xdef#2{\the\toks@@
 \the\toks@}\toks@{}\toks@@{}}
\def\fontlist@{\\\tenrm\\\sevenrm\\\fiverm\\\teni\\\seveni\\\fivei
 \\\tensy\\\sevensy\\\fivesy\\\tenex\\\tenbf\\\sevenbf\\\fivebf
 \\\tensl\\\tenit}
\def\font@#1=#2 {\rightadd@#1\to\fontlist@\font#1=#2 }
\def\ismember@#1#2{\global\let\Next@ F\let\next@= #2%
 {\def\\##1{\let\nextii@##1\ifx\nextii@\next@\global\let\Next@ T\fi}#1}%
 \test@false\ifx\Next@ T\test@true\fi\let\next@\relax}
\def\FNSS@#1{\let\FNSS@@#1\FN@\FNSS@@@}
\def\FNSS@@@{\ifx\next\space@\def\FNSS@@@@. {\FN@\FNSS@@@}\else
 \def\FNSS@@@@.{\FNSS@@}\fi\FNSS@@@@.}
\atdef@"{\unskip
 \DN@{\ifx\next`\DN@`{\FN@\nextii@}%
  \else\ifx\next\lq\DN@\lq{\FN@\nextii@}%
  \else\DN@####1{\FN@\nextiii@}\fi\fi
  \next@}%
 \DNii@{\ifx\next`\DN@`{\sldl@``}%
  \else\ifx\next\lq\DN@\lq{\sldl@``}%
  \else\DN@{\dlsl@`}\fi\fi\next@}%
 \def\nextiii@{\ifx\next'\DN@'{\srdr@''}%
  \else\ifx\next\rq\DN@\rq{\srdr@''}%
  \else\DN@{\drsr@'}\fi\fi\next@}%
 \FNSS@\next@}
\def\root{%
  \DN@{\ifx\next\uproot\let\next@\nextii@\else
   \ifx\next\leftroot\let\next@\nextiii@\else
   \let\next@\plainroot@\fi\fi\next@}%
  \DNii@\uproot##1{\uproot@##1\relax\FNSS@\nextiv@}%
  \def\nextiv@{\ifx\next\leftroot\let\next@\nextv@\else
   \let\next@\plainroot@\fi\next@}%
  \def\nextv@\leftroot##1{\leftroot@##1\relax\plainroot@}%
  \def\nextiii@\leftroot##1{\leftroot@##1\relax\FNSS@\nextvi@}%
  \def\nextvi@{\ifx\next\uproot\let\next@\nextvii@\else
   \let\next@\plainroot@\fi\next@}%
  \def\nextvii@\uproot##1{\uproot@##1\relax\plainroot@}%
  \bgroup\uproot@\z@\leftroot@\z@
 \FNSS@\next@}
\def\loop#1\repeat{\def\iterate{#1\relax\expandafter\iterate\fi}%
 \iterate\let\iterate\relax}
\def\gloop@#1\repeat{\gdef\iterate@{#1\relax\expandafter\iterate@\fi}%
 \iterate@\global\let\iterate@\relax}
\def\printoptions{\W@{Do you want S(yntax check),
  G(alleys) or P(ages)?^^JType S, G or P, follow by <return>: }\loop
 \read\m@ne to\ans@
 \edef\next@{\def\noexpand\Ans@{\ans@}}\uppercase\expandafter{\next@}%
 \ifx\Ans@\S@\test@true\syntax\else
 \ifx\Ans@\G@\test@true\galleys\else
 \ifx\Ans@\P@\test@true\else
 \test@false\fi\fi\fi
 \iftest@\else\W@{Type S, G or P, follow by <return>: }%
 \repeat}
\expandafter\let\csname A@;\endcsname;
\expandafter\let\csname A@:\endcsname:
\expandafter\let\csname A@?\endcsname?
\expandafter\let\csname A@!\endcsname!
\def\APdef#1{\def\next@{\expandafter\let\csname A@\string#1\endcsname#1}%
 \afterassignment\next@\def#1}
\let\fextra@\,
\def\tdots@{\unskip
 \DN@{$\m@th\mathinner{\ldotp\ldotp\ldotp}\,
   \ifx\next,\,$\else\ifx\next.\,$\else
   \ifx\next;\,$\else
   \expandafter\ifx\csname A@\string;\endcsname\next\fextra@$\else
   \ifx\next:\,$\else
   \expandafter\ifx\csname A@\string:\endcsname\next\fextra@$\else
   \ifx\next?\,$\else
   \expandafter\ifx\csname A@\string?\endcsname\next\fextra@$\else
   \ifx\next!\,$\else
   \expandafter\ifx\csname A@\string!\endcsname\next\fextra@$\else
   $ \fi\fi\fi\fi\fi\fi\fi\fi\fi\fi}%
 \ \FN@\next@}
\def\extrap@#1{%
 \ifx\next,\DN@{#1\,}\else
 \ifx\next;\DN@{#1\,}\else
 \expandafter\ifx\csname A@\string;\endcsname\next\DN@{#1\fextra@}\else
 \ifx\next.\DN@{#1\,}\else\extra@
 \ifextra@\DN@{#1\,}\else
 \let\next@#1\fi\fi\fi\fi\fi\next@}
\def\dotsc{\DN@{\ifx\next;\plainldots@\,\else
 \expandafter\ifx\csname A@\string;\endcsname\next\plainldots@\fextra@\else
 \ifx\next.\plainldots@\,\else\extra@\plainldots@
 \ifextra@\,\fi\fi\fi\fi}%
 \FN@\next@}
\def\keybin@{\keybin@true
 \ifx\next+\else\ifx\next=\else\ifx\next<\else\ifx\next>\else\ifx\next-\else
 \ifx\next*\else\ifx\next:\else
 \expandafter\ifx\csname A@\string;\endcsname\next\else
 \keybin@false\fi\fi\fi\fi\fi\fi\fi\fi}
\def\boldkey#1{\ifcat\noexpand#1A%
  \ifcmmibloaded@{\fam\cmmibfam#1}\else
   \Err@{First bold symbol font not loaded}\fi
 \else
 \let\next=#1%
 \ifx#1!\mathchar"5\bffam@21 \else
 \expandafter\ifx\csname A@\string!\endcsname\next\mathchar"5\bffam@21 \else
 \ifx#1(\mathchar"4\bffam@28 \else\ifx#1)\mathchar"5\bffam@29 \else
 \ifx#1+\mathchar"2\bffam@2B \else\ifx#1:\mathchar"3\bffam@3A \else
 \expandafter\ifx\csname A@\string:\endcsname\next\mathchar"3\bffam@3A \else
 \ifx#1;\mathchar"6\bffam@3B \else
 \expandafter\ifx\csname A@\string;\endcsname\next\mathchar"6\bffam@3B \else
 \ifx#1=\mathchar"3\bffam@3D \else
 \ifx#1?\mathchar"5\bffam@3F \else
 \expandafter\ifx\csname A@\string?\endcsname\next\mathchar"5\bffam@3F \else
 \ifx#1[\mathchar"4\bffam@5B \else
 \ifx#1]\mathchar"5\bffam@5D \else
 \ifx#1,\mathchari@63B \else
 \ifx#1-\mathcharii@200 \else
 \ifx#1.\mathchari@03A \else
 \ifx#1/\mathchari@03D \else
 \ifx#1<\mathchari@33C \else
 \ifx#1>\mathchari@33E \else
 \ifx#1*\mathcharii@203 \else
 \ifx#1|\mathcharii@06A \else
 \ifx#10\bold0\else\ifx#11\bold1\else\ifx#12\bold2\else\ifx#13\bold3\else
 \ifx#14\bold4\else\ifx#15\bold5\else\ifx#16\bold6\else\ifx#17\bold7\else
 \ifx#18\bold8\else\ifx#19\bold9\else
  \Err@{\noexpand\boldkey can't be used with #1}%
 \fi\fi\fi\fi\fi\fi\fi\fi\fi\fi\fi\fi\fi\fi\fi
 \fi\fi\fi\fi\fi\fi\fi\fi\fi\fi\fi\fi\fi\fi\fi\fi\fi\fi}
\def\arabic#1{#1}
\def\alph#1{\count@#1\relax\advance\count@96 \ifnum\count@>122
 \Err@{\noexpand\alph invalid for numbers > 26}\else\char\count@\fi}
\def\Alph#1{\count@#1\relax\advance\count@64 \ifnum\count@>90
 \Err@{\noexpand\Alph invalid for numbers > 26}\else\char\count@\fi}

\def\Roman#1{\uppercase\expandafter{\romannumeral#1}}
\def\fnsymbol#1{\count@#1\relax
 \count@@\count@
 \advance\count@\m@ne\divide\count@7
 \count@@@\count@\advance\count@@@\@ne
 \multiply\count@7 \advance\count@@-\count@
 \count@\count@@@
 {\loop
  \ifcase\count@@\or*\or\dag\or\ddag\or\P\or\S\or\text{$\|$}\or\#\fi
  \advance\count@\m@ne\ifnum\count@>\z@\repeat}}
\def\cardnine@#1{\ifcase#1\or one\or two\or three\or four\or five\or
 six\or seven\or eight\or nine\fi}
\let\alloc@\alloc@@
\newcount\ten@
\ten@10
\def\cardinal#1{\count@#1\relax
 \ifnum\count@>99 \number\count@
 \else
  \ifnum\count@=\z@ zero%
  \else
   \ifnum\count@<\ten@\cardnine@\count@
   \else
    \ifnum\count@<20
     \advance\count@-\ten@
     \ifcase\count@ ten\or eleven\or twelve\or thirteen\or fourteen\or
      fifteen\or sixteen\or seventeen\or eighteen\or nineteen\fi
    \else
     \count@@\count@\count@@@\count@@
     \divide\count@\ten@\multiply\count@\ten@
     \advance\count@@@-\count@\divide\count@\ten@
     \ifcase\count@\or\or twenty\or thirty\or forty\or fifty\or sixty\or
      seventy\or eighty\or ninety\fi
     \ifnum\count@@@=\z@\else-\cardnine@\count@@@\fi
    \fi
   \fi
  \fi
 \fi}
\def\ordnine@#1{\ifcase#1\or first\or second\or third\or fourth\or fifth\or
 sixth\or seventh\or eighth\or ninth\fi}
\newcount\count@@@@
\def\ordsuffix@{\count@@@@\count@
 \divide\count@\ten@
 \count@@@\count@\count@@\count@
 \divide\count@@\ten@\multiply\count@@\ten@
 \advance\count@@@-\count@@
 \ifnum\count@@@=\@ne th%
 \else
  \count@@@\count@@@@
  \count@@\count@@@@
  \divide\count@@\ten@\multiply\count@@\ten@
  \advance\count@@@-\count@@
  \ifcase\count@@@ th\or st\or nd\or rd\else th\fi
 \fi}
\def\nordinal#1{\count@#1\relax\number\count@\ordsuffix@}
\def\spordinal#1{\count@#1\relax\number\count@$^{\text{\ordsuffix@}}$}
\def\ordinal#1{\count@#1\relax
 \ifnum\count@>99 \number\count@\ordsuffix@
 \else
   \ifnum\count@=\z@ zeroth%
  \else
    \ifnum\count@<\ten@\ordnine@\count@
    \else
     \ifnum\count@<20 \advance\count@-\ten@
      \ifcase\count@ tenth\or eleventh\or twelfth\or thirteenth\or
       fourteenth\or fifteenth\or sixteenth\or seventeenth\or eighteenth\or
       nineteenth\fi
     \else
      \count@@\count@
      \divide\count@\ten@\multiply\count@\ten@
      \count@@@\count@@\advance\count@@@-\count@
      \divide\count@\ten@
      \ifcase\count@\or\or twent\or thirt\or fort\or fift\or sixt\or sevent\or
       eight\or ninet\fi
      \ifnum\count@@@=\z@ ieth\else y-\ordnine@\count@@@\fi
     \fi
    \fi
  \fi
 \fi}
\font@\tensmc=cmcsc10
\textonlyfont@\smc\tensmc
\newtoks\noexpandtoks@
\noexpandtoks@{\let\arabic\relax\let\alph\relax\let\Alph\relax
 \let\Roman\relax\let\fnsymbol\relax\let\rm\relax
 \let\it\relax\let\bf\relax\let\sl\relax\let\smc\relax
 \let\/\relax\let\null\relax}
\def\noexpands@{\the\noexpandtoks@}
\def\Nonexpanding#1{\global\noexpandtoks@
 \expandafter{\the\noexpandtoks@\let#1\relax}}
\def\prevanish@{\saveskip@\z@\ifhmode\saveskip@\lastskip\unskip\fi}
\def\postvanish@{\ifdim\saveskip@>\z@\hskip\saveskip@\fi\FN@\postvanish@@}
\def\postvanish@@{\DN@.{}%
 \ifx\next\space@\ifdim\saveskip@>\z@\DN@. {}\fi\fi\next@.}
\def\invisible#1{\prevanish@\ignorespaces#1\unskip\postvanish@}
\def\vanishlist@{\\\invisible}
\let\noindent@\noindent
\def\noindent{\par\noindent@\FN@\pretendspace@}
\def\pretendspace@{\ismember@\vanishlist@\next
 \iftest@\nobreak\hskip-\p@\hskip\p@\fi}

\newtoks\everypartoks@
\def\noindent@@{\par\everypartoks@\expandafter{\the\everypar}\everypar{}%
 \noindent@\everypar\expandafter{\the\everypartoks@}}
\def\page{\Err@{\noexpand\page has no meaning by itself}}
\let\page@C\pageno
\let\page@P\empty
\let\page@Q\empty
\def\page@S#1{#1\/}
\def\page@F{\rm}
\def\page@N{\arabic}   
\newif\ifindexing@
\def\indexfile{\ifindexing@\else
 \alloc@@7\write\chardef\sixt@@n\ndx@
 \immediate\openout\ndx@=\jobname.ndx
 \global\indexing@true\fi}
\global\advance\insc@unt\m@ne
\ch@ck0\insc@unt\count
\ch@ck1\insc@unt\dimen
\ch@ck2\insc@unt\skip
\ch@ck4\insc@unt\box
\allocationnumber\insc@unt
\global\chardef\margin@\allocationnumber
\dimen\margin@\maxdimen
\count\margin@\z@
\skip\margin@\z@
\newif\ifindexproofing@
\def\indexproofing{\indexproofing@true}
\def\noindexproofing{\indexproofing@false}
\def\unmacro@#1:#2->#3\unmacro@{\def\macpar@{#2}\def\macdef@{#3}}
\def\starparts@#1{\def\stari@{#1}\def\starii@{#1}\let\stariii@\empty
 \test@false
 \DN@##1*##2##3\next@{\ifx\starparts@##2\test@false\else\test@true\fi}%
 \next@#1*\starparts@\next@
 \iftest@\DN@{\starparts@@#1\starparts@@}\else\let\next@\relax\fi\next@}
\def\starparts@@#1*#2\starparts@@{\def\starii@{#1}\def\stariii@{*#2}}
\def\windex@{\ifindexing@
 \expandafter\unmacro@\meaning\stari@\unmacro@
 \edef\macdef@{\string"\macdef@\string"}%
 \edef\next@{\write\ndx@{\macdef@}}\next@
 \write\ndx@{{\number\pageno}{\page@N}{\page@P}{\page@Q}}%
 \fi
 \ifindexproofing@
  \ifx\stariii@\empty\else
   \expandafter\unmacro@\meaning\stariii@\unmacro@\fi
  \insert\margin@{\hbox{\rm\vrule\height9\p@\depth2\p@\width\z@\starii@
  \ifx\stariii@\empty\else\tt\macdef@\fi}}\fi}
\catcode`\"=\active
\def"{\FN@\quote@}
\def\quote@{\ifx\next"\expandafter\quote@@\else\expandafter\quote@@@\fi}
\def\quote@@@#1"{\starparts@{#1}\starii@\windex@}
\def\quote@@"#1"{\prevanish@\starparts@{#1}\windex@\FN@\quote@@@@}
\def\quote@@@@{\ifx\next"\DN@"{\postvanish@}\else
 \let\next@\postvanish@\fi\next@}
\rightadd@"\to\vanishlist@
\def\idefine#1{\DN@{#1}\DNii@{\noexpand#1}%
 \afterassignment\idefine@\def\nextiii@}
\def\idefine@{\ifindexing@
 \expandafter\let\next@\nextiii@
 \expandafter\unmacro@\meaning\nextiii@\unmacro@
 \immediate\write\ndx@{\noexpand\define\nextii@\macpar@{\macdef@}}\fi}
\def\iabbrev*#1#2{\ifindexing@\toks@{#2}%
 \immediate\write\ndx@{\noexpand\abbrev*\noexpand#1{\the\toks@}}\fi}
\newread\laxread@
\newwrite\laxwrite@
\let\fnpages@\empty
\def\Finit@#1#2\Finit@{\let\nextii@#1\def\nextiii@{#2}}
\catcode`\~=11
\def\getparts@ @#1~#2~#3~#4~#5~#6{\def\nextiv@{#1}%
 \def\nextiii@{#2~#3~#4~#5~}\count@#6\relax}
\newif\ifdocument@
\def\document{\ifdocument@\else\global\document@true
 \let\fontlist@\empty
 \immediate\openin\laxread@=\jobname.lax\relax
 {\endlinechar\m@ne\noexpands@\catcode`\@=11 \catcode`\~=11
  \loop\ifeof\laxread@\else
   \read\laxread@ to\next@
   \ifx\next@\empty
   \else
    \expandafter\Finit@\next@\Finit@
    \if\nextii@ F%
     \expandafter\rightadd@\nextiii@\to\fnpages@
    \else
     \expandafter\getparts@\next@
     \edef\next@{\gdef\csname\nextiv@ @L\endcsname{\nextiii@\number\count@}}%
     \next@
    \fi
   \fi
  \repeat}%
 \immediate\closein\laxread@
 \immediate\openout\laxwrite@=\jobname.lax\relax\fi}
\let\thelabel@\relax
\def\thelabels@{\thelabel@ ~\thelabel@@ ~\thelabel@@@ ~\thelabel@@@@ ~}
\def\label#1{\prevanish@
 \ifx\thelabel@\relax
  \Err@{There's nothing here to be labelled}%
 \else
  {\noexpands@
  \expandafter\ifx\csname#1@L\endcsname\relax
   \expandafter\xdef\csname#1@L\endcsname{\thelabels@0}%
   \immediate\write\laxwrite@{@#1~\thelabels@1}%
  \else
   \edef\next@{@~\csname#1@L\endcsname}%
    \expandafter\getparts@\next@
    \ifodd\count@
    \expandafter\xdef\csname#1@L\endcsname{\thelabels@0}%
    \immediate\write\laxwrite@{@#1~\thelabels@1}%
   \else
    \Err@{Label #1 already used}%
   \fi
  \fi
  }%
 \fi
 \postvanish@}
\rightadd@\label\to\vanishlist@
\def\thepages@{\page@N{\number\page@C}~%
 \page@S{\page@P\page@N{\number\page@C}\page@Q}~%
 \number\page@C ~\page@P\page@N{\number\page@C}\page@Q ~}
\def\pagelabel#1{\prevanish@
 \expandafter\ifx\csname#1@L\endcsname\relax
  {\noexpands@
  \expandafter\xdef\csname#1@L\endcsname{\thepages@2}}%
  \write\laxwrite@{@#1~\thepages@3}%
 \else
  {\noexpands@
  \edef\next@{@~\csname#1@L\endcsname}%
  \expandafter\getparts@\next@
  \ifodd\count@
   \ifnum\count@=\@ne
    \expandafter\xdef\csname#1@L\endcsname{\thelabels@2}%
   \fi
   \write\laxwrite@{@#1~\thepages@3}%
  \else
   \Err@{Label #1 already used}%
  \fi
  }%
 \fi
 \postvanish@}
\rightadd@\pagelabel\to\vanishlist@
\newif\ifreferr@
\referr@true
\def\RefErrors{\global\referr@true}
\def\RefWarnings{\global\referr@false}
\setbox\z@\hbox{\global\count@=`^^30}
\ifnum\count@=48 \let\versionthree@\relax\fi
\def\nolabel@#1#2#3{\expandafter\ifx\csname#2@L\endcsname\relax
 \ifreferr@\Err@{No \noexpand\label found for #2}\else
 \W@{Warning: No \noexpand\label found for #2.}%
 \ifx\versionthree@\relax\W@{l.\number\inputlineno\space ... \string#1{#2}}\fi
 \fi#3\else}
\def\csL@#1{{\noexpands@\xdef\Next@{\csname#1@L\endcsname}}}
\def\ref#1{\nolabel@\ref{#1}\relax
 \DNii@##1~##2\nextii@{##1}%
 \csL@{#1}\expandafter\nextii@\Next@\nextii@\fi}
\def\Ref#1{\nolabel@\Ref{#1}\relax
 \DNii@##1~##2~##3\nextii@{##2}%
 \csL@{#1}\expandafter\nextii@\Next@\nextii@\fi}
\def\nref#1{\nolabel@\nref{#1}\relax
 \DNii@##1~##2~##3~##4\nextii@{##3}%
 \csL@{#1}\expandafter\nextii@\Next@\nextii@\fi}
\def\pref#1{\nolabel@\pref{#1}\relax
 \DNii@##1~##2~##3~##4~##5\nextii@{##4}%
 \csL@{#1}\expandafter\nextii@\Next@\nextii@\fi}
\let\pref@\pref
\def\Evaluatenref#1{\nolabel@\Evaluatenref{#1}{\gdef\Nref{-10000 }}%
 \DNii@##1~##2~##3~##4\nextii@{\DNii@{##3}}%
 \csL@{#1}\expandafter\nextii@\Next@\nextii@
 \xdef\Nref{\nextii@}\fi}
\def\Evaluatepref#1{\nolabel@\Evaluatepref{#1}{\global\let\Pref\empty}%
 \DNii@##1~##2~##3~##4~##5\nextii@{\DNii@{##4}}%
 \csL@{#1}\expandafter\nextii@\Next@\nextii@
 \xdef\Pref{\nextii@}\fi}
\def\readlax#1{\immediate\openin\laxread@=#1.lax\relax
 \ifeof\laxread@\W@{}\W@{File #1.lax not found.}\W@{}\fi
 {\endlinechar\m@ne\noexpands@\catcode`\@=11 \catcode`\~=11
  \loop\ifeof\laxread@\else
   \read\laxread@ to\nextv@
   \ifx\nextv@\empty
   \else
    \expandafter\Finit@\nextv@\Finit@
    \ifx\nextii@ F%
    \else
     \expandafter\getparts@\nextv@
     \expandafter\ifx\csname\nextiv@ @L\endcsname\relax
      \edef\next@{\gdef\csname\nextiv@ @L\endcsname
       {\nextiii@\ifnum\count@=\@ne0\else2\fi}}%
      \next@
     \else
      \Err@{Label \nextiv@\space in #1.lax already used}%
     \fi
    \fi
   \fi
  \repeat}%
 \immediate\closein\laxread@}
\catcode`\~=\active

\def\FNSSP@{\FNSS@\pretendspace@}
\everydisplay{\csname displaymath \endcsname}
\expandafter\def\csname displaymath \endcsname#1$${#1$$\FNSSP@}
\def\locallabel@{\let\thelabel@\Thelabel@\let\thelabel@@\Thelabel@@
 \let\thelabel@@@\Thelabel@@@\let\thelabel@@@@\Thelabel@@@@}
\newcount\tag@C
\tag@C\z@
\let\tag@P\empty
\let\tag@Q\empty
\def\tag@S#1{{\rm(}{#1\/}{\rm)}}
\let\tag@N\arabic
\def\tag@F{\rm}
\def\maketag@{\FN@\maketag@@}
\def\maketag@@{\ifx\next\relax\DN@\relax{\FN@\maketag@@}\else
 \ifx\next"\let\next@\maketag@@@\else
 \let\next@\maketag@@@@\fi\fi\next@}
\def\xdefThelabel@#1{\xdef\Thelabel@{#1{\Thelabel@@@}}}
\def\xdefThelabel@@#1{\xdef\Thelabel@@{#1{\Thelabel@@@@}}}
\def\maketag@@@@#1\maketag@{\global\advance\tag@C\@ne
 {\noexpands@
  \xdef\Thelabel@@@{\number\tag@C}%
  \xdefThelabel@\tag@N
  \xdef\Thelabel@@@@{\ifmathtags@$\tag@P\Thelabel@\tag@Q$\else
   \tag@P\Thelabel@\tag@Q\fi}%
  \xdefThelabel@@\tag@S
  }%
 \locallabel@
 \hbox{\tag@F\thelabel@@}%
 #1}
\def\Qlabel@#1{{\noexpands@\xdef\Thelabel@@{#1}%
 \let\style\empty\xdef\Thelabel@@@@{#1}%
 \let\pre\empty\let\post\empty\xdef\Thelabel@{#1}%
 \let\numstyle\empty\xdef\Thelabel@@@{#1}}}
\def\maketag@@@"#1"#2\maketag@{%
 {\let\pre\tag@P\let\post\tag@Q\let\style\tag@S\let\numstyle\tag@N
  \hbox{\tag@F#1}%
  \noexpands@
  \Qlabel@{#1}%
  }%
 \locallabel@
 #2}
\def\align@{\inalign@true\inany@true
 \vspace@\allowdisplaybreak@\displaybreak@\intertext@
 \def\tag{\global\tag@true\ifnum\and@=\z@
  \DN@{&\omit\global\rwidth@\z@&\relax}\else
  \DN@{&\relax}\fi\next@}%
 \iftagsleft@\DN@{\csname align \endcsname}\else
  \DN@{\csname align \space\endcsname}\fi\next@}
\def\noset@{\def\Offset##1##2{\prevanish@\postvanish@}%
 \def\Reset##1##2{\prevanish@\postvanish@}}
\def\measure@#1\endalign{\global\lwidth@\z@\global\rwidth@\z@
 \global\maxlwidth@\z@\global\maxrwidth@\z@
 \global\and@\z@
 \setbox\z@\vbox
  {\noset@\everycr{\noalign{\global\tag@false\global\and@\z@}}\Let@
  \halign{\setboxz@h{$\m@th\displaystyle{\@lign##}$}%
   \global\lwidth@\wdz@
   \ifdim\lwidth@>\maxlwidth@\global\maxlwidth@\lwidth@\fi
   \global\advance\and@\@ne
   &\setboxz@h{$\m@th\displaystyle{{}\@lign##}$}\global\rwidth@\wdz@
   \ifdim\rwidth@>\maxrwidth@\global\maxrwidth@\rwidth@\fi
   \global\advance\and@\@ne
   &\Tag@\eat@{##}\crcr#1\crcr}}%
 \totwidth@\maxlwidth@\advance\totwidth@\maxrwidth@}
\def\prepost@{\global\let\tag@P@\tag@P\global\let\tag@Q@\tag@Q}
\def\reprepost@{\let\tag@P\tag@P@\let\tag@Q\tag@Q@}
\expandafter\def\csname align \space\endcsname#1\endalign
 {\measure@#1\endalign\global\and@\z@
 \ifingather@\everycr{\noalign{\global\and@\z@}}\else\displ@y@\fi
 \Let@\tabskip\centering@
 \halign to\displaywidth
  {\hfil\strut@\setboxz@h{$\m@th\displaystyle{\@lign##\prepost@}$}%
  \boxz@\global\advance\and@\@ne
  \tabskip\z@skip
  &\setboxz@h{$\m@th\displaystyle{{}\@lign##\prepost@}$}%
  \global\rwidth@\wdz@\boxz@\hfil\global\advance\and@\@ne
  \tabskip\centering@
  &\setboxz@h{\@lign\strut@\reprepost@\maketag@##\maketag@}%
  \dimen@\displaywidth\advance\dimen@-\totwidth@
  \divide\dimen@\tw@\advance\dimen@\maxrwidth@\advance\dimen@-\rwidth@
  \ifdim\dimen@<\tw@\wdz@\llap{\vtop{\normalbaselines\null\boxz@}}%
  \else\llap{\boxz@}\fi
  \tabskip\z@skip
  \crcr#1\crcr
  \black@\totwidth@}}
\expandafter\def\csname align \endcsname#1\endalign{\measure@#1\endalign
 \global\and@\z@
 \ifdim\totwidth@>\displaywidth\let\displaywidth@\totwidth@\else
  \let\displaywidth@\displaywidth\fi
 \ifingather@\everycr{\noalign{\global\and@\z@}}\else\displ@y@\fi
 \Let@\tabskip\centering@\halign to\displaywidth
  {\hfil\strut@\setboxz@h{$\m@th\displaystyle{\@lign##\prepost@}$}%
  \global\lwidth@\wdz@\global\lineht@\ht\z@
  \boxz@\global\advance\and@\@ne
  \tabskip\z@skip&\setboxz@h{$\m@th\displaystyle{{}\@lign##\prepost@}$}%
  \ifdim\ht\z@>\lineht@\global\lineht@\ht\z@\fi
  \boxz@\hfil\global\advance\and@\@ne
  \tabskip\centering@&\kern-\displaywidth@
  \setboxz@h{\@lign\strut@\reprepost@\maketag@##\maketag@}%
  \dimen@\displaywidth\advance\dimen@-\totwidth@
  \divide\dimen@\tw@\advance\dimen@\maxlwidth@\advance\dimen@-\lwidth@
  \ifdim\dimen@<\tw@\wdz@
   \rlap{\vbox{\normalbaselines\boxz@\vbox to\lineht@{}}}\else
   \rlap{\boxz@}\fi
  \tabskip\displaywidth@\crcr#1\crcr\black@\totwidth@}}
\def\attag@#1{\let\Maketag@\maketag@\let\TAG@\Tag@
 \let\Prepost@\prepost@\let\Reprepost@\reprepost@
 \let\Tag@\relax\let\maketag@\relax
 \let\prepost@\relax\let\reprepost@\relax
 \ifmeasuring@
  \def\llap@##1{\setboxz@h{##1}\hbox to\tw@\wdz@{}}%
  \def\rlap@##1{\setboxz@h{##1}\hbox to\tw@\wdz@{}}%
 \else\let\llap@\llap\let\rlap@\rlap\fi
 \toks@{\hfil\strut@
  $\m@th\displaystyle{\@lign\the\hashtoks@\prepost@}$%
  \tabskip\z@skip\global\advance\and@\@ne&
  $\m@th\displaystyle{{}\@lign\the\hashtoks@\prepost@}$\hfil
  \ifxat@\tabskip\centering@\fi\global\advance\and@\@ne}%
 \iftagsleft@
  \toks@@{\tabskip\centering@&\Tag@\kern-\displaywidth
   \rlap@{\@lign\reprepost@\maketag@\the\hashtoks@\maketag@}%
   \global\advance\and@\@ne\tabskip\displaywidth}\else
  \toks@@{\tabskip\centering@&\Tag@\llap@{\@lign\reprepost@\maketag@
   \the\hashtoks@\maketag@}\global\advance\and@\@ne\tabskip\z@skip}\fi
 \atcount@#1\relax\advance\atcount@\m@ne
 \loop\ifnum\atcount@>\z@
  \toks@\expandafter{\the\toks@&\hfil$\m@th\displaystyle{\@lign
  \the\hashtoks@\prepost@}$\global\advance\and@\@ne
  \tabskip\z@skip
  &$\m@th\displaystyle{{}\@lign\the\hashtoks@\prepost@}$\hfil\ifxat@
  \tabskip\centering@\fi\global\advance\and@\@ne}\advance\atcount@\m@ne
 \repeat
 \edef\preamble@{\the\toks@\the\toks@@}%
 \edef\preamble@@{\preamble@}%
 \let\maketag@\Maketag@\let\Tag@\TAG@
 \let\prepost@\Prepost@\let\reprepost@\Reprepost@}
\def\unlabel@{\def\label##1{\prevanish@\postvanish@}%
 \def\pagelabel##1{\prevanish@\postvanish@}}
\newcount\tag@CC
\expandafter\def\csname alignat \endcsname#1#2\endalignat
 {\inany@true\xat@false
 \def\tag{\global\tag@true
  \count@#1\relax\multiply\count@\tw@\advance\count@\m@ne
  \gdef\tag@{&}%
  \loop\ifnum\count@>\and@\xdef\tag@{&\omit\tag@}%
  \advance\count@\m@ne\repeat
  \tag@\relax}%
 \vspace@\allowdisplaybreak@\displaybreak@\intertext@
 \displ@y@\measuring@true\tag@CC\tag@C
 \setbox\savealignat@\hbox{\noset@\unlabel@$\m@th\displaystyle\Let@
  \attag@{#1}\vbox{\halign{\span\preamble@@\crcr#2\crcr}}$}%
 \measuring@false
 \Let@\attag@{#1}\tag@C\tag@CC
 \tabskip\centering@\halign to\displaywidth
  {\span\preamble@@\crcr#2\crcr\black@{\wd\savealignat@}}}
\expandafter\def\csname xalignat \endcsname#1#2\endxalignat
 {\inany@true\xat@true
 \def\tag{\global\tag@true
  \count@#1\relax\multiply\count@\tw@\advance\count@\m@ne
  \gdef\tag@{&}%
  \loop\ifnum\count@>\and@\xdef\tag@{&\omit\tag@}%
  \advance\count@\m@ne\repeat
  \tag@\relax}%
 \vspace@\allowdisplaybreak@\displaybreak@\intertext@
 \displ@y@\measuring@true\tag@CC\tag@C
 \setbox\savealignat@\hbox{\noset@\unlabel@$\m@th\displaystyle\Let@
  \attag@{#1}\vbox{\halign{\span\preamble@@\crcr#2\crcr}}$}%
 \measuring@false\Let@\attag@{#1}\tag@C\tag@CC
 \tabskip\centering@\halign to\displaywidth
 {\span\preamble@@\crcr#2\crcr\black@{\wd\savealignat@}}}
\def\gather{\RIfMIfI@\DN@{\onlydmatherr@\gather}\else
 \ingather@true\inany@true\def\tag{&\relax}%
 \vspace@\allowdisplaybreak@\displaybreak@\intertext@
 \displ@y\Let@
 \iftagsleft@\DN@{\csname gather \endcsname}\else
  \DN@{\csname gather \space\endcsname}\fi\fi
 \else\DN@{\onlydmatherr@\gather}\fi\next@}
\def\exstring@{\expandafter\eat@\string}
\def\newcounter#1{\define#1{}%
 \edef\next@{\def\noexpand#1{\futurelet\noexpand\next
  \csname\exstring@#1@Z\endcsname}}\next@
 \edef\next@{\def\csname\exstring@#1@Z\endcsname
  {\global\advance\csname\exstring@#1@C\endcsname\@ne
  {\csname\exstring@#1@F\endcsname\csname\exstring@#1@S\endcsname
   {\csname\exstring@#1@P\endcsname\csname\exstring@#1@N\endcsname
   {\noexpand\number\csname\exstring@#1@C\endcsname}%
   \csname\exstring@#1@Q\endcsname}}%
  \noexpand\ifx\noexpand\next\noexpand\label
   \def\noexpand\next@\noexpand\label########1{{\noexpand\noexpands@
    \xdef\noexpand\Thelabel@{\csname\exstring@#1@N\endcsname
     {\noexpand\number\csname\exstring@#1@C\endcsname}}%
    \xdef\noexpand\Thelabel@@@{\noexpand\number
     \csname\exstring@#1@C\endcsname}%
    \xdef\noexpand\Thelabel@@{\csname\exstring@#1@S\endcsname
     {\csname\exstring@#1@P\endcsname
     \csname\exstring@#1@N\endcsname
     {\noexpand\number\csname\exstring@#1@C\endcsname}%
     \csname\exstring@#1@Q\endcsname}}%
    \xdef\noexpand\Thelabel@@@@{\csname\exstring@#1@P\endcsname
     \csname\exstring@#1@N\endcsname
     {\noexpand\number\csname\exstring@#1@C\endcsname}%
     \csname\exstring@#1@Q\endcsname}}%
    {\noexpand\locallabel@\noexpand\label{########1}}}%
   \noexpand\else\let\noexpand\next@\relax\noexpand\fi\noexpand\next@}}\next@
 \expandafter\newcount@\csname\exstring@#1@C\endcsname
 \expandafter\let\csname\exstring@#1@N\endcsname\arabic
 \expandafter\def\csname\exstring@#1@S\endcsname##1{##1\/}%
 \expandafter\let\csname\exstring@#1@P\endcsname\empty
 \expandafter\let\csname\exstring@#1@Q\endcsname\empty
 \expandafter\def\csname\exstring@#1@F\endcsname{\rm}%
 }
\def\HASH@#1#2{\ifnum#2=\z@\else
 \edef\next@{\toks@{\the\toks@\the\hashtoks@#2}%
 \toks@@{\the\toks@@{\the\hashtoks@#2}}}\next@\expandafter\HASH@\fi}
\def\HASH@@{\toks@{}\toks@@{}\expandafter\HASH@\macpar@00}
\def\usecounter#1#2{\expandafter\ifx\csname\exstring@#1@Z\endcsname
 \relax\Err@{\noexpand#1not created with \string\newcounter}\fi
 \expandafter\let\csname\exstring@#1@@Z\endcsname\relax
 \expandafter\let\csname\exstring@#1@@Z@\endcsname\relax
 \expandafter\let\csname\exstring@#1@@Z@@\endcsname\relax
 \edef\next@{\def\noexpand#2{\futurelet\noexpand\next
  \csname\exstring@#1@@Z\endcsname}}\next@
 \edef\next@{\def\csname\exstring@#1@@Z\endcsname{\noexpand\ifx
  \noexpand\next\noexpand\label\def\noexpand\next@\noexpand\label
   ########1{\csname\exstring@#1@@Z@\endcsname
   {\noexpand#1\noexpand\label{########1}}}%
   \noexpand\else\noexpand\ifx\noexpand\next
   \noexpand"\def\noexpand\next@\noexpand"########1\noexpand"%
   {\csname\exstring@#1@@Z@\endcsname{{\expandafter\noexpand
   \csname\exstring@#1@F\endcsname
   \let\noexpand\pre\expandafter\noexpand\csname\exstring@#1@P\endcsname
   \let\noexpand\post\expandafter\noexpand\csname\exstring@#1@Q\endcsname
   \let\noexpand\style\expandafter\noexpand\csname\exstring@#1@S\endcsname
   \let\noexpand\numstyle\expandafter\noexpand\csname\exstring@#1@N\endcsname
   ########1}}}\noexpand\else
   \def\noexpand\next@{\csname\exstring@#1@@Z@\endcsname{\noexpand#1}}%
   \noexpand\fi\noexpand\fi\noexpand\next@}}\next@
 \def\next@{\expandafter\expandafter\expandafter\unmacro@\expandafter
  \meaning\csname\exstring@#1@@Z@@\endcsname\unmacro@
  \HASH@@
  \edef\next@{\def\csname\exstring@#1@@Z@\endcsname\the\toks@{%
   \expandafter\noexpand\csname\exstring@#1@@Z@@\endcsname\the\toks@@
   \noexpand\FNSSP@}}\next@}%
 \afterassignment\next@
 \expandafter\def\csname\exstring@#1@@Z@@\endcsname}
\def\listbi@{\penalty50 \medskip}
\def\listbii@{\penalty100 \smallskip}
\let\listbiii@\relax
\let\listbiv@\relax
\let\listbv@\relax
\def\listmi@{\advance\leftskip30\p@\relax}
\let\listmii@\listmi@
\let\listmiii@\listmi@
\let\listmiv@\listmi@
\let\listmv@\listmi@
\def\itemi@#1{\noindent@@\llap{#1\hskip5\p@}}
\let\itemii@\itemi@
\let\itemiii@\itemi@
\let\itemiv@\itemi@
\let\itemv@\itemi@
\def\liste@{\penalty-50 \medskip}
\def\listei@{\penalty-100 \smallskip}
\let\listeii@\relax
\let\listeiii@\relax
\let\listeiv@\relax
\expandafter\newcount\csname list@C1\endcsname
\csname list@C1\endcsname\z@
\expandafter\newcount\csname list@C2\endcsname
\csname list@C2\endcsname\z@
\expandafter\newcount\csname list@C3\endcsname
\csname list@C3\endcsname\z@
\expandafter\newcount\csname list@C4\endcsname
\csname list@C4\endcsname\z@
\expandafter\newcount\csname list@C5\endcsname
\csname list@C5\endcsname\z@
\expandafter\let\csname list@P1\endcsname\empty
\expandafter\let\csname list@P2\endcsname\empty
\expandafter\let\csname list@P3\endcsname\empty
\expandafter\let\csname list@P4\endcsname\empty
\expandafter\let\csname list@P5\endcsname\empty
\expandafter\let\csname list@Q1\endcsname\empty
\expandafter\let\csname list@Q2\endcsname\empty
\expandafter\let\csname list@Q3\endcsname\empty
\expandafter\let\csname list@Q4\endcsname\empty
\expandafter\let\csname list@Q5\endcsname\empty
\expandafter\def\csname list@S1\endcsname#1{{\rm(}{#1\/}{\rm)}}
\expandafter\def\csname list@S2\endcsname#1{{\rm(}{#1\/}{\rm)}}
\expandafter\def\csname list@S3\endcsname#1{{\rm(}{#1\/}{\rm)}}
\expandafter\def\csname list@S4\endcsname#1{{\rm(}{#1\/}{\rm)}}
\expandafter\def\csname list@S5\endcsname#1{{\rm(}{#1\/}{\rm)}}
\expandafter\let\csname list@N1\endcsname\arabic
\expandafter\let\csname list@N2\endcsname\arabic
\expandafter\let\csname list@N3\endcsname\arabic
\expandafter\let\csname list@N4\endcsname\arabic
\expandafter\let\csname list@N5\endcsname\arabic
\expandafter\def\csname list@F1\endcsname{\rm}
\expandafter\def\csname list@F2\endcsname{\rm}
\expandafter\def\csname list@F3\endcsname{\rm}
\expandafter\def\csname list@F4\endcsname{\rm}
\expandafter\def\csname list@F5\endcsname{\rm}
\newcount\listlevel@
\listlevel@\z@
\def\list@@C{\csname list@C\number\listlevel@\endcsname}
\def\list@@P{\csname list@P\number\listlevel@\endcsname}
\def\list@@Q{\csname list@Q\number\listlevel@\endcsname}
\def\list@@S{\csname list@S\number\listlevel@\endcsname}
\def\list@@N{\csname list@N\number\listlevel@\endcsname}
\def\list@@F{\csname list@F\number\listlevel@\endcsname}
\newif\iffirstitemi@
\newif\iffirstitemii@
\newif\iffirstitemiii@
\newif\iffirstitemiv@
\newif\iffirstitemv@
\def\Firstitem@true{\csname firstitem\romannumeral\listlevel@
 @true\endcsname}
\def\Firstitem@false{\csname firstitem\romannumeral\listlevel@
 @false\endcsname}
\def\Listm@{\csname listm\romannumeral\listlevel@ @\endcsname}
\def\Item@{\csname item\romannumeral\listlevel@ @\endcsname}
\def\Liste@{\csname liste\romannumeral\listlevel@ @\endcsname}
\newif\iflistcontinue@
\def\keepitem{\listcontinue@true}
\newcount\list@C@
\def\list{%
 \iflistcontinue@\csname list@C1\endcsname\csname list@C@\endcsname\fi
 \global\csname list@C2\endcsname\z@
 \global\csname list@C3\endcsname\z@
 \global\csname list@C4\endcsname\z@
 \global\csname list@C5\endcsname\z@
 \begingroup
 \firstitemi@true
 \listlevel@\@ne
 \def\item{\FN@\item@}%
 \FN@\list@}
\Invalid@\runinitem
\def\list@{\ifx\next\par
 \DN@\par{\FN@\list@}\else
 \ifx\next\runinitem
  \DN@\runinitem{\FN@\runinitem@}\else
  \DN@{\par\dimen@\parskip\parskip\dimen@}\fi\fi\next@}
\newif\ifoutlevel@
\newif\ifrunin@
\def\item@{%
 \ifoutlevel@\Liste@\outlevel@false\fi
 \ifrunin@\runin@false\par
  \dimen@\parskip\parskip\dimen@
  \Listm@\fi
 \iffirstitemi@\listbi@\listmi@\firstitemi@false\else\par\fi
 \iffirstitemii@\listbii@\listmii@\firstitemii@false\else\par\fi
 \iffirstitemiii@\listbiii@\listmiii@\firstitemiii@false\else\par\fi
 \iffirstitemiv@\listbiv@\listmiv@\firstitemiv@false\else\par\fi
 \iffirstitemv@\listbv@\listmv@\firstitemv@false\else\par\fi
 \DN@"##1"{{\let\pre\list@@P\let\post\list@@Q
  \let\style\list@@S\let\numstyle\list@@N
  \vskip-\parskip
  \Item@{\list@@F##1}%
  \noexpands@
  \Qlabel@{##1}}%
  \locallabel@
  \FNSSP@}%
 \DNii@{\global\advance\list@@C\@ne
  {\noexpands@
   \xdef\Thelabel@@@{\number\list@@C}%
   \xdefThelabel@\list@@N
   \xdef\Thelabel@@@@{\list@@P\Thelabel@\list@@Q}%
   \xdefThelabel@@\list@@S
  }%
  \locallabel@
  \vskip-\parskip
  \Item@{\list@@F\thelabel@@}%
  \FN@\pretendspace@}%
 \ifx\next"\expandafter\next@\else\expandafter\nextii@\fi}
\def\runinitem@{%
  \runin@true
  \Firstitem@false
  \DN@"##1"{{\let\pre\list@@P\let\post\list@@Q
   \let\style\list@@S\let\numstyle\list@@N
   \unskip\space{\list@@F##1} %
   \noexpands@
   \Qlabel@{##1}}%
   \locallabel@
   \ignorespaces}%
  \DNii@{\global\advance\list@@C\@ne
   {\noexpands@
    \xdef\Thelabel@@@{\number\list@@C}%
    \xdefThelabel@\list@@N
    \xdef\Thelabel@@@@{\list@@P\Thelabel@\list@@Q}%
    \xdefThelabel@@\list@@S
   }%
   \locallabel@
   \unskip\space{\list@@F\thelabel@@} }%
  \ifx\next"\expandafter\next@\else\expandafter\nextii@\fi}
\def\inlevel{\ifnum\listlevel@=5
 \DN@{\Err@{Already 5 levels down}}\else
 \DN@{\begingroup\advance\listlevel@\@ne
 \Firstitem@true\FN@\inlevel@}\fi\next@}
\def\inlevel@{\ifx\next\par
 \DN@\par{\FN@\inlevel@}\else
 \ifx\next\runinitem
  \DN@\runinitem{\FN@\runinitem@}\else
  \let\next@\relax\fi\fi\next@}
\def\outlevel{\ifnum\listlevel@=\@ne
 \Err@{At top level}\else
 \par\global\list@@C\z@\endgroup\outlevel@true\fi}
\def\endlist{%
 \expandafter\global\csname list@C@\endcsname\csname list@C1\endcsname
 \par
 \global\toks\@ne{}\count@\listlevel@
 {\loop
  \ifnum\count@>\z@\global\toks\@ne\expandafter{\the\toks\@ne\endgroup}%
  \advance\count@\m@ne
  \repeat}%
 \the\toks\@ne
 \liste@
 \listcontinue@false\global\csname list@C1\endcsname\z@
 \vskip-\parskip
 \noindent@@
 \FN@\pretendspace@}
\newif\iffirstdescribe@
\def\describe{\par
 \begingroup\firstdescribe@true
 \def\item##1{%
  \iffirstdescribe@\penalty50 \medskip\vskip-\parskip
  \firstdescribe@false\else\par\fi
  \noindent@@\hangindent2pc\hangafter\@ne
  {\bf##1}\hskip.5em}}

\Invalid@\pullin
\Invalid@\pullinmore
\newif\iffirstpull@
\def\margins{\par\begingroup\firstpull@true
 \def\pullin##1##2{\par
  \iffirstpull@\firstpull@false\else\endgroup\fi
  \begingroup\DN@{##1}%
  \ifx\next@\empty\leftskip\z@\else\ifx\next@\space\leftskip\z@
  \else\leftskip##1\fi\fi
  \DN@{##2}\ifx\next@\empty\rightskip\z@\else\ifx\next@\space
  \rightskip\z@\else\rightskip##2\fi\fi\ignorespaces}%
 \def\pullinmore##1##2{\par
  \xdef\Next@{\leftskip\the\leftskip\relax\rightskip\the\rightskip\relax}%
  \iffirstpull@\firstpull@false\else\endgroup\fi
  \begingroup\Next@
  \DN@{##1}%
  \ifx\next@\empty\else\ifx\next@\space\else\advance\leftskip##1\fi\fi
  \DN@{##2}\ifx\next@\empty\else\ifx\next@\space\else
  \advance\rightskip##2\fi\fi\ignorespaces}}

\newif\ifnopunct@
\newif\ifnospace@
\newif\ifoverlong@
\let\nofrillslist@\empty
\let\overlonglist@\empty
\def\nopunct{\nopunct@true\FN@\nopunct@}
\def\nospace{\nospace@true\FN@\nospace@}
\def\overlong{\overlong@true\FN@\overlong@}
\def\nopunct@{\ifx\next\nospace
 \DN@\nospace{\nospace@true\FN@\nopnos@}\else\ifx\next\overlong
 \DN@\overlong{\overlong@true\FN@\nopol@}\else
 \let\next@\nopunct@@\fi\fi\next@}
\def\nopunct@@#1{\ismember@\nofrillslist@#1%
 \iftest@\let\next@#1\else
 \DN@{\nopunct@false\Err@{\noexpand\nopunct can't be used with
 \string#1}#1}\fi\next@}
\def\nospace@{\ifx\next\nopunct
 \DN@\nopunct{\nopunct@true\FN@\nopnos@}\else\ifx\next\overlong
 \DN@\overlong{\overlong@true\FN@\nosol@}\else
 \let\next@\nospace@@\fi\fi\next@}
\def\nospace@@#1{\ismember@\nofrillslist@#1%
 \iftest@\let\next@#1\else
 \DN@{\nospace@false\Err@{\noexpand\nospace can't be used with
 \string#1}#1}\fi\next@}
\def\overlong@{\ifx\next\nopunct
 \DN@\nopunct{\nopunct@true\FN@\nopol@}\else\ifx\next\nospace
 \DN@\nospace{\nospace@true\FN@\nosol@}\else
 \let\next@\overlong@@\fi\fi\next@}
\def\overlong@@#1{\ismember@\overlonglist@#1%
 \iftest@\let\next@#1\else
 \DN@{\overlong@false\Err@{\noexpand\overlong can't be used with
 \string#1}#1}\fi\next@}
\def\nopnos@{\ifx\next\overlong
 \DN@\overlong{\overlong@true\nopnosol@}\else
 \let\next@\nopnos@@\fi\next@}
\def\nopol@{\ifx\next\nospace
 \DN@\nospace{\nospace@true\nopnosol@}\else
 \let\next@\nopol@@\fi\next@}
\def\nosol@{\ifx\next\nopunct
 \DN@\nopunct{\nopunct@true\nopnosol@}\else
 \let\next@\nosol@@\fi\next@}
\def\nopnos@@#1{\ismember@\nofrillslist@#1%
 \iftest@\let\next@#1\else
 \DN@{\nopunct@false\nospace@false
  \Err@{\noexpand\nopunct\noexpand\nospace
   can't be used with \string#1}#1}\fi\next@}
\def\testii@#1{\ismember@\nofrillslist@#1%
 \iftest@\let\nextiii@ T\else\let\nextiii@ F\fi
 \ismember@\overlonglist@#1%
 \iftest@\let\nextiv@ T\else\let\nextiv@ F\fi
 \test@false\if\nextiii@ T\if\nextiv@ T\test@true\fi\fi}
\def\nopol@@#1{\testii@{#1}%
 \iftest@\let\next@#1%
 \else\DN@{\if\nextiii@ T\else\nopunct@false\fi
  \if\nextiv@ T\else\overlong@false\fi
  \Err@{\if\nextiii@ T\else\noexpand\nopunct\fi
  \if\nextiv@ T\else\noexpand\overlong\fi can't be used
  with \string#1}#1}\fi\next@}
\def\nosol@@#1{\testii@{#1}%
 \iftest@\let\next@#1%
 \else\DN@{\if\nextiii@ T\else\nospace@false\fi
  \if\nextiv@ T\else\overlong@false\fi
  \Err@{\if\nextiii@ T\else\noexpand\nospace\fi
  \if\nextiv@ T\else\noexpand\overlong\fi can't be used
  with \string#1}#1}\fi\next@}
\def\nopnosol@#1{\testii@{#1}%
 \iftest@\let\next@#1%
 \else\DN@{\if\nextiii@ T\else\nopunct@false\nospace@false\fi
  \if\nextiv@ T\else\overlong@false\fi
  \Err@{\if\nextiii@ T\else\noexpand\nopunct\noexpand\nospace\fi
  \if\nextiv@ T\else\noexpand\overlong\fi can't be used
  with \string#1}#1}\fi\next@}
\def\punct@#1{\ifnopunct@\else#1\fi}
\def\addspace@#1{\ifnospace@\else#1\fi}
\def\hss@{\ifoverlong@\z@ plus\@m\p@ minus\@m\p@
 \else \z@ plus\@m\p@\fi}
\rightadd@\demo\to\nofrillslist@
\newif\ifclaim@
\def\exxx@{\expandafter\expandafter\expandafter\eat@\expandafter\string}
\let\colon@:
\def\demo#1{\ifclaim@
 \Err@{Previous \expandafter\noexpand\claimtype@ has
  no matching \string\end\exxx@\claimtype@}%
 \let\next@\relax
 \else
  \par
  \ifdim\lastskip<\smallskipamount\removelastskip\smallskip\fi
  \begingroup
  \noindent@@{\smc\ignorespaces#1\unskip
   \punct@{\null\colon@}\addspace@\enspace}%
  \nopunct@false\nospace@false
  \rm
  \DN@{\FNSSP@}%
 \fi
 \next@}
\def\enddemo{\par\endgroup\nopunct@false\nospace@false\smallskip}
\rightadd@\claim\to\nofrillslist@
\def\claim@F{\smc}
\def\claim@@@F{\csname\exxx@\claimtype@ @F\endcsname}
\def\claimformat@#1#2#3{%
 \medbreak\noindent@@{\smc#1 {\claim@@@F#2} #3%
 \punct@{\null.}\addspace@\enspace}\sl}
\def\claimformat@@#1#2{\claimformat@{\ignorespaces#1\unskip}%
 {\ifx\thelabel@@\empty\unskip\else\thelabel@@\fi}%
 {\ignorespaces#2\unskip}%
 \let\Claimformat@@\claimformat@@\FNSSP@}
\let\Claimformat@@\claimformat@@
\def\claim@@@P{\csname\exxx@\claimtype@ @P\endcsname}
\def\claim@@@Q{\csname\exxx@\claimtype@ @Q\endcsname}
\def\claim@@@S{\csname\exxx@\claimtype@ @S\endcsname}
\def\claim@@@N{\csname\exxx@\claimtype@ @N\endcsname}
\def\claim@@@C{\csname claim@C\claimclass@\endcsname}
\newcount\claim@C
\claim@C\z@
\let\claim@P\empty
\let\claim@Q\empty
\def\claim@S#1{#1\/}
\let\claim@N\arabic
\def\claim{\claim@true\let\claimclass@\empty
 \def\claimtype@{\claim}\FN@\claim@}
\def\claim@{%
 \ifx\next\c
  \let\next@\claim@c
 \else
  \ifx\next"%
   \let\next@\claim@q
  \else
   \begingroup\global\advance\claim@C\@ne
   {\noexpands@
    \xdef\Thelabel@@@{\number\claim@C}%
    \xdefThelabel@\claim@N
    \xdef\Thelabel@@@@{\claim@P\Thelabel@\claim@Q}%
    \xdefThelabel@@\claim@S
   }%
   \locallabel@
   \let\next@\Claimformat@@
  \fi
 \fi
 \next@}
\def\claim@c\c#1{\claim@true\begingroup
 \expandafter
 \ifx\csname claim@C#1\endcsname\relax
  \expandafter\newcount@\csname claim@C#1\endcsname
  \global\csname claim@C#1\endcsname\@ne
 \else
  \global\advance\csname claim@C#1\endcsname\@ne
 \fi
 \def\claimclass@{#1}%
 {\noexpands@
  \xdef\Thelabel@@@{\number\claim@@@C}%
  \xdefThelabel@\claim@@@N
  \xdef\Thelabel@@@@{\claim@@@P\Thelabel@\claim@@@Q}%
  \xdefThelabel@@\claim@@@S
 }%
 \locallabel@
 \FNSS@\claim@c@}
\def\claim@q"#1"{\begingroup
 {\let\pre\claim@@@P\let\post\claim@@@Q
  \let\style\claim@@@S\let\numstyle\claim@@@N
  \noexpands@
  \Qlabel@{#1}}%
 \locallabel@
 \FNSS@\claim@q@}
\def\claim@c@{\ifx\next"%
 \global\advance\claim@@@C\m@ne\let\next@\claim@cq
 \else\let\next@\Claimformat@@\fi\next@}
\def\claim@cq"#1"{{\let\pre\claim@@@P\let\post\claim@@@Q
 \let\style\claim@@@S\let\numstyle\claim@@@N
 \noexpands@
 \Qlabel@{#1}}%
 \locallabel@
 \FNSS@\Claimformat@@}
\def\claim@q@{\ifx\next\c\expandafter\claim@qc
 \else\expandafter\Claimformat@@\fi}
\def\claim@qc\c#1{\expandafter\ifx\csname claim@C#1\endcsname\relax
 \expandafter\newcount@\csname claim@C#1\endcsname
 \global\csname claim@C#1\endcsname\z@\fi
 \FNSS@\Claimformat@@}
\def\endclaim{\endgroup\claim@false\nopunct@false\nospace@false
 \let\Claimformat@@\claimformat@@\medbreak}
\Invalid@\claimclause
\def\newclaim{\FN@\newclaim@}
\def\newclaim@{\ifx\next\claimclause
 \DN@\claimclause##1{\newclaim@@{##1}}\else
 \DN@{\newclaim@@\relax}\fi\next@}
\def\claimlist@{\\\claim}
\newtoks\claim@i
\newtoks\claim@v
\let\noclaimclause@=F
\def\newclaim@@#1#2#3\c#4#5{\define#2{}%
 \rightadd@#2\to\claimlist@\rightadd@#2\to\nofrillslist@%
 \expandafter\def\csname\exstring@#2@P\endcsname{\claim@P}%
 \expandafter\def\csname\exstring@#2@Q\endcsname{\claim@Q}%
 \expandafter\def\csname\exstring@#2@S\endcsname{\claim@S}%
 \expandafter\def\csname\exstring@#2@N\endcsname{\claim@N}%
 \expandafter\def\csname\exstring@#2@F\endcsname{\claim@F}%
 \expandafter\def\csname end\exstring@#2\endcsname{\endclaim}%
 \expandafter\ifx\csname claim@C#4\endcsname\relax
  \expandafter\newcount@\csname claim@C#4\endcsname
  \global\csname claim@C#4\endcsname\z@\fi
 \edef\next@{\let\csname\exstring@#2@C\endcsname
   \csname claim@C#4\endcsname}\next@
 \def#2{\ifx\noclaimclause@ T\else#1\fi
  \global\claim@i{#1}\gdef\claim@iv{#4}\global\claim@v{#5}%
  \def\claimtype@{#2}\def\Claimformat@@{\claimformat@@{#5}}\claim@c\c{#4}}}
\def\shortenclaim#1#2{\define#2{}%
 \ismember@\claimlist@#1%
 \iftest@
  \rightadd@#2\to\nofrillslist@%
  \expandafter\def\csname\exstring@#2@P\endcsname
   {\csname\exstring@#1@P\endcsname}%
  \expandafter\def\csname\exstring@#2@Q\endcsname
   {\csname\exstring@#1@Q\endcsname}%
  \expandafter\def\csname\exstring@#2@S\endcsname
   {\csname\exstring@#1@S\endcsname}%
  \expandafter\def\csname\exstring@#2@N\endcsname
   {\csname\exstring@#1@N\endcsname}%
  \expandafter\def\csname\exstring@#2@F\endcsname
   {\csname\exstring@#1@F\endcsname}%
  \expandafter\def\csname end\exstring@#2\endcsname{\endclaim}%
  \edef\next@{\let\csname\exstring@#2@C\endcsname
    \csname claim\exstring@#1C\endcsname}\next@
  \setbox\z@\vbox{\let\noclaimclause@ T#1""\relax\endgroup}%
  \edef#2{\the\claim@i
   \def\noexpand\claimtype@{\noexpand#2}%
   \def\noexpand\Claimformat@@{\noexpand\claimformat@@{\the\claim@v}\relax}%
   \noexpand\claim@c\noexpand\c{\claim@iv}}%
 \else
  \Err@{\noexpand#1not yet created by \string\newclaim}%
 \fi}
\def\classtest@#1{\DN@{#1}\ifx\next@\claimclass@
 \test@true\else\test@false\fi}
\def\typetest@#1{\DN@{#1}\ifx\next@\claimtype@\test@true\else
  \test@false\fi}
\newif\iftoc@
\def\tocfile{\iftoc@\else\alloc@@7\write\chardef\sixt@@n\toc@
 \immediate\openout\toc@=\jobname.toc
 \alloc@@7\write\chardef\sixt@@n\tic@
 \immediate\openout\tic@=\jobname.tic
 \global\toc@true\fi}
\rightadd@\hl\to\nofrillslist@
\rightadd@\HL\to\overlonglist@
\def\HL@@C{\csname HL@C\HLlevel@\endcsname}
\def\HL@@P{\csname HL@P\HLlevel@\endcsname}
\def\HL@@Q{\csname HL@Q\HLlevel@\endcsname}
\def\HL@@S{\csname HL@S\HLlevel@\endcsname}
\def\HL@@N{\csname HL@N\HLlevel@\endcsname}
\def\HL@@F{\csname HL@F\HLlevel@\endcsname}
\def\HL@@@C{\csname\exxx@\HLtype@ @C\endcsname}
\def\HL@@@P{\csname\exxx@\HLtype@ @P\endcsname}
\def\HL@@@Q{\csname\exxx@\HLtype@ @Q\endcsname}
\def\HL@@@S{\csname\exxx@\HLtype@ @S\endcsname}
\def\HL@@@N{\csname\exxx@\HLtype@ @N\endcsname}
\def\HL#1{\expandafter
 \ifx\csname HL@C#1\endcsname\relax
  \DN@{\Err@{\string\HL#1 not defined in this style}}%
 \else
  \DN@{\gdef\HLlevel@{#1}\def\HLname@{\HL{#1}}\let\HLtype@\relax\FNSS@\HL@}%
 \fi
 \next@}%
\newif\ifquoted@
\let\aftertoc@\relax
\def\HL@{%
 \DN@"##1"##2\endHL{\def\entry@{##2}\quoted@true
  {\noexpands@
  \ifx\HLtype@\relax
   \let\pre\HL@@P\let\post\HL@@Q\let\style\HL@@S\let\numstyle\HL@@N
  \else
   \let\pre\HL@@@P\let\post\HL@@@Q\let\style\HL@@@S\let\numstyle\HL@@@N
  \fi
  \Qlabel@{##1}\let\style\relax\xdef\Qlabel@@@@{##1}%
  \xdef\Thepref@{\Thelabel@@@@}}%
  \csname HL@\HLlevel@\endcsname##2\endHL
  \let\pref\Thepref@
  \csname HL@I\HLlevel@\endcsname
  \csname HL@J\HLlevel@\endcsname
  \let\pref\pref@
  \HLtoc@	
  \aftertoc@
  \let\aftertoc@\relax\overlong@false}%
 \DNii@##1\endHL{\def\entry@{##1}\quoted@false
  {\noexpands@
  \ifx\HLtype@\relax
   \global\advance\HL@@C\@ne
   \xdef\Thelabel@@@{\number\HL@@C}%
   \xdefThelabel@{\HL@@N}%
   \xdef\Thelabel@@@@{\HL@@P\Thelabel@\HL@@Q}%
   \xdefThelabel@@{\HL@@S}%
  \else
   \global\advance\HL@@@C\@ne
   \xdef\Thelabel@@@{\number\HL@@@C}%
   \xdefThelabel@{\HL@@@N}%
   \xdef\Thelabel@@@@{\HL@@@P\Thelabel@\HL@@@Q}%
   \xdefThelabel@@{\HL@@@S}%
  \fi
  \xdef\Thepref@{\Thelabel@@@@}}%
  \csname HL@\HLlevel@\endcsname##1\endHL
  \let\pref\Thepref@
  \csname HL@I\HLlevel@\endcsname
  \csname HL@J\HLlevel@\endcsname
  \let\pref\pref@
  \HLtoc@
  \aftertoc@
  \let\aftertoc@\relax\overlong@false}%
 \ifx\next"\expandafter\next@\else\expandafter\nextii@\fi}%
\Invalid@\endHL
\def\hl@@C{\csname hl@C\hllevel@\endcsname}
\def\hl@@P{\csname hl@P\hllevel@\endcsname}
\def\hl@@Q{\csname hl@Q\hllevel@\endcsname}
\def\hl@@S{\csname hl@S\hllevel@\endcsname}
\def\hl@@N{\csname hl@N\hllevel@\endcsname}
\def\hl@@F{\csname hl@F\hllevel@\endcsname}
\def\hl@@@C{\csname\exxx@\hltype@ @C\endcsname}
\def\hl@@@P{\csname\exxx@\hltype@ @P\endcsname}
\def\hl@@@Q{\csname\exxx@\hltype@ @Q\endcsname}
\def\hl@@@S{\csname\exxx@\hltype@ @S\endcsname}
\def\hl@@@N{\csname\exxx@\hltype@ @N\endcsname}
\def\hl#1{\expandafter
 \ifx\csname hl@C#1\endcsname\relax
  \DN@{\Err@{\string\hl#1 not defined in this style}}%
 \else
  \DN@{\gdef\hllevel@{#1}\def\hlname@{\hl{#1}}\let\hltype@\relax\FNSS@\hl@}%
 \fi
 \next@}
\def\hl@{%
 \DN@"##1"##2{\def\entry@{##2}\quoted@true
  {\noexpands@
  \ifx\hltype@\relax
   \let\pre\hl@@P\let\post\hl@@Q\let\style\hl@@S\let\numstyle\hl@@N
  \else
   \let\pre\hl@@@P\let\post\hl@@@Q\let\style\hl@@@S\let\numstyle\hl@@@N
  \fi
  \Qlabel@{##1}\let\style\relax\xdef\Qlabel@@@@{##1}%
  \xdef\Thepref@{\Thelabel@@@@}}%
  \csname hl@\hllevel@\endcsname{##2}%
  \let\pref\Thepref@
  \csname hl@I\hllevel@\endcsname
  \csname hl@J\hllevel@\endcsname
  \let\pref\pref@
  \hltoc@
  \aftertoc@
  \let\aftertoc@\relax\nopunct@false\nospace@false\FNSSP@}%
 \DNii@##1{\def\entry@{##1}\quoted@false
  {\noexpands@
  \ifx\hltype@\relax
   \global\advance\hl@@C\@ne
   \xdef\Thelabel@@@{\number\hl@@C}%
   \xdefThelabel@{\hl@@N}%
   \xdef\Thelabel@@@@{\hl@@P\Thelabel@\hl@@Q}%
   \xdefThelabel@@{\hl@@S}%
  \else
   \global\advance\hl@@@C\@ne
   \xdef\Thelabel@@@{\number\hl@@@C}%
   \xdefThelabel@{\hl@@@N}%
   \xdef\Thelabel@@@@{\hl@@@P\Thelabel@\hl@@@Q}%
   \xdefThelabel@@{\hl@@@S}%
  \fi
  \xdef\Thepref@{\Thelabel@@@@}}%
  \csname hl@\hllevel@\endcsname{##1}%
  \let\pref\Thepref@
  \csname hl@I\hllevel@\endcsname
  \csname hl@J\hllevel@\endcsname
  \let\pref\pref@
  \hltoc@
  \aftertoc@
  \let\aftertoc@\relax\nopunct@false\nospace@false\FNSSP@}%
 \ifx\next"\expandafter\next@\else\expandafter\nextii@\fi}%
\def\six@#1#2 #3 #4 #5 #6 #7 {\DN@{#2}\ifx\next@\empty
 \DN@##1\six@{}\else
 \write#1{ #2 #3 #4 #5 #6 #7}\DN@{\six@#1}\fi
 \next@}
\def\Sixtoc@{\ifx\macdef@\empty\else
 \DN@##1##2\next@{\def\macdef@{##1##2}}%
 \expandafter\next@\macdef@\next@
 \edef\next@
  {\noexpand\six@\toc@\macdef@
  \space\space\space\space\space\space\space\space\space\space\space\space
  \noexpand\six@}%
 \next@\let\macdef@\relax\fi}
\def\QorThelabel@@@@{\ifquoted@
 \noexpand\noexpand\noexpand"\Qlabel@@@@\noexpand\noexpand\noexpand"\else
 \Thelabel@@@@\fi}
\def\HLtoc@{%
 \iftoc@
 \expandafter\expandafter\expandafter\unmacro@
  \expandafter\meaning\csname HL@W\HLlevel@\endcsname\unmacro@
  {\noexpands@\let\style\relax
   \edef\next@{\write\toc@{\noexpand\noexpand\expandafter\noexpand\HLname@
   {\macdef@}{\QorThelabel@@@@}}}%
  \next@}%
  \expandafter\unmacro@\meaning\entry@\unmacro@
  \Sixtoc@
  \write\toc@{\noexpand\Page{\number\pageno}{\page@N}%
   {\page@P}{\page@Q}^^J}%
 \fi}
\def\hltoc@{%
 \iftoc@
 \expandafter\expandafter\expandafter\unmacro@
  \expandafter\meaning\csname hl@W\hllevel@\endcsname\unmacro@
  {\noexpands@\let\style\relax
  \edef\next@{\write\toc@{%
   \ifnopunct@\noexpand\noexpand\noexpand\nopunct\fi
   \ifnospace@\noexpand\noexpand\noexpand\nospace\fi
   \noexpand\noexpand\expandafter\noexpand\hlname@
   {\macdef@}{\QorThelabel@@@@}}}%
  \next@}%
  \expandafter\unmacro@\meaning\entry@\unmacro@
  \Sixtoc@
  \write\toc@{\noexpand\Page{\number\pageno}{\page@N}%
   {\page@P}{\page@Q}^^J}%
 \fi}
\def\mainfile#1{\def\mainfile@{#1}}
\def\checkmainfile@{\ifx\mainfile@\undefined
 \Err@{No \noexpand\mainfile specified}\fi}
\expandafter\newcount@\csname HL@C1\endcsname
\csname HL@C1\endcsname\z@
\expandafter\def\csname HL@S1\endcsname#1{#1\null.}
\expandafter\let\csname HL@N1\endcsname\arabic
\expandafter\let\csname HL@P1\endcsname\empty
\expandafter\let\csname HL@Q1\endcsname\empty
\expandafter\def\csname HL@F1\endcsname{\bf}
\expandafter\let\csname HL@W1\endcsname\empty
\expandafter\newcount@\csname hl@C1\endcsname
\csname hl@C1\endcsname\z@
\expandafter\def\csname hl@S1\endcsname#1{#1\/}
\expandafter\let\csname hl@N1\endcsname\arabic
\expandafter\let\csname hl@P1\endcsname\empty
\expandafter\let\csname hl@Q1\endcsname\empty
\expandafter\def\csname hl@F1\endcsname{\bf}
\expandafter\let\csname hl@W1\endcsname\empty
\expandafter\def\csname HL@1\endcsname#1\endHL{\bigbreak
 {\locallabel@
  \global\setbox\@ne\vbox{\Let@\tabskip\hss@
  \halign to\hsize{\bf\hfil\ignorespaces##\unskip\hfil\cr
  \expandafter\ifx\csname HL@W1\endcsname\empty\else
   \csname HL@W1\endcsname\space\fi
  {\HL@@F\ifx\thelabel@@\empty\else\thelabel@@\space\fi}%
  \ignorespaces#1\crcr}}%
  }%
 \unvbox\@ne\nobreak\medskip}
\expandafter\def\csname hl@1\endcsname#1{\medbreak\noindent@@
 {\locallabel@
 \bf{\hl@@F\ifx\thelabel@@\empty\else\thelabel@@\space\fi}%
 \ignorespaces#1\unskip\punct@{\null.}\addspace@\enspace}}
\expandafter\def\csname HL@I1\endcsname{\Reset\hl1{1}%
 \ifx\pref\empty\newpre\hl1{}\else\newpre\hl1{\pref.}\fi}
\def\NameHL#1#2{\define#2{}%
 \expandafter\ifx\csname HL@R#1\endcsname\relax
 \else
  \def\nextiv@{\let\nextiii@}%
  \expandafter\nextiv@\csname HL@R#1\endcsname
  \expandafter\let\nextiii@\undefined
  \expandafter\let\csname\exxx@\nextiii@ @C\endcsname\relax
  \expandafter\let\csname\exxx@\nextiii@ @P\endcsname\relax
  \expandafter\let\csname\exxx@\nextiii@ @Q\endcsname\relax
  \expandafter\let\csname\exxx@\nextiii@ @S\endcsname\relax
  \expandafter\let\csname\exxx@\nextiii@ @N\endcsname\relax
  \expandafter\let\csname\exxx@\nextiii@ @F\endcsname\relax
  \expandafter\let\csname\exxx@\nextiii@ @W\endcsname\relax
  \expandafter\let\csname end\exxx@\nextiii@\endcsname\undefined
 \fi
 \expandafter\gdef\csname HL@R#1\endcsname{#2}%
 \expandafter\gdef\csname\exstring@#2@R\endcsname{{HL}{#1}}%
 \iftoc@\write\toc@{\noexpand\NameHL#1\noexpand#2^^J}\fi
 \rightadd@#2\to\overlonglist@
 \edef\next@{\let\csname\exstring@#2@C\endcsname\expandafter\noexpand
  \csname HL@C#1\endcsname}\next@
 \edef\next@{\let\csname\exstring@#2@P\endcsname\expandafter\noexpand
  \csname HL@P#1\endcsname}\next@
 \edef\next@{\let\csname\exstring@#2@Q\endcsname\expandafter\noexpand
  \csname HL@Q#1\endcsname}\next@
 \edef\next@{\let\csname\exstring@#2@S\endcsname\expandafter\noexpand
  \csname HL@S#1\endcsname}\next@
 \edef\next@{\let\csname\exstring@#2@N\endcsname\expandafter\noexpand
  \csname HL@N#1\endcsname}\next@
 \edef\next@{\let\csname\exstring@#2@F\endcsname\expandafter\noexpand
  \csname HL@F#1\endcsname}\next@
 \edef\next@{\let\csname\exstring@#2@W\endcsname\expandafter\noexpand
  \csname HL@W#1\endcsname}\next@
 \edef\next@{\def\noexpand#2####1\expandafter\noexpand
  \csname end\exstring@#2\endcsname
  {\def\noexpand\HLtype@{\noexpand#2}%
   \def\noexpand\HLname@{\noexpand#2}%
   \gdef\noexpand\HLlevel@{#1}%
   \noexpand\FNSS@\noexpand\HL@####1\noexpand\endHL}}%
  \next@
 \edef\next@{\noexpand\Invalid@\expandafter\noexpand
  \csname end\exstring@#2\endcsname}%
 \next@}
\def\Namehl#1#2{\define#2{}%
 \expandafter\ifx\csname hl@R#1\endcsname\relax
 \else
  \def\nextiv@{\let\nextiii@}%
  \expandafter\nextiv@\csname hl@R#1\endcsname
  \expandafter\let\nextiii@\undefined
  \expandafter\let\csname\exxx@\nextiii@ @C\endcsname\relax
  \expandafter\let\csname\exxx@\nextiii@ @P\endcsname\relax
  \expandafter\let\csname\exxx@\nextiii@ @Q\endcsname\relax
  \expandafter\let\csname\exxx@\nextiii@ @S\endcsname\relax
  \expandafter\let\csname\exxx@\nextiii@ @N\endcsname\relax
  \expandafter\let\csname\exxx@\nextiii@ @F\endcsname\relax
  \expandafter\let\csname\exxx@\nextiii@ @W\endcsname\relax
 \fi
 \expandafter\gdef\csname hl@R#1\endcsname{#2}%
 \expandafter\gdef\csname\exstring@#2@R\endcsname{{hl}{#1}}%
 \iftoc@\write\toc@{\noexpand\Namehl#1\noexpand#2^^J}\fi
 \rightadd@#2\to\nofrillslist@%
 \edef\next@{\let\csname\exstring@#2@C\endcsname\expandafter\noexpand
  \csname hl@C#1\endcsname}\next@
 \edef\next@{\let\csname\exstring@#2@P\endcsname\expandafter\noexpand
  \csname hl@P#1\endcsname}\next@
 \edef\next@{\let\csname\exstring@#2@Q\endcsname\expandafter\noexpand
  \csname hl@Q#1\endcsname}\next@
 \edef\next@{\let\csname\exstring@#2@S\endcsname\expandafter\noexpand
  \csname hl@S#1\endcsname}\next@
 \edef\next@{\let\csname\exstring@#2@N\endcsname\expandafter\noexpand
  \csname hl@N#1\endcsname}\next@
 \edef\next@{\let\csname\exstring@#2@F\endcsname\expandafter\noexpand
  \csname hl@F#1\endcsname}\next@
 \edef\next@{\let\csname\exstring@#2@W\endcsname\expandafter\noexpand
  \csname hl@W#1\endcsname}\next@
 \edef\next@{\def\noexpand#2{%
  \def\noexpand\hltype@{\noexpand#2}%
  \def\noexpand\hlname@{\noexpand#2}%
  \gdef\noexpand\hllevel@{#1}%
  \noexpand\FNSS@\noexpand\hl@}}%
 \next@}%
\def\Initialize{\FN@\Init@}
\def\Init@{\ifx\next\HL\let\next@\InitH@\else\ifx\next\hl\let\next@\InitH@
  \else\let\next@\InitS@\fi\fi\next@}
\def\InitH@#1#2{\expandafter\ifx\csname\exstring@#1@C#2\endcsname\relax
 \DN@{\Err@{\noexpand#1level #2 not defined in this style}}\else
 \DN@{\expandafter\gdef\csname\exstring@#1@J#2\endcsname}\fi\next@}
\def\InitC@#1#2{\edef\nextii@{\expandafter\noexpand\csname#1\endcsname{#2}}}
\def\InitS@#1{\expandafter\ifx\csname\exstring@#1@R\endcsname\relax
 \Err@{\noexpand#1not defined in this style}\let\next@\relax\else
 \DN@{\let\next@}\expandafter\next@\csname\exstring@#1@R\endcsname
 \expandafter\InitC@\next@
 \DN@{\expandafter\InitH@\nextii@}\fi\next@}
\def\value#1{\expandafter
 \ifx\csname\exstring@#1@C\endcsname\relax
  \expandafter\ifx\csname\exstring@#1@C1\endcsname\relax
   \DN@{\Err@{\noexpand\value can't be used with \string#1}}%
  \else
   \DN@{\value@#1}%
  \fi
 \else
  \DN@{\number\csname\exstring@#1@C\endcsname\relax}%
 \fi
 \next@}
\def\value@#1#2{\expandafter
 \ifx\csname\exstring@#1@C#2\endcsname\relax
  \DN@{\Err@{\string\value\string#1 can't be followed by \string#2}}%
 \else
  \DN@{\number\csname\exstring@#1@C#2\endcsname\relax}%
 \fi
 \next@}
\newcount\Value
\def\Evaluate#1{\expandafter
 \ifx\csname\exstring@#1@C\endcsname\relax
  \expandafter\ifx\csname\exstring@#1@C1\endcsname\relax
   \DN@{\Err@{\noexpand\Evaluate can't be used with \string#1}}%
  \else
   \DN@{\Evaluate@#1}%
  \fi
 \else
  \DN@{\global\Value\csname\exstring@#1@C\endcsname}%
 \fi
 \next@}
\def\Evaluate@#1#2{\expandafter
 \ifx\csname\exstring@#1@C#2\endcsname\relax
  \DN@{\Err@{\string\Evaluate\string#1 can't be followed by \string#2}}%
 \else
  \DN@{\global\Value\csname\exstring@#1@C#2\endcsname}%
 \fi\next@}
\def\pre#1{\expandafter
 \ifx\csname\exstring@#1@P\endcsname\relax
  \expandafter\ifx\csname\exstring@#1@P1\endcsname\relax
   \DN@{\Err@{\noexpand\pre can't be used with \string#1}}%
  \else
   \DN@{\pre@#1}%
  \fi
 \else
  \DN@{{\csname\exstring@#1@P\endcsname}}%
 \fi
 \next@}
\def\pre@#1#2{\expandafter
 \ifx\csname\exstring@#1@P#2\endcsname\relax
  \DN@{\Err@{\string\pre\string#1 can't be followed by \string#2}}%
 \else
  \DN@{{\csname\exstring@#1@P#2\endcsname}}%
 \fi
 \next@}
\def\post#1{\expandafter
 \ifx\csname\exstring@#1@Q\endcsname\relax
  \expandafter\ifx\csname\exstring@#1@Q1\endcsname\relax
   \DN@{\Err@{\noexpand\post can't be used with \string#1}}%
  \else
   \DN@{\post@#1}%
  \fi
 \else
  \DN@{{\csname\exstring@#1@Q\endcsname}}%
 \fi
 \next@}
\def\post@#1#2{\expandafter
 \ifx\csname\exstring@#1@Q#2\endcsname\relax
  \DN@{\Err@{\string\post\string#1 can't be followed by \string#2}}%
 \else
  \DN@{{\csname\exstring@#1@Q#2\endcsname}}%
 \fi
 \next@}
\def\style#1{\expandafter
 \ifx\csname\exstring@#1@S\endcsname\relax
  \expandafter\ifx\csname\exstring@#1@S1\endcsname\relax
   \DN@{\Err@{\noexpand\style can't be used with \string#1}}%
  \else
   \DN@{\style@#1}%
  \fi
 \else
  \DN@{\csname\exstring@#1@S\endcsname}%
 \fi
 \next@}
\def\style@#1#2{\expandafter
 \ifx\csname\exstring@#1@S#2\endcsname\relax
  \DN@{\Err@{\string\style\string#1 can't be followed by \string#2}}%
 \else
  \DN@{\csname\exstring@#1@S#2\endcsname}%
 \fi
 \next@}
\def\fontstyle#1{\expandafter
 \ifx\csname\exstring@#1@F\endcsname\relax
  \expandafter\ifx\csname\exstring@#1@F1\endcsname\relax
   \DN@{\Err@{\noexpand\fontstyle can't be used with \string#1}}%
  \else
   \DN@{\fontstyle@#1}%
  \fi
 \else
  \DN@##1{{\csname\exstring@#1@F\endcsname##1}}%
 \fi
 \next@}
\def\fontstyle@#1#2{\expandafter
 \ifx\csname\exstring@#1@F#2\endcsname\relax
  \DN@{\Err@{\string\fontstyle\string#1 can't be followed by \string#2}}%
 \else
  \DN@##1{{\csname\exstring@#1@F#2\endcsname##1}}%
 \fi
 \next@}
\def\Reset#1{\expandafter
 \ifx\csname\exstring@#1@C\endcsname\relax
  \expandafter\ifx\csname\exstring@#1@C1\endcsname\relax
   \DN@{\Err@{\noexpand\Reset can't be used with \string#1}}%
  \else
   \DN@{\Reset@#1}%
  \fi
 \else
  \DN@##1{\count@##1\relax\ifx#1\page\else\advance\count@\m@ne\fi
   \global\csname\exstring@#1@C\endcsname\count@}%
 \fi
 \next@}
\def\Reset@#1#2{\expandafter
 \ifx\csname\exstring@#1@C#2\endcsname\relax
  \DN@{\Err@{\string\Reset\string#1 can't be followed by \string#2}}%
 \else
  \DN@##1{\count@##1\relax\advance\count@\m@ne
   \global\csname\exstring@#1@C#2\endcsname\count@}%
 \fi
 \next@}
\def\Offset#1{\expandafter
 \ifx\csname\exstring@#1@C\endcsname\relax
  \expandafter\ifx\csname\exstring@#1@C1\endcsname\relax
   \DN@{\Err@{\noexpand\Offset can't be used with \string#1}}%
  \else
   \DN@{\Offset@#1}%
  \fi
 \else
  \DN@##1{\count@##1\relax\advance\count@\m@ne\global\advance
   \csname\exstring@#1@C\endcsname\count@}%
 \fi
 \next@}
\def\Offset@#1#2{\expandafter
 \ifx\csname\exstring@#1@C#2\endcsname\relax
  \DN@{\Err@{\string\Offset\string#1 can't be followed by \string#2}}%
 \else
  \DN@##1{\count@##1\relax\advance\count@\m@ne
   \global\advance\csname\exstring@#1@C#2\endcsname\count@}%
 \fi
 \next@}
\def\getR@#1#2{\def\nextiv@{\let\nextiii@}\expandafter\nextiv@
 \csname\exstring@#1@R#2\endcsname}
\def\letR@#1#2#3{\expandafter\let\csname#1@#3#2\endcsname\Next@}
\def\letR@@#1#2{\expandafter\let\csname\exstring@#1@#2\endcsname\Next@}
\def\newpre#1{\expandafter
 \ifx\csname\exstring@#1@P\endcsname\relax
  \expandafter\ifx\csname\exstring@#1@P1\endcsname\relax
   \DN@{\Err@{\noexpand\newpre can't be used with \string#1}}%
  \else
   \DN@{\newpre@#1}%
  \fi
 \else
  \DN@{%
   \DNii@{%
    \endgroup
    \expandafter\let\csname\exstring@#1@P\endcsname\Next@
    \expandafter\ifx\csname\exstring@#1@R\endcsname\relax\else
    \getR@#1{}\expandafter\letR@\nextiii@ P\fi
    }%
   \begingroup\noexpands@\afterassignment\nextii@\xdef\Next@}%
 \fi
 \next@}
\def\newpre@#1#2{\expandafter
 \ifx\csname\exstring@#1@P#2\endcsname\relax
  \DN@{\Err@{\string\newpre\string#1 can't be followed by \string#2}}%
 \else
  \DN@{%
   \DNii@{%
    \endgroup
    \expandafter\let\csname\exstring@#1@P#2\endcsname\Next@
    \expandafter\ifx\csname\exstring@#1@R#2\endcsname\relax\else
    \getR@#1{#2}\expandafter\letR@@\nextiii@ P\fi
    }%
   \begingroup\noexpands@\afterassignment\nextii@\xdef\Next@}%
 \fi
 \next@}
\def\newpost#1{\expandafter
 \ifx\csname\exstring@#1@Q\endcsname\relax
  \expandafter\ifx\csname\exstring@#1@Q1\endcsname\relax
   \DN@{\Err@{\noexpand\newpost can't be used with \string#1}}%
  \else
   \DN@{\newpost@#1}%
  \fi
 \else
  \DN@{%
   \DNii@{%
    \endgroup
    \expandafter\let\csname\exstring@#1@Q\endcsname\Next@
    \expandafter\ifx\csname\exstring@#1@R\endcsname\relax\else
    \getR@#1{}\expandafter\letR@\nextiii@ Q\fi
    }%
   \begingroup\noexpands@\afterassignment\nextii@\xdef\Next@}%
 \fi
 \next@}
\def\newpost@#1#2{\expandafter
 \ifx\csname\exstring@#1@Q#2\endcsname\relax
  \DN@{\Err@{\string\newpost\string#1 can't be followed by \string#2}}%
 \else
  \DN@{%
   \DNii@{%
    \endgroup
    \expandafter\let\csname\exstring@#1@Q#2\endcsname\Next@
    \expandafter\ifx\csname\exstring@#1@R#2\endcsname\relax\else
    \getR@#1{#2}\expandafter\letR@@\nextiii@ Q\fi
    }%
   \begingroup\noexpands@\afterassignment\nextii@\xdef\Next@}%
 \fi
 \next@}
\def\newstyle#1{\expandafter
 \ifx\csname\exstring@#1@S\endcsname\relax
  \expandafter\ifx\csname\exstring@#1@S1\endcsname\relax
   \DN@{\Err@{\noexpand\newstyle can't be used
    with \string#1}}%
  \else
   \DN@{\newstyle@#1}%
  \fi
 \else
  \DN@{%
   \DNii@{%
    \expandafter\let\csname\exstring@#1@S\endcsname\Next@
    \expandafter\ifx\csname\exstring@#1@R\endcsname\relax\else
    \getR@#1{}\expandafter\letR@\nextiii@ S\fi
    }%
   \afterassignment\nextii@\gdef\Next@}%
 \fi
 \next@}
\def\newstyle@#1#2{\expandafter
 \ifx\csname\exstring@#1@S#2\endcsname\relax
  \DN@{\Err@{\string\newstyle\string#1 can't be followed by
   \string#2}}%
 \else
  \DN@{%
   \DNii@{%
    \expandafter\let\csname\exstring@#1@S#2\endcsname\Next@
    \expandafter\ifx\csname\exstring@#1@R#2\endcsname\relax\else
    \getR@#1{#2}\expandafter\letR@@\nextiii@ S\fi
    }%
   \afterassignment\nextii@\gdef\Next@}%
 \fi
 \next@}
\def\newnumstyle#1{\expandafter
 \ifx\csname\exstring@#1@N\endcsname\relax
  \expandafter\ifx\csname\exstring@#1@N1\endcsname\relax
   \DN@{\Err@{\noexpand\newnumstyle can't be used with
    \string#1}}%
  \else
   \DN@{\newnumstyle@#1}%
  \fi
 \else
  \DN@##1{%
   \gdef\Next@{##1}%
    \expandafter\let\csname\exstring@#1@N\endcsname\Next@
    \expandafter\ifx\csname\exstring@#1@R\endcsname\relax\else
    \getR@#1{}\expandafter\letR@\nextiii@ N\fi
    }%
 \fi
 \next@}
\def\newnumstyle@#1#2{\expandafter
 \ifx\csname\exstring@#1@N#2\endcsname\relax
  \DN@{\Err@{\string\newnumstyle\string#1 can't be followed by
   \string#2}}%
 \else
  \DN@##1{%
   \gdef\Next@{##1}%
    \expandafter\let\csname\exstring@#1@N#2\endcsname\Next@
    \expandafter\ifx\csname\exstring@#1@R#2\endcsname\relax\else
    \getR@#1{#2}\expandafter\letR@@\nextiii@ N\fi
    }%
  \fi
 \next@}
\def\newfontstyle#1{\expandafter
 \ifx\csname\exstring@#1@F\endcsname\relax
  \expandafter\ifx\csname\exstring@#1@F1\endcsname\relax
   \DN@{\Err@{\noexpand\newfontstyle can't be used with
    \string#1}}%
  \else
   \DN@{\newfontstyle@#1}%
  \fi
 \else
  \DN@##1{%
   \gdef\Next@{##1}%
    \expandafter\let\csname\exstring@#1@F\endcsname\Next@
    \expandafter\ifx\csname\exstring@#1@R\endcsname\relax\else
    \getR@#1{}\expandafter\letR@\nextiii@ F\fi
    }%
 \fi
 \next@}
\def\newfontstyle@#1#2{\expandafter
 \ifx\csname\exstring@#1@F#2\endcsname\relax
  \DN@{\Err@{\string\newfontstyle\string#1 can't be followed by
   \string#2}}%
 \else
  \DN@##1{%
   \gdef\Next@{##1}%
    \expandafter\let\csname\exstring@#1@F#2\endcsname\Next@
    \expandafter\ifx\csname\exstring@#1@R#2\endcsname\relax\else
    \getR@#1{#2}\expandafter\letR@@\nextiii@ F\fi
    }%
 \fi
 \next@}
\def\word#1{\expandafter
 \ifx\csname\exstring@#1@W\endcsname\relax
  \expandafter\ifx\csname\exstring@#1@W1\endcsname\relax
   \DN@{\Err@{\noexpand\word can't be used with \string#1}}%
  \else
   \DN@{\word@#1}%
  \fi
 \else
  \DN@{{\csname\exstring@#1@W\endcsname}}%
 \fi
 \next@}
\def\word@#1#2{\expandafter
 \ifx\csname\exstring@#1@W#2\endcsname\relax
  \DN@{\Err@{\string\word\noexpand#1can't be followed by \string#2}}%
 \else
  \DN@{{\csname\exstring@#1@W#2\endcsname}}%
 \fi
 \next@}
\def\newword#1{\expandafter
 \ifx\csname\exstring@#1@W\endcsname\relax
  \expandafter\ifx\csname\exstring@#1@W1\endcsname\relax
   \DN@{\Err@{\noexpand\newword can't be used  with \string#1}}%
  \else
   \DN@{\newword@#1}%
  \fi
 \else
  \DN@{%
   \DNii@{%
    \expandafter\let\csname\exstring@#1@W\endcsname\Next@
    \expandafter\ifx\csname\exstring@#1@R\endcsname\relax\else
     \getR@#1{}\expandafter\letR@\nextiii@ W\fi
    }%
   \afterassignment\nextii@\gdef\Next@}%
 \fi
 \next@}
\def\newword@#1#2{\expandafter
 \ifx\csname\exstring@#1@W#2\endcsname\relax
  \DN@{\Err@{\string\newword\noexpand#1can't be followed by \string#2}}%
 \else
  \DN@{%
   \DNii@{%
    \expandafter\let\csname\exstring@#1@W#2\endcsname\Next@
    \expandafter\ifx\csname\exstring@#1@R#2\endcsname\relax\else
     \getR@#1{#2}\expandafter\letR@@\nextiii@ W\fi
    }%
   \afterassignment\nextii@\gdef\Next@}%
 \fi
 \next@}
\newif\iffn@
\newcount\footmark@C
\footmark@C\z@
\def\footmark@S#1{$^{#1}$}
\let\footmark@N\arabic
\def\footmark@F{\rm}
\def\foottext@S#1{$^{#1}$}
\def\foottext@F{\rm}
\let\modifyfootnote@\relax
\def\modifyfootnote#1{\def\modifyfootnote@{#1}}
\def\vfootnote@#1{\insert\footins
 \bgroup
 \floatingpenalty\@MM\interlinepenalty\interfootnotelinepenalty
 \leftskip\z@\rightskip\z@\spaceskip\z@\xspaceskip\z@
 \rm\splittopskip\ht\strutbox\splitmaxdepth\dp\strutbox
 \locallabel@\noindent@@{\foottext@F#1}\modifyfootnote@
 \footstrut\FN@\fo@t}
\def\fo@t{\ifcat\bgroup\noexpand\next\expandafter\f@@t\else
 \expandafter\f@t\fi}
\def\f@t#1{#1\@foot}
\def\f@@t{\bgroup\aftergroup\@foot\afterassignment\FNSSP@\let\next@}
\def\@foot{\unskip\lower\dp\strutbox\vbox to\dp\strutbox{}\egroup
 \iffn@\expandafter\fn@false\else
 \expandafter\postvanish@\fi}
\newif\ifplainfn@
\plainfn@true
\def\fancyfootnotes{\plainfn@false}
\newcount\fancyfootmarkcount@
\fancyfootmarkcount@\z@
\newcount\lastfnpage@
\lastfnpage@-\@M
\let\justfootmarklist@\empty
\def\footmark{\let\@sf\empty
 \ifhmode\edef\@sf{\spacefactor\the\spacefactor}\/\fi
 \DN@{\ifx"\next\expandafter\nextii@\else\expandafter\footmark@\fi}%
 \DNii@"##1"{%
  \iffirstchoice@
   {\let\style\footmark@S\let\numstyle\footmark@N
   \footmark@F##1%
   \noexpands@
   \let\style\foottext@S
   \Qlabel@{##1}%
   }%
   \iffn@\else
    {\noexpands@
    \xdef\Next@{{\Thelabel@}{\Thelabel@@}{\Thelabel@@@}{\Thelabel@@@@}}%
    }%
    \expandafter\rightappend@\Next@\to\justfootmarklist@
   \fi
  \fi
  \@sf\relax}%
 \FN@\next@}
\def\footmark@{%
 \iffirstchoice@
  \global\advance\footmark@C\@ne
  \ifplainfn@
   \xdef\adjustedfootmark@{\number\footmark@C}%
  \else
   {\let\\\or\xdef\Next@{\ifcase\number\footmark@C\fnpages@\else
     -\@M\fi}}%
   \ifnum\Next@=-\@M
    \xdef\adjustedfootmark@{\number\footmark@C}%
   \else
    \ifnum\Next@=\lastfnpage@
     \global\advance\fancyfootmarkcount@\@ne
    \else
     \global\fancyfootmarkcount@\@ne
     \global\lastfnpage@\Next@
    \fi
    \xdef\adjustedfootmark@{\number\fancyfootmarkcount@}%
   \fi
  \fi
  {\noexpands@
  \xdef\Thelabel@@@{\adjustedfootmark@}%
  \xdefThelabel@\footmark@N
  \xdef\Thelabel@@@@{\Thelabel@}%
  \xdefThelabel@@\foottext@S
  }%
  \iffn@\else
   {\noexpands@
   \xdef\Next@{{\Thelabel@}{\Thelabel@@}{\Thelabel@@@}{\Thelabel@@@@}}%
   }%
   \expandafter\rightappend@\Next@\to\justfootmarklist@
  \fi
  \ifplainfn@
  \else
   \edef\next@{\write\laxwrite@{F\noexpand\the\pageno}}\next@
  \fi
 \fi
 \footmark@S{\footmark@N{\adjustedfootmark@}}%
 \@sf\relax}
\def\foottext{\prevanish@
 \ifx\justfootmarklist@\empty
  \Err@{There is no \noexpand\footmark for this \string\foottext}\fi
 \DN@\\##1##2\next@{\DN@{##1}\gdef\justfootmarklist@{##2}}%
 \expandafter\next@\justfootmarklist@\next@
 \expandafter\foottext@\next@}
\def\foottext@#1#2#3#4{{\noexpands@
  \xdef\Thelabel@{#1}\xdef\Thelabel@@{#2}%
  \xdef\Thelabel@@@{#3}\xdef\Thelabel@@@@{#4}}%
  \vfootnote@{\thelabel@@}}
\rightadd@\foottext\to\vanishlist@
\def\footnote{\fn@true
 \let\@sf\empty
 \ifhmode\edef\@sf{\spacefactor\the\spacefactor}\/\fi
 \DN@{\ifx"\next\expandafter\nextii@\else\expandafter\nextiii@\fi}%
 \DNii@"##1"{\footmark"##1"\vfootnote@{\let\style\foottext@S
  \let\numstyle\footmark@N##1}}%
 \def\nextiii@{\footmark\vfootnote@{\foottext@S{\footmark@N
  {\adjustedfootmark@}}}}%
 \FN@\next@}
\newdimen\litindent
\litindent20\p@
\newbox\litbox@
\newbox\Litbox@
\newcount\interlitpenalty@
\interlitpenalty@\@M
\newcount\litlines@
{\obeyspaces\gdef\defspace@{\def {\allowbreak\hskip.5emminus.15em}}}
{\obeylines\gdef\letM@{\let^^M\CtrlM@}}
\def\CtrlM@{\egroup
 \ifcase\litlines@\advance\litlines@\@ne\or
 \box\litbox@\advance\litlines@\@ne\else
 \penalty\interlitpenalty@\box\litbox@\fi
 \Lit@}
\def\Lit@{\setbox\litbox@\hbox\bgroup\litdefs@\hskip\litindent}
\newcount\littab@
\littab@8
\def\littab#1{\littab@#1\relax}
{\catcode`\^^I=\active\gdef\letTAB@{\let^^I\TAB@}}
\def\TAB@{\egroup
 \dimen@\wd\litbox@
 \advance\dimen@-\litindent
 \setboxz@h{\tt0}%
 \dimen@ii\littab@\wdz@
 \divide\dimen@\dimen@ii
 \multiply\dimen@\dimen@ii
 \advance\dimen@\littab@\wdz@
 \advance\dimen@\litindent
 \setbox\litbox@\hbox\bgroup\litdefs@\hbox to\dimen@{\unhbox\litbox@\hfil}}
{\catcode`\`=\active\gdef`{\relax\lq}}
\let\litbs@\relax
\let\litbs@@\relax
\def\litbackslash#1{%
 \edef\litbs@{\catcode`\string#1=\z@
 \def\noexpand\litbs@@{\def\expandafter\noexpand\csname\string#1\endcsname
  {\char`\string#1}}}}
\def\litcodes@{\catcode`\\=12
 \catcode`\{=12 \catcode`\}=12
 \catcode`\$=12 \catcode`\&=12
 \catcode`\#=12
 \catcode`\^=12 \catcode`\_=12
 \catcode`\@=12 \catcode`\~=12 \catcode`\"=12
 \catcode`\;=12 \catcode`\:=12 \catcode`\!=12 \catcode`\?=12
 \catcode`\%=12 \litbs@\catcode`\`=\active\obeyspaces\defspace@}
\def\activate@#1#2{{\lccode`\~=`#2%
 \lowercase{%
  \if0#1%
  \gdef\Next@{\def~{\egroup\endgroup\bigskip\vskip-\parskip
   \def\next@{\noindent@@\FN@\pretendspace@}\FNSS@\next@}}\else
  \gdef\Next@{\def~{\egroup\egroup\endgroup}}\fi
  }%
 }}
\def\litdefs@{\let\0\empty\let\1\litdelim@\def\ {\char32 }\litbs@@}%
\def\litdelimiter#1{%
 \edef\litdelim@{\char`#1}%
 \def\lit#1{\leavevmode\begingroup\litcodes@\litdefs@
  \tt\hyphenchar\tentt\m@ne\lit@}%
 \def\lit@##1#1{##1\endgroup\null}%
 \def\Lit#1{\ifhmode$$\abovedisplayskip\bigskipamount
  \abovedisplayshortskip\bigskipamount
  \belowdisplayskip\z@\belowdisplayshortskip\z@
  \postdisplaypenalty\@M
  $$\vskip-\baselineskip\else\bigskip\fi
  \begingroup\litlines@\z@
  \catcode`#1=\active\activate@0#1\Next@
  \def\displaybreak{\egroup\break\litlines@\z@\Lit@}%
  \def\allowdisplaybreak{\egroup\allowbreak\litlines@\z@\Lit@}%
  \def\allowdisplaybreaks{\egroup\allowbreak\interlitpenalty@\z@
   \litlines@\z@\Lit@}%
  \litcodes@\tt\catcode`\^^I=\active\letTAB@
  \obeylines\letM@\Lit@}%
 \def\Litbox##1=#1{\begingroup\ifodd##1\relax\aftergroup\global\fi
  \aftergroup\setbox\aftergroup##1\aftergroup\box\aftergroup\Litbox@
  \def\allowdisplaybreak{\egroup\allowbreak\litlines@\z@\Lit@}%
  \def\allowdisplaybreaks{\egroup\allowbreak\interlitpenalty@\z@
   \litlines@\z@\Lit@}%
  \catcode`#1=\active\activate@1#1\Next@
  \litcodes@\tt\catcode`\^^I=\active\letTAB@
  \obeylines\letM@\global\setbox\Litbox@\vbox\bgroup\litindent\z@%
  \litlines@\z@\Lit@}%
}
\newbox\titlebox@
\setbox\titlebox@\vbox{}
\rightadd@\title\to\overlonglist@
\def\title{\begingroup\Let@
 \global\setbox\titlebox@\vbox\bgroup\tabskip\hss@
 \halign to\hsize\bgroup\bf\hfil\ignorespaces##\unskip\hfil\cr}
\def\endtitle{\crcr\egroup\egroup\endgroup\overlong@false}
\newbox\authorbox@
\rightadd@\author\to\overlonglist@
\def\author{\begingroup\Let@
 \global\setbox\authorbox@\vbox\bgroup\tabskip\hss@
 \halign to\hsize\bgroup\rm\hfil\ignorespaces##\unskip\hfil\cr}
\def\endauthor{\crcr\egroup\egroup\endgroup\overlong@false}
\newbox\affilbox@
\def\affil{\begingroup\Let@
 \global\setbox\affilbox@\vbox\bgroup\tabskip\hss@
 \halign to\hsize\bgroup\rm\hfil\ignorespaces##\unskip\hfil\cr}%
\def\endaffil{\crcr\egroup\egroup\endgroup\overlong@false}
\let\date@\relax
\def\date#1{\gdef\date@{\ignorespaces#1\unskip}}
\def\today{\ifcase\month\or January\or February\or March\or April\or May\or
 June\or July\or August\or September\or October\or November\or December\fi
 \space\number\day, \number\year}
\def\maketitle{\hrule\height\z@\vskip-\topskip
 \vskip24\p@ plus12\p@ minus12\p@
 \unvbox\titlebox@
 \ifvoid\authorbox@\else\vskip12\p@ plus6\p@ minus3\p@\unvbox\authorbox@\fi
 \ifvoid\affilbox@\else\vskip10\p@ plus5\p@ minus2\p@\unvbox\affilbox@\fi
 \ifx\date@\relax\else\vskip6\p@ plus2\p@ minus\p@\centerline{\rm\date@}\fi
 \vskip18\p@ plus12\p@ minus6\p@}
\def\cite{%
 \DNii@(##1)##2{{\rm[}{##2}, {##1\/}{\rm]}}%
 \def\nextiii@##1{{\rm[}{##1\/}{\rm]}}%
 \DN@{\ifx\next(\expandafter\nextii@\else\expandafter\nextiii@\fi}%
 \FN@\next@}
\def\makebib@W{Bibliography}
\def\makebib{\begingroup\rm\bigbreak\centerline{\smc\makebib@W}%
 \nobreak\medskip
 \sfcode`\.=\@m\everypar{}\parindent\z@
 \def\nopunct{\nopunct@true}\def\nospace{\nospace@true}%
 \nopunct@false\nospace@false
 \def\lkerns@{\null\kern\m@ne sp\kern\@ne sp}%
 \def\nkerns@{\null\kern-\tw@ sp\kern\tw@ sp}%
}

\newif\ifnoprepunct@
\newif\ifnoprespace@
\newif\ifnoquotes@
\def\noprepunct{\noprepunct@true}
\def\noprespace{\noprespace@true}
\def\noquotes{\noquotes@true}
\newbox\nobox@
\newbox\keybox@
\newbox\bybox@
\newbox\paperbox@
\newbox\paperinfobox@
\newbox\jourbox@
\newbox\volbox@
\newbox\issuebox@
\newbox\yrbox@
\newbox\pgbox@
\newbox\ppbox@
\newbox\bookbox@
\newbox\inbookbox@
\newbox\bookinfobox@
\newbox\publbox@
\newbox\publaddrbox@
\newbox\edbox@
\newbox\edsbox@
\newbox\langbox@
\newbox\translbox@
\newbox\finalinfobox@
\def\setbibinfo@#1{\edef\next@{\ifnopunct@1\else0\fi
 \ifnospace@1\else0\fi\ifnoprepunct@1\else0\fi\ifnoprespace@1\else0\fi
 \ifnoquotes@1\else0\fi}%
 \DNii@{00000}%
 \ifx\next@\nextii@\else\xdef\bibinfo@{\bibinfo@\the#1,\next@}%
 \fi}
\def\getbibinfo@#1{\ifx\bibinfo@\empty
 \let\next@0\let\nextii@0\let\nextiii@0\let\nextiv@0\let\nextv@0\else
 \edef\next@{\def
  \noexpand\next@####1\the#1,####2####3####4####5####6####7\noexpand\next@
  {\let\noexpand\next@####2\let\noexpand\nextii@####3%
  \let\noexpand\nextiii@####4\let\noexpand\nextiv@####5%
  \let\noexpand\nextv@####6}%
  \noexpand\next@\bibinfo@\the#1,00000\noexpand\next@}\next@
 \fi}
\newif\ifbookinquotes@
\def\bookinquotes{\bookinquotes@true}
\newif\ifpaperinquotes@
\def\paperinquotes{\paperinquotes@true}
\newif\ifininbook@
\def\ininbook{\ininbook@true}
\newif\ifopenquotes@
\def\closequotes@{\ifopenquotes@''\openquotes@false\fi}
\newif\ifbeginbib@
\newif\ifendbib@
\newif\ifprevjour@
\newif\ifprevbook@
\newdimen\bibindent@
\bibindent@20\p@
\def\bib{\global\let\bibinfo@\empty\global\let\translinfo@\relax\beginbib@true
 \begingroup\noindent@
 \hangindent\bibindent@\hangafter\@ne\bib@}
\def\v@id#1{\setbox#1\box\voidb@x}
\def\bib@{\v@id\nobox@\v@id\keybox@\v@id\bybox@\v@id\paperbox@
 \v@id\paperinfobox@\v@id\jourbox@\v@id\volbox@\v@id\issuebox@
 \v@id\yrbox@\v@id\pgbox@\v@id\ppbox@\v@id\bookbox@\v@id\inbookbox@
 \v@id\bookinfobox@\v@id\publbox@\v@id\publaddrbox@\v@id\edbox@
 \v@id\edsbox@\v@id\langbox@\v@id\translbox@\v@id\finalinfobox@
 \bgroup}
\def\Setnonemptybox@#1#2{\unskip\setbibinfo@#1\egroup#2%
 \def\aftergroup@{\ifdim\wd#1=\z@\setbox#1\box\voidb@x\fi}%
 \setbox#1\vbox\bgroup\aftergroup\aftergroup@\hsize\maxdimen\leftskip\z@
 \rightskip\z@\hbadness\@M\hfuzz\maxdimen\noindent}
\def\setnonemptybox@#1{\Setnonemptybox@#1\relax}
\def\no{\setnonemptybox@\nobox@}
\def\key{\setnonemptybox@\keybox@\bf}
\def\by{\setnonemptybox@\bybox@}
\def\bysame{\setnonemptybox@\bybox@\leaders\hrule\hskip3em\null}
\def\paper{\setnonemptybox@\paperbox@
 \ifpaperinquotes@\getbibinfo@\paperbox@
 \if\nextv@1\else``\fi\else\it\fi}
\def\paperinfo{\setnonemptybox@\paperinfobox@}
\def\jour{\Setnonemptybox@\jourbox@\prevjour@true}
\def\vol{\setnonemptybox@\volbox@\bf}
\def\issue{\setnonemptybox@\issuebox@}
\def\yr{\setnonemptybox@\yrbox@}
\def\toappear{\noprepunct\finalinfo(to appear)}
\def\pg{\setnonemptybox@\pgbox@}
\def\pp{\setnonemptybox@\ppbox@}
\def\book{\Setnonemptybox@\bookbox@\prevbook@true
 \ifbookinquotes@\getbibinfo@\bookbox@
 \if\nextv@1\else``\fi\else\it\fi}
\def\inbook{\Setnonemptybox@\inbookbox@\prevbook@true
 \ifininbook@ in \fi\ifbookinquotes@\getbibinfo@\inbookbox@
 \if\nextv@1\else``\fi\fi}
\def\bookinfo{\setnonemptybox@\bookinfobox@}
\def\publ{\setnonemptybox@\publbox@}
\def\publaddr{\setnonemptybox@\publaddrbox@}
\def\ed{\setnonemptybox@\edbox@}
\def\eds{\setnonemptybox@\edsbox@}
\def\lang{\setnonemptybox@\langbox@}
\def\finalinfo{\setnonemptybox@\finalinfobox@}
\def\setboxzl@{\setbox\z@\lastbox}
\def\getbox@#1{\setbox\z@\vbox{\vskip-\@M\p@
 \unvbox#1%
 \setboxzl@
 \global\setbox\@ne\hbox{\unhbox\z@\unskip\unskip\unpenalty}%
 \ifdim\lastskip=-\@M\p@\else
 \loop\ifdim\lastskip=-\@M\p@
 \else\unskip\unpenalty\setboxzl@
 \global\setbox\@ne\hbox{\unhbox\z@\unhbox\@ne}%
 \repeat\fi}%
 \unhbox\@ne}
\def\adjustpunct@#1{\count@\lastkern
 \ifnum\count@=\z@#1\closequotes@\else
 \ifnum\count@>\tw@#1\closequotes@\else
 \ifnum\count@<-\tw@#1\closequotes@\else
  \unkern\unkern\setboxzl@
  \skip@\lastskip\unskip
  \count@@\lastpenalty\unpenalty
  \ifnum\count@=\tw@\unskip\setboxzl@\fi
  \ifdim\skip@=\z@\else\hskip\skip@\fi
  #1\closequotes@
  \ifnum\count@=\tw@\null\hfill\fi
  \penalty\count@@
 \fi\fi\fi}
\def\prepunct@#1#2{\getbibinfo@#2%
 \ifnopunct@
 \else
  \if\nextiii@0\adjustpunct@#1\fi
 \fi
 \closequotes@
 \ifnospace@
 \else
  \if\nextiv@0\space\else\fi
 \fi
 \nopunct@false\nospace@false
 \if\next@1\nopunct@true\fi
 \if\nextii@1\nospace@true\fi}
\def\ppunbox@#1#2{\prepunct@{#1}#2%
 \getbox@#2}
\let\semicolon@;
\def\endbib@{%
 \ifbeginbib@
  \ifvoid\nobox@
   \ifvoid\keybox@\else\hbox to\bibindent@{[\getbox@\keybox@]\hss}\fi
  \else\hbox to\bibindent@{\hss\getbox@\nobox@. }\fi
  \ifvoid\bybox@\else\getbox@\bybox@\fi
 \else
  \nopunct@true
  \ifvoid\bybox@\else\ppunbox@\relax\bybox@\fi
 \fi
 \ifvoid\translbox@\else\ppunbox@,\translbox@\fi
 \ifvoid\paperbox@\else\ppunbox@,\paperbox@\ifpaperinquotes@
  \if\nextv@1\else\openquotes@true\fi\fi
 \fi
 \ifvoid\paperinfobox@\else\ppunbox@,\paperinfobox@\fi
 \test@false
 \ifvoid\jourbox@\else\test@true\ppunbox@,\jourbox@\fi
 \ifprevjour@\test@true\fi
 \iftest@
  \ifvoid\volbox@\else\ppunbox@\relax\volbox@\fi
  \ifvoid\issuebox@
   \else\prepunct@\relax\issuebox@ no.~\getbox@\issuebox@\fi
  \ifvoid\yrbox@\else\prepunct@\relax\yrbox@(\getbox@\yrbox@)\fi
  \ifvoid\ppbox@\else\ppunbox@,\ppbox@\fi
  \ifvoid\pgbox@\else\prepunct@,\pgbox@ p.~\getbox@\pgbox@\fi
 \fi
 \test@false
 \ifvoid\bookbox@\else\test@true\ppunbox@,\bookbox@\ifbookinquotes@
  \if\nextv@1\else\openquotes@true\fi\fi\fi
 \ifvoid\inbookbox@\else\test@true\ppunbox@,\inbookbox@\ifbookinquotes@
  \if\nextv@1\else\openquotes@true\fi\fi\fi
 \ifprevbook@\test@true\fi
 \iftest@
  \ifvoid\edbox@\else\prepunct@\relax\edbox@(\getbox@\edbox@, ed.)\fi
  \ifvoid\edsbox@\else\prepunct@\relax\edsbox@(\getbox@\edsbox@, eds.)\fi
  \ifvoid\bookinfobox@\else\ppunbox@,\bookinfobox@\fi
  \ifvoid\publbox@\else\ppunbox@,\publbox@\fi
  \ifvoid\publaddrbox@\else\ppunbox@,\publaddrbox@\fi
  \ifvoid\yrbox@\else\ppunbox@,\yrbox@\fi
  \ifvoid\ppbox@\else\prepunct@,\ppbox@ pp.~\getbox@\ppbox@\fi
  \ifvoid\pgbox@\else\prepunct@,\pgbox@ p.~\getbox@\pgbox@\fi
 \fi
 \ifvoid\finalinfobox@
  \ifendbib@
   \ifnopunct@\else.\closequotes@\fi
  \else
  \ifvoid\langbox@\else\space(\getbox@\langbox@)\fi
   \/\semicolon@\closequotes@
  \fi
 \else
  \ifendbib@
   \ppunbox@{.\spacefactor3000\relax}\finalinfobox@
    \ifnopunct@\else.\fi
  \else
   \ppunbox@,\finalinfobox@\/\semicolon@\fi
 \fi
 \ifvoid\langbox@\else\space(\getbox@\langbox@)\fi
}
\def\endbib{\unskip\egroup\endbib@true\endbib@\par\endgroup}
\def\morebib{\unskip\egroup
 \endbib@false\endbib@
 \global\let\bibinfo@\empty\beginbib@false
 \bib@}
\def\anotherbib{\unskip\egroup
 \endbib@false\endbib@
 \global\let\bibinfo@\empty\beginbib@false
 \prevjour@false\prevbook@false\bib@}
\def\transl{\unskip
 \xdef\translinfo@{\the\translbox@,\ifnopunct@1\else0\fi
 \ifnospace@1\else0\fi\ifnoprepunct@1\else0\fi\ifnoprespace@1\else0\fi0}%
 \egroup\endbib@false\endbib@
 \global\let\bibinfo@\translinfo@\beginbib@false
 \bib@
 \egroup
 \def\aftergroup@{\ifdim\wd\translbox@=\z@\setbox\translbox@\box\voidb@x\fi}%
 \setbox\translbox@\vbox\bgroup\aftergroup\aftergroup@
 \hsize\maxdimen\leftskip\z@\rightskip\z@\hbadness\@M\hfuzz\maxdimen
 \noindent}
\newwrite\auxwrite@
\newread\bbl@
\def\UseBibTeX{\immediate\openout\auxwrite@=\jobname.aux
 \let\cite\BTcite@
 \def\nocite##1{\immediate\write\auxwrite@{\string\citation{##1}}}%
 \def\bibliographystyle##1{\immediate\write\auxwrite@{\string
  \bibstyle{##1}}}%
 \def\bibliography@W{Bibliography}%
 \def\bibliography##1{\immediate\write\auxwrite@{\string\bibdata{##1}}%
  \immediate\openin\bbl@=\jobname.bbl
  \ifeof\bbl@
   \W@{No .bbl file}%
  \else
   \immediate\closein\bbl@
   \begingroup\input bibtex \input\jobname.bbl \endgroup
  \fi}%
 }
\def\BTcite@{%
 \DNii@(##1)##2{{\rm[}\BTcite@@##2,\BTcite@@{\rm, }{##1\/}{\rm]}%
  \immediate\write\auxwrite@{\string\citation{##2}}}%
 \def\nextiii@##1{{\rm[}\BTcite@@##1,\BTcite@@\/{\rm]}%
  \immediate\write\auxwrite@{\string\citation{##1}}}%
 \DN@{\ifx\next(\expandafter\nextii@\else\expandafter\nextiii@\fi}%
 \FN@\next@}%
\def\BTcite@@#1,{\BTcite@@@{#1}\FN@\BTcite@@@@}
\def\BTcite@@@@{\ifx\next\BTcite@@
 \expandafter\eat@\else{\rm, }\expandafter\BTcite@@\fi}
\catcode`\~=11
\def\BTcite@@@#1{\nolabel@\cite{#1}\relax
 \DNii@##1~##2\nextii@{##1}%
 \csL@{#1}\expandafter\nextii@\Next@\nextii@\fi}
\catcode`\~=\active

\def\beginthebibliography@#1{\rm\setboxz@h{#1\ }\bibindent@\wdz@
 \bigbreak\centerline{\smc\bibliography@W}\nobreak\medskip
 \sfcode`\.=\@m\everypar{}\parindent\z@}

\newif\iffigproofing@
\def\Figureproofing{\figproofing@true}
\def\noFigureproofing{\figproofing@false}
\newif\ifHby@
\def\Hbyw#1{\global\Hby@true\hbyw\vsize{#1}}
\def\hbyw#1#2{%
 \hbox{%
  \ifHby@
  \else
   \iffigproofing@
    \setbox\z@\vbox{\hrule\width5\p@}\ht\z@\z@
    \vbox to#1{\hrule\height5\p@\width.4\p@\vfil\hrule\height5\p@\width.4\p@}%
    \kern-.4\p@\rlap{\copy\z@}\raise#1\hbox{\rlap{\copy\z@}}%
   \fi
  \fi
  \vbox to#1{\hbox to#2{}\vfil}%
  \ifHby@
  \else
   \iffigproofing@
    \vbox to#1{\hrule\height5\p@\width.4\p@\vfil\hrule\height5\p@\width.4\p@}%
    \kern-.4\p@\llap{\copy\z@}\raise#1\hbox{\llap{\boxz@}}%
   \fi
  \fi}}
\newcount\island@C
\let\island@P\empty
\let\island@Q\empty
\def\island@S#1{#1\null.}
\let\island@N\arabic
\def\island@F{\rm}
\def\island@@@P{\csname\exxx@\islandtype@ @P\endcsname}
\def\island@@@Q{\csname\exxx@\islandtype@ @Q\endcsname}
\def\island@@@S{\csname\exxx@\islandtype@ @S\endcsname}
\def\island@@@N{\csname\exxx@\islandtype@ @N\endcsname}
\def\island@@@F{\csname\exxx@\islandtype@ @F\endcsname}
\def\island@@@C{\csname island@C\islandclass@\endcsname}
\newif\ifplace@
\newif\ifisland@
\def\island{%
 \ifplace@
  \DN@{\let\islandclass@\empty\def\islandtype@{\island}\FN@\island@}%
 \else
  \long\DN@##1\endisland{\Err@{\noexpand\island must be used after some
   type of \string\...place}}%
 \fi
 \next@}
\def\island@{\ifx\next\c\let\next@\island@c\else
 \DN@{\FN@\island@@}\fi\next@}
\def\island@@{\ifcat\bgroup\noexpand\next\let\next@\island@@@\else
 \DN@{\Err@{\noexpand\island must be followed by a {prefix} for
 \string\caption's}}\fi\next@}
\newbox\islandbox@
\newcount\captioncount@
\def\island@@@#1{\def\captionprefix@{#1}\captioncount@\z@
 \global\setbox\islandbox@\vbox\bgroup}
\def\island@c\c#1{%
 \ifplace@
 \DN@{\def\islandclass@{#1}%
  \expandafter\ifx\csname island@C#1\endcsname\relax
  \expandafter\newcount@\csname island@C#1\endcsname
   \global\csname island@C#1\endcsname\z@\fi
  \FNSS@\island@c@}%
 \else
 \DN@{\edef\next@{\long\def\noexpand\next@########1\expandafter\noexpand
  \csname end\exxx@\islandtype@\endcsname{\noexpand\Err@{\noexpand\noexpand
  \expandafter\noexpand
  \islandtype@ must be used after some type of \noexpand\string
   \noexpand\...place}}}\next@\next@}%
 \fi
 \next@}
\def\island@c@{%
 \ifcat\bgroup\noexpand\next
  \let\next@\island@c@@
 \else
  \DN@{\Err@{\noexpand\island\string\c{\expandafter\string\islandclass@} must
   be followed by a {prefix} for \string\caption's}}%
 \fi\next@}
\def\island@c@@#1{\def\captionprefix@{#1}%
 \captioncount@\z@\global\setbox\islandbox@\vbox\bgroup}
\rightadd@\caption\to\nofrillslist@
\newbox\captionbox@
\newbox\Captionbox@
\def\caption{%
 \ifnum\captioncount@=\z@
  \ifnopunct@
   \DN@{\egroup\nopunct@true}%
  \else
   \let\next@\egroup
  \fi
 \else
  \let\next@\relax
 \fi
 \next@
 \advance\captioncount@\@ne
 \FN@\caption@}
\def\caption@{\ifx\next"\expandafter\caption@q\else\expandafter\caption@@\fi}
\def\caption@q"#1"{\quoted@true
 {\noexpands@
 \let\pre\island@@@P\let\post\island@@@Q
 \let\style\island@@@S\let\numstyle\island@@@N
 \Qlabel@{#1}\let\style\relax\xdef\Qlabel@@@@{#1}}%
 \finishcaption@}
\def\caption@@{\quoted@false
 \global\advance\island@@@C\@ne
 {\noexpands@
 \xdef\Thelabel@@@{\number\island@@@C}%
 \xdefThelabel@\island@@@N
 \xdef\Thelabel@@@@{\island@@@P\Thelabel@\island@@@Q}%
 \xdefThelabel@@\island@@@S
 \xdef\Thepref@{\Thelabel@@@@}}%
 \finishcaption@}
\long\def\captionformat@#1#2#3{\rm\strut#1 {\island@@@F#2} #3%
 \punct@.\strut}
\long\def\widerthanisland@#1#2#3{\test@true\setbox\z@\vbox{\hsize\maxdimen
 \noindent@@\captionformat@{#1}{#2}{#3}\par\setboxzl@}%
 \ifdim\wdz@=\z@
  \global\setbox\captionbox@\hbox{\noset@\unlabel@
   \captionformat@{#1}{#2}{#3}}%
  \ifdim\wd\captionbox@>\wd\islandbox@\else\test@false\fi
 \fi}
\long\def\captionformat@@#1#2#3{\widerthanisland@{#1}{#2}{#3}%
 \iftest@
  \global\setbox\captionbox@\vbox{\hsize\wd\islandbox@
   \vskip-\parskip\noindent@@\noset@\unlabel@
   \captionformat@{#1}{#2}{#3}\par}%
 \else
  \global\setbox\captionbox@
   \hbox to\wd\islandbox@{\hfil\box\captionbox@\hfil}%
 \fi}
\long\def\finishcaption@#1{\def\entry@{#1}%
 {\locallabel@
 \captionformat@@
  {\expandafter\ignorespaces\captionprefix@\unskip}%
  {\ifx\thelabel@@\empty\unskip\else\thelabel@@\fi}%
  {\ignorespaces#1\unskip}%
 \ifnum\captioncount@=\@ne
  \global\setbox\islandbox@\vbox{\ticwrite@\vbox{\box\islandbox@}}%
  \global\setbox\Captionbox@\vbox{\box\captionbox@}%
 \else
  \global\setbox\islandbox@\vbox{\unvbox\islandbox@\setboxzl@
   \ticwrite@\boxz@}%
  \global\setbox\Captionbox@\vbox{\unvbox\Captionbox@
   \smallskip\box\captionbox@}%
 \fi}%
 \nopunct@false\nospace@false\ignorespaces}
\def\Sixtic@{\ifx\macdef@\empty\else
 \DN@##1##2\next@{\def\macdef@{##1##2}}%
 \expandafter\next@\macdef@\next@
 \edef\next@
  {\noexpand\six@\tic@\macdef@
  \space\space\space\space\space\space\space\space\space\space\space\space
  \noexpand\six@}%
 \next@\let\macdef@\relax\fi}
\def\ticwrite@{%
 \iftoc@
  {\noexpands@\let\style\relax
  \DN@{\island}%
  \edef\next@{\write\tic@{%
   \ifnopunct@\noexpand\noexpand\noexpand\nopunct\fi
   \ifx\islandtype@\next@\noexpand\noexpand\noexpand\island
    \noexpand\string\noexpand\c{\islandclass@}{\captionprefix@}%
     {\QorThelabel@@@@}\else\noexpand\noexpand\expandafter\noexpand
     \islandtype@{\QorThelabel@@@@}}\fi}%
  \next@}%
  \expandafter\unmacro@\meaning\entry@\unmacro@
  \Sixtic@
  \write\tic@{\noexpand\Page{\number\pageno}{\page@N}{\page@P}{\page@Q}^^J}%
 \fi}
\def\Htrim@#1{%
 \ifHby@
  \dimen@\vsize
  \ifnum\captioncount@=\z@
  \else
   \advance\dimen@-\ht\Captionbox@
   \advance\dimen@-#1%
  \fi
  \global\Hby@false
  \dimen@ii\wd\islandbox@
  \global\setbox\islandbox@\vbox
   {\unvbox\islandbox@\setboxzl@
   \vbox to\z@{\vss\boxz@}\nointerlineskip\hbyw\dimen@\dimen@ii}%
  \global\Hby@true
 \fi}
\newif\ifdata@
\def\iclasstest@#1{\DN@{#1}\ifx\next@\islandclass@
 \test@true\else\test@false\fi}
\skipdef\skipi@=1
\def\endisland{\ifnum\captioncount@=\z@\expandafter\egroup\fi
 \ifdata@
 \else
  \iclasstest@{T}%
  \iftest@
   {\rm\global\skipi@-\dp\strutbox}\global\advance\skipi@\bigskipamount
   \Htrim@\skipi@
   \global\setbox\islandbox@\vbox
    {\ifnum\captioncount@=\z@\else
     \box\Captionbox@
     \nointerlineskip
     \vskip\skipi@\fi
     \box\islandbox@}%
  \else
   {\rm\global\skipi@\dp\strutbox}\global\advance\skipi@\medskipamount
   \Htrim@\skipi@
   \global\setbox\islandbox@\vbox
    {\box\islandbox@
     \ifnum\captioncount@=\z@\else
     \nointerlineskip
     \vskip\skipi@
     \box\Captionbox@
     \fi}%
  \fi
  \ifHby@
  \else
   \dimen@\ht\islandbox@\advance\dimen@\dp\islandbox@
   \ifdim\dimen@>\vsize
    \DN@{\island}%
    \Err@{%
     \ifx\islandtype@\next@\noexpand\island\else
      \expandafter\noexpand\islandtype@\fi
     \ifnum\captioncount@=\z@\else
       with \noexpand\caption\fi
      is larger than page}%
     \ht\islandbox@=\vsize
   \fi
  \fi
 \fi
 \global\Hby@false\island@true}
\def\newisland#1\c#2#3{\define#1{}%
 \iftoc@\immediate\write\tic@{\noexpand\newisland\noexpand#1%
  \string\c{#2}{#3}^^J}\fi
 \expandafter\def\csname\exstring@#1@S\endcsname{\island@S}%
 \expandafter\def\csname\exstring@#1@N\endcsname{\island@N}%
 \expandafter\def\csname\exstring@#1@P\endcsname{\island@P}%
 \expandafter\def\csname\exstring@#1@Q\endcsname{\island@Q}%
 \expandafter\def\csname\exstring@#1@F\endcsname{\island@F}%
 \expandafter\def\csname end\exstring@#1\endcsname{\endisland}%
 \expandafter
 \ifx\csname island@C#2\endcsname\relax
  \expandafter\newcount@\csname island@C#2\endcsname
  \global\csname island@C#2\endcsname\z@
 \fi
 \edef\next@{\noexpand\expandafter\noexpand\let\noexpand
  \csname\exstring@#1@C\noexpand\endcsname
  \csname island@C#2\endcsname}%
 \next@
 \def#1{\def\islandtype@{#1}\island@c\c{#2}{#3}}}
\newisland\Figure\c{F}{Figure}
\newisland\Table\c{T}{Table}
\newbox\islandboxi
\newbox\islandboxii
\newbox\islandboxiii
\newbox\captionboxi
\newbox\captionboxii
\newbox\captionboxiii
\long\def\islandpairdata#1#2{{\data@true
 \place@true
 #1%
 \global\setbox\islandboxi\box\islandbox@
 \global\setbox\captionboxi\box\Captionbox@
 #2%
 \global\setbox\islandboxii\box\islandbox@
 \global\setbox\captionboxii\box\Captionbox@
 }}
\long\def\islandpairbox#1#2{\islandpairdata{#1}{#2}%
 \dimen@\ht\captionboxi
 \ifdim\ht\captionboxii>\dimen@\dimen@\ht\captionboxii\fi
 \ifdim\dimen@>\z@
  \ifdim\ht\captionboxi<\dimen@
   \global\setbox\captionboxi\vbox to\dimen@{\unvbox\captionboxi\vfil}\fi
  \ifdim\ht\captionboxii<\dimen@
   \global\setbox\captionboxii\vbox to\dimen@{\unvbox\captionboxii\vfil}\fi
 \fi
 \global\setbox\islandbox@\vbox
 {\hbox to\hsize{\hfil\box\islandboxi\hfil\box\islandboxii\hfil}%
 \ifdim\dimen@>\z@\nointerlineskip
 {\rm\global\skipi@\dp\strutbox}\global\advance\skipi@\medskipamount
  \vskip\skipi@
  \hbox to\hsize{\hfil\box\captionboxi\hfil\box\captionboxii\hfil}\fi}}	
\long\def\islandpairboxa#1#2{\islandpairdata{#1}{#2}%
 \dimen@\ht\captionboxi
 \ifdim\ht\captionboxii>\dimen@\dimen@\ht\captionboxii\fi
 \ifdim\dimen@>\z@
  \ifdim\ht\captionboxi<\dimen@
   \global\setbox\captionboxi\vbox to\dimen@{\vfil\unvbox\captionboxi}\fi
  \ifdim\ht\captionboxii<\dimen@
   \global\setbox\captionboxii\vbox to\dimen@{\vfil\unvbox\captionboxii}\fi
 \fi
 \dimen@ii\ht\islandboxi
 \ifdim\ht\islandboxii>\dimen@ii \dimen@ii\ht\islandboxii\fi
 \ifdim\dimen@ii>\z@
  \ifdim\ht\islandboxi<\dimen@ii
   \global\setbox\islandboxi\vbox to\dimen@ii{\box\islandboxi\vfil}\fi
  \ifdim\ht\islandboxii<\dimen@ii
   \global\setbox\islandboxii\vbox to\dimen@ii{\box\islandboxii\vfil}\fi
 \fi
 \global\setbox\islandbox@\vbox{\ifdim\dimen@>\z@
  \hbox to\hsize{\hfil\box\captionboxi\hfil\box\captionboxii\hfil}%
  \nointerlineskip{\rm\global\skipi@-\dp\strutbox}%
  \global\advance\skipi@\bigskipamount\vskip\skipi@\fi
  \hbox to\hsize{\hfil\box\islandboxi\hfil\box\islandboxii\hfil}}}
\long\def\islandtripledata#1#2#3{{\data@true\place@true
 #1%
 \global\setbox\islandboxi\box\islandbox@
 \global\setbox\captionboxi\box\Captionbox@
 #2%
 \global\setbox\islandboxii\box\islandbox@
 \global\setbox\captionboxii\box\Captionbox@
 #3%
 \global\setbox\islandboxiii\box\islandbox@
 \global\setbox\captionboxiii\box\Captionbox@
 }}
\long\def\islandtriplebox#1#2#3{\islandtripledata{#1}{#2}{#3}%
 \dimen@\ht\captionboxi
 \ifdim\ht\captionboxii>\dimen@ \dimen@\ht\captionboxii\fi
 \ifdim\ht\captionboxiii>\dimen@ \dimen@\ht\captionboxiii\fi
 \ifdim\dimen@>\z@
  \ifdim\ht\captionboxi<\dimen@
   \global\setbox\captionboxi\vbox to\dimen@{\unvbox\captionboxi\vfil}\fi
  \ifdim\ht\captionboxii<\dimen@
   \global\setbox\captionboxii\vbox to\dimen@{\unvbox\captionboxii\vfil}\fi
  \ifdim\ht\captionboxiii<\dimen@
   \global\setbox\captionboxiii\vbox to\dimen@{\unvbox\captionboxiii\vfil}\fi
 \fi
 \global\setbox\islandbox@\vbox
  {\hbox to\hsize{\hfil\box\islandboxi\hfil\box\islandboxii\hfil
   \box\islandboxiii\hfil}%
 \ifdim\dimen@>\z@\nointerlineskip
  {\rm\global\skipi@\dp\strutbox}\global\advance\skipi@\medskipamount
  \vskip\skipi@
  \hbox to\hsize{\hfil\box\captionboxi\hfil\box\captionboxii\hfil
   \box\captionboxiii\hfil}\fi}}
\def\islandtripleboxa#1#2#3{\islandtripledata{#1}{#2}{#3}%
 \dimen@\ht\captionboxi
 \ifdim\ht\captionboxii>\dimen@ \dimen@\ht\captionboxii\fi
 \ifdim\ht\captionboxiii>\dimen@ \dimen@\ht\captionboxiii\fi
 \ifdim\dimen@>\z@
  \ifdim\ht\captionboxi<\dimen@
   \global\setbox\captionboxi\vbox to\dimen@{\vfil\unvbox\captionboxi}\fi
  \ifdim\ht\captionboxii<\dimen@
   \global\setbox\captionboxii\vbox to\dimen@{\vfil\unvbox\captionboxii}\fi
  \ifdim\ht\captionboxiii<\dimen@
   \global\setbox\captionboxiii\vbox to\dimen@{\vfil\unvbox\captionboxiii}\fi
 \fi
 \dimen@ii\ht\islandboxi
 \ifdim\ht\islandboxii>\dimen@ii \dimen@ii\ht\islandboxii\fi
 \ifdim\ht\islandboxiii>\dimen@ii \dimen@ii\ht\islandboxiii\fi
 \ifdim\dimen@ii>\z@
  \ifdim\ht\islandboxi<\dimen@ii
   \global\setbox\islandboxi\vbox to\dimen@ii{\box\islandboxi\vfil}\fi
  \ifdim\ht\islandboxii<\dimen@ii
   \global\setbox\islandboxii\vbox to\dimen@ii{\box\islandboxii\vfil}\fi
  \ifdim\ht\islandboxiii<\dimen@ii
   \global\setbox\islandboxiii\vbox to\dimen@ii{\box\islandboxiii\vfil}\fi
 \fi
 \global\setbox\islandbox@\vbox
  {\ifdim\dimen@>\z@
  \hbox to\hsize{\hfil\box\captionboxi\hfil\box\captionboxii\hfil
   \box\captionboxiii\hfil}%
  \nointerlineskip{\rm\global\skipi@-\dp\strutbox}%
  \global\advance\skipi@\bigskipamount\vskip\skipi@\fi
  \hbox to\hsize{\hfil\box\islandboxi\hfil\box\islandboxii\hfil
   \box\islandboxiii\hfil}}}
\def\Figurepair#1\and#2\endFigurepair{\island@true
 \islandpairbox{\Figure#1\endFigure}{\Figure#2\endFigure}}
\def\Figuretriple#1\and#2\and#3\endFiguretriple{\island@true
 \islandtriplebox{\Figure#1\endFigure}{\Figure#2\endFigure}%
  {\Figure#3\endFigure}}
\def\Tablepair#1\and#2\endTablepair{\island@true
 \islandpairboxa{\Table#1\endTable}{\Table#2\endTable}}
\def\Tabletriple#1\and#2\and#3\endTabletriple{\island@true
 \islandtripleboxa{\Table#1\endTable}{\Table#2\endTable}%
 {\Table#3\endTable}}
\def\place#1{\place@true\island@false
 #1%
 \ifisland@
  \box\islandbox@
 \else
  \Err@{Whoa ... there's no \string\Figure, \string\Table,
   etc., here}%
 \fi
 \place@false}
\newskip\belowtopfigskip
\belowtopfigskip 15\p@ plus 5\p@ minus5\p@
\newskip\abovebotfigskip
\abovebotfigskip 18\p@ plus 6\p@ minus6\p@
\newdimen\minpagesize
\minpagesize 5pc
\dimen@\belowtopfigskip
\advance\dimen@-\abovebotfigskip
\skip\topins\dimen@
\dimen\topins\z@
\newcount\topinscount@
\newbox\topinsdims@
\def\storedim@{\global\setbox\topinsdims@
 \vbox{\hbox to\dimen@{}\unvbox\topinsdims@}}
\def\advancedimtopins@{%
 \ifnum\pageno=\@ne
 \else
   \advance\dimen@\dimen\topins
   \global\dimen\topins\dimen@
 \fi}
\newcount\flipcount@
\def\fliptopins@{%
 \global\flipcount@\z@
 \ifvoid\topins\else
 \setbox\z@\vbox
  {\vskip\p@
   \unvbox\topins
   \global\setbox\topins\vbox{}%
   \loop
    \test@false
    \ifdim\lastskip=\z@\unskip
     \ifdim\lastskip=\z@
      \test@true\fi\fi
    \iftest@
    \global\advance\flipcount@\@ne
    \setboxzl@
    \global\setbox\topins\vbox{\unvbox\topins\boxz@}%
    \unpenalty
   \repeat}\fi}
\newif\ifPar@
\newcount\Parcount@
\newbox\Parbox@
\expandafter\newbox\csname Parfigbox1\endcsname
\expandafter\newbox\csname Parfigbox2\endcsname
\expandafter\newbox\csname Parfigbox3\endcsname
\expandafter\newbox\csname Parfigbox4\endcsname
\expandafter\newbox\csname Parfigbox5\endcsname
\expandafter\newdimen\csname Parprev1\endcsname
\expandafter\newdimen\csname Parprev2\endcsname
\expandafter\newdimen\csname Parprev3\endcsname
\expandafter\newdimen\csname Parprev4\endcsname
\expandafter\newdimen\csname Parprev5\endcsname
\expandafter\newdimen\csname Parprev6\endcsname
\def\Par{\par\global\csname Parprev1\endcsname\prevdepth
 \global\Parcount@\@ne
 \global\Par@true\global\let\Parlist@\empty
 \global\setbox\Parbox@\vbox\bgroup\break}
\def\place@#1#2{%
 \ifisland@
  \ifhmode
   \ifPar@
    \ifnum\Parcount@>5
     \Err@{Only 5 \string\place's allowed per
      \string\Par...\noexpand\endPar paragraph}%
    \else
     \expandafter\expandafter\expandafter
      \global\expandafter\setbox
       \csname Parfigbox\number\Parcount@\endcsname\box\islandbox@
     \global\advance\Parcount@\@ne
     \xdef\Parlist@{\Parlist@#1}%
    \fi
   \else
    \vadjust{#2}%
   \fi
  \else
   #2%
  \fi
 \else
  \Err@{Whoa ... there's no \string\Figure,
   \string\Table, etc., here}%
 \fi
 \place@false}
\long\def\Aplace#1{\prevanish@
 \place@true\island@false
 #1%
 \place@ a\Aplace@
 \postvanish@}
\long\def\AAplace#1{\prevanish@\place@true\island@false
 #1%
 \place@ A\AAplace@
 \postvanish@}
\newif\ifAA@
\def\AAplace@{\AA@true\Aplace@\AA@false}
\let\AAlist@\empty
\def\Aplace@{\allowbreak
 \dimen@=\ht\islandbox@
 \advance\dimen@\abovebotfigskip
 \ht\islandbox@\dimen@
 \advance\dimen@\dp\islandbox@
 \storedim@
 \ifAA@
  \xdef\AAlist@{\AAlist@1}%
  \advancedimtopins@
 \else
  \xdef\AAlist@{\AAlist@0}%
  \ifnum\topinscount@>\@ne\else\advancedimtopins@\fi
 \fi
 \insert\topins{\penalty\z@\splittopskip\z@\floatingpenalty\z@
  \box\islandbox@}%
 \global\advance\topinscount@\@ne}
\long\def\Bplace#1{\prevanish@\place@true\island@false
 #1%
 \place@ b\Bplace@
 \postvanish@}
\def\Bplace@{\allowbreak
 \ifnum\topinscount@=\z@
  \setbox\z@\vbox{\vbox to-\belowtopfigskip{}}%
  \dimen@-\skip\topins
  \ht\z@\dimen@
  \storedim@
  \advancedimtopins@
  \insert\topins{\boxz@}%
  \global\advance\topinscount@\@ne
  \xdef\AAlist@{\AAlist@0}%
 \fi
 \dimen@\ht\islandbox@
 \advance\dimen@\abovebotfigskip
 \ht\islandbox@\dimen@
 \advance\dimen@\dp\islandbox@
 \storedim@
 \xdef\AAlist@{\AAlist@0}%
 \ifnum\topinscount@>\@ne\else\advancedimtopins@\fi
 \insert\topins{\penalty\z@\splittopskip\z@
  \floatingpenalty\z@
  \box\islandbox@}%
 \global\advance\topinscount@\@ne}
\def\breakisland@{\global\setbox\@ne\lastbox\global\skipi@\lastskip\unskip
 \global\setbox\thr@@\lastbox}%
\def\printisland@{\centerline{\box\thr@@}\nobreak\nointerlineskip
 \vskip\skipi@
 \ifdim\ht\@ne<\z@\box\@ne\else\centerline{\box\@ne}\fi}
\def\bottomfigs@{%
 \count@\@ne
 \loop
  \ifnum\count@<\flipcount@
  \nointerlineskip
  \vskip\abovebotfigskip
  \global\setbox\topins\vbox{\unvbox\topins\setboxzl@
   \unvbox\z@
   \breakisland@}%
  \printisland@
  \advance\count@\@ne
  \repeat}
\def\resetdimtopins@{%
 \global\advance\topinscount@-\flipcount@
 \global\setbox\topinsdims@\vbox
  {\unvbox\topinsdims@
   \count@\z@
   \DN@##1##2\next@{\gdef\AAlist@{##2}}%
   \loop
    \ifnum\count@<\flipcount@\setboxzl@
    \expandafter\next@\AAlist@\next@
    \advance\count@\@ne
    \repeat
   \dimen@\z@
   \count@\z@
   \setbox\tw@\vbox{}%
   \edef\nextiii@{\AAlist@}%
   \DN@##1##2\next@{\DNii@{##1}\def\nextiii@{##2}}%
   \loop
    \test@false
    \ifnum\count@<\topinscount@
    \expandafter\next@\nextiii@\next@
     \ifnum\count@<\tw@
      \test@true
     \else
      \if\nextii@ 1\test@true\fi
     \fi
    \fi
    \iftest@
     \setboxzl@
     \advance\dimen@\wdz@
     \setbox\tw@\vbox{\boxz@\unvbox\tw@}%
     \advance\count@\@ne
    \repeat
    \unvbox\tw@
    \global\dimen\topins\dimen@}}
\def\Place@#1#2{%
 \ifisland@
  \ifhmode
   \ifPar@
    \ifnum\Parcount@>5
     \Err@{Only 5 \string\place's allowed per
       \string\Par...\noexpand\endPar paragraph}%
    \else
     \expandafter\expandafter\expandafter\global\expandafter\setbox
      \csname Parfigbox\number\Parcount@\endcsname\box\islandbox@
     \global\advance\Parcount@\@ne
     \xdef\Parlist@{\Parlist@#1}%
     \vadjust{\break}%
    \fi
   \else
    \Err@{\noexpand#2allowed only in a \string\Par...\noexpand\endPar
     paragraph}%
   \fi
  \else
   #2%
  \fi
 \else
  \Err@{Who ... there's no \string\Figure, \string\Table,
   etc., here}%
 \fi
 \place@false}
\newif\ifC@
\newdimen\Cdim@
\long\def\Cplace#1{\prevanish@\place@true\island@false
 #1%
 \Place@ c\Cplace@
 \postvanish@}
\def\Cplace@{\allowbreak
 \ifnum\topinscount@>\z@\else
  \global\C@true\global\Cdim@\pagetotal\fi
 \Aplace@}
\long\def\Mplace#1{\prevanish@\place@true\island@false
 #1%
 \Place@ m\Mplace@
 \postvanish@}
\long\def\MXplace#1{\prevanish@\place@true\island@false
 #1%
 \Place@ M\MXplace@
 \postvanish@}
\newif\ifMX@
\def\MXplace@{\MX@true\Mplace@\MX@false}
\def\Mplace@{\allowbreak
 \dimen@\ht\islandbox@\advance\dimen@\dp\islandbox@
 \ifdim\pagetotal=\z@\else
  \ifdim\lastskip<\abovebotfigskip\advance\dimen@\abovebotfigskip
  \advance\dimen@-\lastskip\fi
 \fi
 \advance\dimen@\pagetotal
 \ifdim\dimen@>\pagegoal
  \Aplace@
 \else
  \nointerlineskip
  \ifdim\lastskip<\abovebotfigskip\removelastskip\vskip\abovebotfigskip\fi
  \setbox\z@\vbox{\unvbox\islandbox@
   \breakisland@}%
  \printisland@
  \ifnum\topinscount@=\z@
   \setbox\z@\vbox{\vbox to-\belowtopfigskip{}}%
   \dimen@-\skip\topins
   \ht\z@\dimen@
   \storedim@
   \advancedimtopins@
   \insert\topins{\boxz@}%
   \global\advance\topinscount@\@ne
   \xdef\AAlist@{\AAlist@0}%
  \fi
  \ifMX@
   \ifnum\topinscount@=\@ne
    \setbox\z@\vbox{\vbox to-\abovebotfigskip{}}%
    \ht\z@\z@
    \dimen@\z@
    \storedim@
    \advancedimtopins@
    \insert\topins{\boxz@}%
    \global\advance\topinscount@\@ne
    \xdef\AAlist@{\AAlist@0}%
   \fi
  \fi
  \nointerlineskip
  \vskip\belowtopfigskip
 \fi}
\expandafter\newbox\csname Parbox1\endcsname
\expandafter\newbox\csname Parbox2\endcsname
\expandafter\newbox\csname Parbox3\endcsname
\expandafter\newbox\csname Parbox4\endcsname
\expandafter\newbox\csname Parbox5\endcsname
\def\endPar{\egroup
 \count@\@ne
 {\vbadness\@M\vfuzz\maxdimen\splitmaxdepth\maxdimen\splittopskip\ht\strutbox
 \setbox\z@\vsplit\Parbox@ to\ht\Parbox@
 \loop
  \ifnum\count@<\Parcount@
  \expandafter\expandafter\expandafter\global\expandafter\setbox
   \csname Parbox\number\count@\endcsname\vsplit\Parbox@ to\ht\Parbox@
  \count@@\count@\advance\count@@\@ne
  \global\csname Parprev\number\count@@\endcsname
   \dp\csname Parbox\number\count@\endcsname
  \advance\count@\@ne
  \repeat}%
 \vskip\parskip
 \count@\@ne
 \def\nextv@##1##2\nextv@{\DN@{##1}\gdef\Parlist@{##2}}%
 \loop
  \ifnum\count@<\Parcount@
   \dimen@\csname Parprev\number\count@\endcsname
   \advance\dimen@\ht\strutbox
   \ifdim\dimen@<\baselineskip
    \advance\dimen@-\baselineskip\vskip-\dimen@
   \else
    \vskip\lineskip
   \fi
   \unvbox\csname Parbox\number\count@\endcsname
   \global\setbox\islandbox@\box\csname Parfigbox\number\count@\endcsname
   \expandafter\nextv@\Parlist@\nextv@
   \if a\next@\Aplace@\else
   \if A\next@\AAplace@\else
   \if b\next@\Bplace@\else
   \if c\next@\Cplace@\else
   \if m\next@\Mplace@\else
   \if M\next@\MXplace@\fi\fi\fi\fi\fi\fi
  \advance\count@\@ne
  \repeat
 \global\Par@false
 \ifvoid\Parbox@
  \prevdepth\csname Parprev\number\count@\endcsname
 \else
  \dimen@\csname Parprev\number\count@\endcsname\advance\dimen@\ht\strutbox
  \ifdim\dimen@<\baselineskip
    \advance\dimen@-\baselineskip\vskip-\dimen@
  \else
    \vskip\lineskip
  \fi
  \dimen@\dp\Parbox@
  \unvbox\Parbox@
  \prevdepth\dimen@
 \fi}
\def\folio{{\page@F\page@S{\page@P\page@N{\number\page@C}\page@Q}}}
\def\advancepageno{\global\advance\pageno\@ne}
\newif\ifspecialsplit@
\newbox\outbox@
\let\shipout@\shipout
\def\plainoutput{\specialsplit@false\ifvoid\topins\else\ifdim\ht\topins=\z@
 \specialsplit@true\advance\minpagesize-\skip\topins\fi\fi
 \fliptopins@
 \setbox\outbox@\vbox{\makeheadline\pagebody\makefootline}%
 {\noexpands@\let\style\relax
 \shipout@\box\outbox@}%
 \advancepageno
 \resetdimtopins@
 \ifvoid\@cclv\else\unvbox\@cclv\penalty\outputpenalty\fi
 \ifnum\outputpenalty>-\@MM\else\dosupereject\fi}
\def\pagebody{\vbox to\vsize{\boxmaxdepth\maxdepth
 \ifvoid\margin@\else
 \rlap{\kern\hsize\vbox to\z@{\kern4\p@\box\margin@\vss}}\fi
 \pagecontents}}
\newif\ifonlytop@
\def\pagecontents{%
 \onlytop@false
 \ifdim\ht\@cclv<\minpagesize\ifnum\flipcount@<\tw@\ifvoid\footins
  \onlytop@true\fi\fi\fi
 \test@false
 \ifC@
  \ifnum\flipcount@=\@ne
   \global\multiply\Cdim@\tw@
   \ifdim\Cdim@>\ht\@cclv
    \test@true
   \fi
  \fi
 \fi
 \global\C@false
 \iftest@
  \dimen@\ht\@cclv
  \advance\dimen@\skip\topins
  {\vfuzz\maxdimen\vbadness\@M
  \splitmaxdepth\maxdepth\splittopskip\topskip
  \setbox\z@\vsplit\@cclv to\dimen@
  \unvbox\z@}%
  \global\setbox\topins\vbox{\unvbox\topins
   \global\setbox\@ne\lastbox}%
  \setbox\z@\vbox{\unvbox\@ne
   \breakisland@}%
  \nointerlineskip
  \vskip\abovebotfigskip
  \printisland@
 \else
  \ifnum\flipcount@>\z@
   \global\setbox\topins\vbox{\unvbox\topins\global\setbox\@ne\lastbox}%
   \setbox\z@\vbox{\unvbox\@ne
    \breakisland@}%
   \printisland@
   \ifonlytop@\kern-\prevdepth\vfill\else\vskip\belowtopfigskip\fi
  \fi
 \fi
 \ifdim\ht\@cclv<\minpagesize
  \ifonlytop@\else\vfill\fi
 \else
  \ifspecialsplit@
   {\vfuzz\maxdimen\vbadness\@M
   \splitmaxdepth\maxdepth\splittopskip\topskip
   \dimen@ii\ht\@cclv \advance\dimen@ii\skip\topins
   \setbox\z@\vsplit\@cclv to\dimen@ii
   \unvbox\z@}%
  \else
   \unvbox\@cclv
  \fi
 \fi
 \bottomfigs@
 \ifvoid\footins\else\vskip\skip\footins\footnoterule\unvbox\footins\fi}
\newread\readdata@
\def\readthedata@#1{\expandafter
 \ifx\csname#1@D\endcsname\relax
  \immediate\openin\readdata@=#1.dat
  \ifeof\readdata@
   \Err@{No file #1.dat}%
  \else
   {\endlinechar\m@ne\gdef\Next@{}%
   \DNii@##1 ##2 ##3pt{\global\data@ht##1\global\data@dp##2%
    \global\data@wd##3pt}%
   \loop
    \ifeof\readdata@
    \else
    \read\readdata@ to\next@
    \ifx\next@\empty\else
     \edef\next@{\expandafter\nextii@\next@}%
     \expandafter\rightadd@\next@\to\Next@
    \fi
    \repeat}%
   \immediate\closein\readdata@
   \expandafter\expandafter\expandafter\global\expandafter
    \let\csname#1@D\endcsname\Next@\global\let\Next@\relax
  \fi
 \fi}
\newdimen\data@ht
\newdimen\data@dp
\newdimen\data@wd
\newif\ifgetdata@
\def\getdata@#1#2{\global\getdata@true\count@#2\relax
 {\let\\\or\xdef\Next@{\ifcase\number\count@#1\else
 \global\noexpand\getdata@false\fi}}\Next@}
\def\paste#1#2{\readthedata@{#1}%
 \getdata@{\csname#1@D\endcsname}{#2}%
 \ifgetdata@
 \dimen@\data@ht \advance\dimen@\data@dp
  \hbox{\special{dvipaste: #1 #2}%
   \lower\data@dp\vbox to\dimen@{\hbox to\data@wd{}\vfil}}%
 \else
  {\lccode`\Z=`\#\lccode`\N=`\N\lccode`\F=`\F%
   \lowercase{\Err@{No data for File [#1], Z#2}}}%
 \fi}
\newdimen\httable
\newdimen\dptable
\newdimen\wdtable
\def\measuretable#1#2{\readthedata@{#1}%
 \getdata@{\csname#1@D\endcsname}{#2}%
 \ifgetdata@
  \httable\data@ht \dptable\data@dp \wdtable\data@wd
 \else
  {\lccode`\Z=`\#\lccode`\N=`\N\lccode`\F=`\F%
  \lowercase{\Err@{No data for File [#1], Z#2}}}%
 \fi}
\def\East#1#2{\setboxz@h{$\m@th\ssize\;{#1}\;\;$}%
 \setbox\tw@\hbox{$\m@th\ssize\;{#2}\;\;$}\setbox4=\hbox{$\m@th#2$}%
 \dimen@\minaw@
 \ifdim\wdz@>\dimen@\dimen@\wdz@\fi\ifdim\wd\tw@>\dimen@\dimen@\wd\tw@\fi
 \ifdim\wd4 >\z@
  \mathrel{\mathop{\hbox to\dimen@{\rightarrowfill}}\limits^{#1}_{#2}}%
 \else
  \mathrel{\mathop{\hbox to\dimen@{\rightarrowfill}}\limits^{#1}}%
 \fi}
\def\West#1#2{\setboxz@h{$\m@th\ssize\;\;{#1}\;$}%
 \setbox\tw@\hbox{$\m@th\ssize\;\;{#2}\;$}\setbox4=\hbox{$\m@th#2$}%
 \dimen@\minaw@
 \ifdim\wdz@>\dimen@\dimen@\wdz@\fi\ifdim\wd\tw@>\dimen@\dimen@\wd\tw@\fi
 \ifdim\wd4 >\z@
  \mathrel{\mathop{\hbox to\dimen@{\leftarrowfill}}\limits^{#1}_{#2}}%
 \else
  \mathrel{\mathop{\hbox to\dimen@{\leftarrowfill}}\limits^{#1}}%
 \fi}
\font\arrow@i=lams1
\font\arrow@ii=lams2
\font\arrow@iii=lams3
\font\arrow@iv=lams4
\font\arrow@v=lams5
\newdimen\standardcgap
\standardcgap40\p@
\newdimen\hunit
\hunit\tw@\p@
\newdimen\standardrgap
\standardrgap32\p@
\newdimen\vunit
\vunit1.6\p@
\def\Cgaps#1{\RIfM@
 \standardcgap#1\standardcgap\relax\hunit#1\hunit\relax
 \else\nonmatherr@\Cgaps\fi}
\def\Rgaps#1{\RIfM@
 \standardrgap#1\standardrgap\relax\vunit#1\vunit\relax
 \else\nonmatherr@\Rgaps\fi}
\newdimen\getdim@
\def\getcgap@#1{\ifcase#1\or\getdim@\z@\else\getdim@\standardcgap\fi}
\def\getrgap@#1{\ifcase#1\getdim@\z@\else\getdim@\standardrgap\fi}
\def\cgaps{\RIfM@\expandafter\cgaps@\else\expandafter\nonmatherr@
 \expandafter\cgaps\fi}
\def\cgaps@{\ifnum\catcode`\;=\active\expandafter\cgapsA@\else
 \expandafter\cgapsO@\fi}
\def\cgapsO@#1{\toks@{\ifcase\i@\or\getdim@=\z@}%
 \gaps@@\standardcgap#1;\gaps@@\gaps@@
 \edef\next@{\the\toks@\noexpand\else\noexpand\getdim@\noexpand\standardcgap
  \noexpand\fi}%
 \toks@=\expandafter{\next@}%
 \edef\getcgap@##1{\i@##1\relax\the\toks@}\toks@{}}
{\catcode`\;=\active
 \gdef\cgapsA@#1{\toks@{\ifcase\i@\or\getdim@=\z@}%
 \gaps@@\standardcgap#1;\gaps@@\gaps@@
 \edef\next@{\the\toks@\noexpand\else\noexpand\getdim@\noexpand\standardcgap
  \noexpand\fi}%
 \toks@=\expandafter{\next@}%
 \edef\getcgap@##1{\i@##1\relax\the\toks@}\toks@{}}
}
\def\Gaps@@{\gaps@@}
\def\gaps@@#1#2;#3{\mgaps@#1#2\mgaps@
 \edef\next@{\the\toks@\noexpand\or\noexpand\getdim@
  \noexpand#1\the\mgapstoks@@}%
 \toks@\expandafter{\next@}%
 \DN@{#3}%
 \ifx\next@\Gaps@@\def\next@##1\gaps@@{}\else
  \def\next@{\gaps@@#1#3}\fi\next@}
{\catcode`\;=\active
 \gdef\rgaps#1{\RIfM@{\ifnum\catcode`\;=\active\def;{\string;}\fi
   \xdef\Next@{\noexpand\rgaps@{#1}}}%
  \Next@\edef\getrgap@##1{\i@##1\relax\the\toks@}\toks@{}\else
  \nonmatherr@\rgaps\fi}
}
\def\rgaps@#1{\toks@{\ifcase\i@\getdim@=\z@}%
 \gaps@@\standardrgap#1;\gaps@@\gaps@@
 \edef\next@{\the\toks@\noexpand\else\noexpand\getdim@\noexpand\standardrgap
  \noexpand\fi}%
 \toks@=\expandafter{\next@}}
\newbox\ZER@
\def\mgaps@#1{\let\mgapsnext@#1\FNSS@\mgaps@@}
\def\mgaps@@{\ifx\next\w\expandafter\mgaps@@@\else
 \expandafter\mgaps@@@@\fi}
\newtoks\mgapstoks@@
\def\mgaps@@@@#1\mgaps@{\getdim@\mgapsnext@\getdim@#1\getdim@
 \edef\next@{\noexpand\getdim@\the\getdim@}%
 \mgapstoks@@\expandafter{\next@}}
\def\mgaps@@@\w#1#2\mgaps@{\mgaps@@@@#2\mgaps@
 \setbox\ZER@\hbox{$\m@th\hskip15\p@\tsize@#1$}%
 \dimen@\wd\ZER@
 \ifdim\dimen@>\getdim@\getdim@\dimen@\fi
 \edef\next@{\noexpand\getdim@\the\getdim@}%
 \mgapstoks@@\expandafter{\next@}}
\def\changewidth#1#2{\setbox\ZER@{$\m@th#2}%
 \hbox to\wd\ZER@{\hss$\m@th#1$\hss}}
\atdef@({\FN@\ARROW@}
\def\ARROW@{\ifx\next)\let\next@\OPTIONS@\else
 \DN@{\csname\string @(\endcsname}\fi\next@}
\newif\ifoptions@
\def\OPTIONS@){\ifoptions@\let\next@\relax\else
 \DN@{\global\options@true\begingroup\optioncodes@}\fi\next@}
\newif\ifN@
\newif\ifE@
\newif\ifNESW@
\newif\ifH@
\newif\ifV@
\newif\ifHshort@
\expandafter\def\csname\string @(\endcsname #1,#2){%
 \ifoptions@\expandafter\endgroup\fi
 \N@false\E@false\H@false\V@false\Hshort@false
 \ifnum#1>\z@\E@true\fi
 \ifnum#1=\z@\V@true\global\tX@false\global\tY@false\global\a@false\fi
 \ifnum#2>\z@\N@true\fi
 \ifnum#2=\z@\H@true\global\tX@false\global\tY@false\global\a@false
  \ifshort@\Hshort@true\fi\fi
 \NESW@false
 \ifN@\ifE@\NESW@true\fi\else\ifE@\else\NESW@true\fi\fi
 \arrow@{#1}{#2}%
 \global\options@false
 \global\scount@\z@\global\tcount@\z@\global\arrcount@\z@
 \global\s@false\global\sxdimen@\z@\global\sydimen@\z@
 \global\tX@false\global\tXdimen@i\z@\global\tXdimen@ii\z@
 \global\tY@false\global\tYdimen@i\z@\global\tYdimen@ii\z@
 \global\a@false\global\exacount@\z@
 \global\x@false\global\xdimen@\z@
 \global\X@false\global\Xdimen@\z@
 \global\y@false\global\ydimen@\z@
 \global\Y@false\global\Ydimen@\z@
 \global\p@false\global\pdimen@\z@
 \global\label@ifalse\global\label@iifalse
 \global\dl@ifalse\global\ldimen@i\z@
 \global\dl@iifalse\global\ldimen@ii\z@
 \global\short@false\global\unshort@false}
\newif\iflabel@i
\newif\iflabel@ii
\newcount\scount@
\newcount\tcount@
\newcount\arrcount@
\newif\ifs@
\newdimen\sxdimen@
\newdimen\sydimen@
\newif\iftX@
\newdimen\tXdimen@i
\newdimen\tXdimen@ii
\newif\iftY@
\newdimen\tYdimen@i
\newdimen\tYdimen@ii
\newif\ifa@
\newcount\exacount@
\newif\ifx@
\newdimen\xdimen@
\newif\ifX@
\newdimen\Xdimen@
\newif\ify@
\newdimen\ydimen@
\newif\ifY@
\newdimen\Ydimen@
\newif\ifp@
\newdimen\pdimen@
\newif\ifdl@i
\newif\ifdl@ii
\newdimen\ldimen@i
\newdimen\ldimen@ii
\newif\ifshort@
\newif\ifunshort@
\def\zero@#1{\ifnum\scount@=\z@
 \if#1e\global\scount@\m@ne\else
 \if#1t\global\scount@\tw@\else
 \if#1h\global\scount@\thr@@\else
 \if#1'\global\scount@6 \else
 \if#1`\global\scount@7 \else
 \if#1(\global\scount@8 \else
 \if#1)\global\scount@9 \else
 \if#1s\global\scount@12 \else
 \if#1H\global\scount@13 \else
 \Err@{\Invalid@@ option \string\0}\fi\fi\fi\fi\fi\fi\fi\fi\fi
 \fi}
\def\one@#1{\ifnum\tcount@=\z@
 \if#1e\global\tcount@\m@ne\else
 \if#1h\global\tcount@\tw@\else
 \if#1t\global\tcount@\thr@@\else
 \if#1'\global\tcount@4 \else
 \if#1`\global\tcount@5 \else
 \if#1(\global\tcount@\ten@ \else
 \if#1)\global\tcount@11 \else
 \if#1s\global\tcount@12 \else
 \if#1H\global\tcount@13 \else
 \Err@{\Invalid@@ option \string\1}\fi\fi\fi\fi\fi\fi\fi\fi\fi
 \fi}
\def\a@#1{\ifnum\arrcount@=\z@
 \if#10\global\arrcount@\m@ne\else
 \if#1+\global\arrcount@\@ne\else
 \if#1-\global\arrcount@\tw@\else
 \if#1=\global\arrcount@\thr@@\else
 \Err@{\Invalid@@ option \string\a}\fi\fi\fi\fi
 \fi}
\def\ds@{\ifnum\catcode`\;=\active\expandafter\dsA@\else
 \expandafter\dsO@\fi}
\def\dsO@(#1;#2){\ds@@{#1}{#2}}
\def\ds@@#1#2{\ifs@\else
 \global\s@true
 \global\sxdimen@\hunit\global\sxdimen@#1\sxdimen@\relax
 \global\sydimen@\vunit\global\sydimen@#2\sydimen@\relax
 \fi}
\def\dtX@{\ifnum\catcode`\;=\active\expandafter\dtXA@\else
 \expandafter\dtXO@\fi}
\def\dtXO@(#1;#2){\dtX@@{#1}{#2}}
\def\dtX@@#1#2{\iftX@\else
 \global\tX@true
 \global\tXdimen@i\hunit\global\tXdimen@i#1\tXdimen@i\relax
 \global\tXdimen@ii\vunit\global\tXdimen@ii#2\tXdimen@ii\relax
 \fi}
\def\dtY@{\ifnum\catcode`\;=\active\expandafter\dtYA@\else
 \expandafter\dtYO@\fi}
\def\dtYO@(#1;#2){\dtY@@{#1}{#2}}
\def\dtY@@#1#2{\iftY@\else
 \global\tY@true
 \global\tYdimen@i\hunit\global\tYdimen@i#1\tYdimen@i\relax
 \global\tYdimen@ii\vunit\global\tYdimen@ii#2\tYdimen@ii\relax
 \fi}
{\catcode`\;=\active
 \gdef\dsA@(#1;#2){\ds@@{#1}{#2}}
 \gdef\dtXA@(#1;#2){\dtX@@{#1}{#2}}
 \gdef\dtYA@(#1;#2){\dtY@@{#1}{#2}}
}
\def\da@#1{\ifa@\else\global\a@true\global\exacount@#1\relax\fi}
\def\dx@#1{\ifx@\else
 \global\x@true
 \global\xdimen@\hunit\global\xdimen@#1\xdimen@\relax
 \fi}
\def\dX@#1{\ifX@\else
 \global\X@true
 \global\Xdimen@\hunit\global\Xdimen@#1\Xdimen@\relax
 \fi}
\def\dy@#1{\ify@\else
 \global\y@true
 \global\ydimen@\vunit\global\ydimen@#1\ydimen@\relax
 \fi}
\def\dY@#1{\ifY@\else
 \global\Y@true
 \global\Ydimen@\vunit\global\Ydimen@#1\Ydimen@\relax
 \fi}
\def\p@@#1{\ifp@\else
 \global\p@true
 \global\pdimen@\hunit\global\divide\pdimen@\tw@
 \global\pdimen@#1\pdimen@\relax
 \fi}
\def\L@#1{\iflabel@i\else
 \global\label@itrue\gdef\label@i{#1}%
 \fi}
\def\l@#1{\iflabel@ii\else
 \global\label@iitrue\gdef\label@ii{#1}%
 \fi}
\def\dL@#1{\ifdl@i\else
 \global\dl@itrue\global\ldimen@i\hunit\global\ldimen@i#1\ldimen@i\relax
 \fi}
\def\dl@#1{\ifdl@ii\else
 \global\dl@iitrue\global\ldimen@ii\hunit\global\ldimen@ii#1\ldimen@ii\relax
 \fi}
\def\s@{\ifunshort@\else\global\short@true\fi}
\def\uns@{\ifshort@\else\global\unshort@true\global\short@false\fi}
\def\optioncodes@{\let\0\zero@\let\1\one@\let\a\a@\let\ds\ds@\let\dtX\dtX@
 \let\dtY\dtY@\let\da\da@\let\dx\dx@\let\dX\dX@\let\dY\dY@\let\dy\dy@
 \let\p\p@@\let\L\L@\let\l\l@\let\dL\dL@\let\dl\dl@\let\s\s@\let\uns\uns@}
\def\slopes@{\\161\\152\\143\\134\\255\\126\\357\\238\\349\\45{10}\\56{11}%
 \\11{12}\\65{13}\\54{14}\\43{15}\\32{16}\\53{17}\\21{18}\\52{19}\\31{20}%
 \\41{21}\\51{22}\\61{23}}
\newcount\tan@i
\newcount\tan@ip
\newcount\tan@ii
\newcount\tan@iip
\newdimen\slope@i
\newdimen\slope@ip
\newdimen\slope@ii
\newdimen\slope@iip
\newcount\angcount@
\newcount\extracount@
\def\slope@{{\slope@i\secondy@\advance\slope@i-\firsty@
 \ifN@\else\multiply\slope@i\m@ne\fi
 \slope@ii\secondx@\advance\slope@ii-\firstx@
 \ifE@\else\multiply\slope@ii\m@ne\fi
 \ifdim\slope@ii<\z@
  \global\tan@i6 \global\tan@ii\@ne\global\angcount@23
 \else
  \dimen@\slope@i\multiply\dimen@6
  \ifdim\dimen@<\slope@ii
   \global\tan@i\@ne\global\tan@ii6 \global\angcount@\@ne
  \else
   \dimen@\slope@ii\multiply\dimen@6
   \ifdim\dimen@<\slope@i
    \global\tan@i6 \global\tan@ii\@ne\global\angcount@23
   \else
    \global\tan@ip\z@\global\tan@iip\@ne
    \def\\##1##2##3{\global\angcount@##3\relax
     \slope@ip\slope@i\slope@iip\slope@ii
     \multiply\slope@iip##1\relax\multiply\slope@ip##2\relax
     \ifdim\slope@iip<\slope@ip
      \global\tan@ip##1\relax\global\tan@iip##2\relax
     \else
      \global\tan@i##1\relax\global\tan@ii##2\relax
      \def\\####1####2####3{}%
     \fi}%
    \slopes@
    \slope@i\secondy@\advance\slope@i-\firsty@
    \ifN@\else\multiply\slope@i\m@ne\fi
    \multiply\slope@i\tan@ii\multiply\slope@i\tan@iip\multiply\slope@i\tw@
    \count@\tan@i\multiply\count@\tan@iip
    \extracount@\tan@ip\multiply\extracount@\tan@ii
    \advance\count@\extracount@
    \slope@ii\secondx@\advance\slope@ii-\firstx@
    \ifE@\else\multiply\slope@ii\m@ne\fi
    \multiply\slope@ii\count@
    \ifdim\slope@i<\slope@ii
     \global\tan@i\tan@ip\global\tan@ii\tan@iip
     \global\advance\angcount@\m@ne
    \fi
   \fi
  \fi
 \fi}%
}
\def\slope@a#1{{\def\\##1##2##3{\ifnum##3=#1\global\tan@i##1\relax
 \global\tan@ii##2\relax\fi}\slopes@}}
\newcount\i@
\newcount\j@
\newcount\colcount@
\newcount\Colcount@
\newcount\tcolcount@
\newdimen\rowht@
\newdimen\rowdp@
\newcount\rowcount@
\newcount\Rowcount@
\newcount\maxcolrow@
\newtoks\colwidthtoks@
\newtoks\Rowheighttoks@
\newtoks\Rowdepthtoks@
\newtoks\widthtoks@
\newtoks\Widthtoks@
\newtoks\heighttoks@
\newtoks\Heighttoks@
\newtoks\depthtoks@
\newtoks\Depthtoks@
\newif\iffirstCDcr@
\def\dotoks@i{%
 \global\widthtoks@\expandafter{\the\widthtoks@\else\getdim@\z@\fi}%
 \global\heighttoks@\expandafter{\the\heighttoks@\else\getdim@\z@\fi}%
 \global\depthtoks@\expandafter{\the\depthtoks@\else\getdim@\z@\fi}}
\def\dotoks@ii{%
 \global\widthtoks@{\ifcase\j@}%
 \global\heighttoks@{\ifcase\j@}%
 \global\depthtoks@{\ifcase\j@}}
\def\preCD@#1\endCD{\setbox\ZER@
 \vbox{%
  \def\arrow@##1##2{{}}%
  \global\rowcount@\m@ne\global\colcount@\z@\global\Colcount@\z@
  \global\firstCDcr@true\toks@{}%
  \global\widthtoks@{\ifcase\j@}%
  \global\Widthtoks@{\ifcase\i@}%
  \global\heighttoks@{\ifcase\j@}%
  \global\Heighttoks@{\ifcase\i@}%
  \global\depthtoks@{\ifcase\j@}%
  \global\Depthtoks@{\ifcase\i@}%
  \global\Rowheighttoks@{\ifcase\i@}%
  \global\Rowdepthtoks@{\ifcase\i@}%
  \Let@
  \everycr{%
   \noalign{%
    \global\advance\rowcount@\@ne
    \ifnum\colcount@<\Colcount@
    \else
     \global\Colcount@\colcount@\global\maxcolrow@\rowcount@
    \fi
    \global\colcount@\z@
    \iffirstCDcr@
     \global\firstCDcr@false
    \else
     \edef\next@{\the\Rowheighttoks@\noexpand\or\noexpand\getdim@\the\rowht@}%
      \global\Rowheighttoks@\expandafter{\next@}%
     \edef\next@{\the\Rowdepthtoks@\noexpand\or\noexpand\getdim@\the\rowdp@}%
      \global\Rowdepthtoks@\expandafter{\next@}%
     \global\rowht@\z@\global\rowdp@\z@
     \dotoks@i
     \edef\next@{\the\Widthtoks@\noexpand\or\the\widthtoks@}%
      \global\Widthtoks@\expandafter{\next@}%
     \edef\next@{\the\Heighttoks@\noexpand\or\the\heighttoks@}%
      \global\Heighttoks@\expandafter{\next@}%
     \edef\next@{\the\Depthtoks@\noexpand\or\the\depthtoks@}%
      \global\Depthtoks@\expandafter{\next@}%
     \dotoks@ii
    \fi}}%
  \tabskip\z@
  \halign{&\setbox\ZER@\hbox{\vrule\height\ten@\p@\width\z@\depth\z@     
   $\m@th\displaystyle{##}$}\copy\ZER@
   \ifdim\ht\ZER@>\rowht@\global\rowht@\ht\ZER@\fi
   \ifdim\dp\ZER@>\rowdp@\global\rowdp@\dp\ZER@\fi
   \global\advance\colcount@\@ne
   \edef\next@{\the\widthtoks@\noexpand\or\noexpand\getdim@\the\wd\ZER@}%
    \global\widthtoks@\expandafter{\next@}%
   \edef\next@{\the\heighttoks@\noexpand\or\noexpand\getdim@\the\ht\ZER@}%
    \global\heighttoks@\expandafter{\next@}%
   \edef\next@{\the\depthtoks@\noexpand\or\noexpand\getdim@\the\dp\ZER@}%
    \global\depthtoks@\expandafter{\next@}%
   \cr#1\crcr}}%
 \Rowcount@\rowcount@
 \global\Widthtoks@\expandafter{\the\Widthtoks@\fi\relax}%
 \edef\Width@##1##2{\i@##1\relax\j@##2\relax\the\Widthtoks@}%
 \global\Heighttoks@\expandafter{\the\Heighttoks@\fi\relax}%
 \edef\Height@##1##2{\i@##1\relax\j@##2\relax\the\Heighttoks@}%
 \global\Depthtoks@\expandafter{\the\Depthtoks@\fi\relax}%
 \edef\Depth@##1##2{\i@##1\relax\j@##2\relax\the\Depthtoks@}%
 \edef\next@{\the\Rowheighttoks@\noexpand\fi\relax}%
 \global\Rowheighttoks@\expandafter{\next@}%
 \edef\Rowheight@##1{\i@##1\relax\the\Rowheighttoks@}%
 \edef\next@{\the\Rowdepthtoks@\noexpand\fi\relax}%
 \global\Rowdepthtoks@\expandafter{\next@}%
 \edef\Rowdepth@##1{\i@##1\relax\the\Rowdepthtoks@}%
 \global\colwidthtoks@{\fi}%
 \setbox\ZER@\vbox{%
  \unvbox\ZER@
  \count@\rowcount@
  \loop
   \unskip\unpenalty
   \setbox\ZER@\lastbox
   \ifnum\count@>\maxcolrow@\advance\count@\m@ne
   \repeat
  \hbox{%
   \unhbox\ZER@
   \count@\z@
   \loop
    \unskip
    \setbox\ZER@\lastbox
    \edef\next@{\noexpand\or\noexpand\getdim@\the\wd\ZER@\the\colwidthtoks@}%
     \global\colwidthtoks@\expandafter{\next@}%
    \advance\count@\@ne
    \ifnum\count@<\Colcount@
    \repeat}}%
 \edef\next@{\noexpand\ifcase\noexpand\i@\the\colwidthtoks@}%
  \global\colwidthtoks@\expandafter{\next@}%
 \edef\Colwidth@##1{\i@##1\relax\the\colwidthtoks@}%
 \global\colwidthtoks@{}\global\Rowheighttoks@{}\global\Rowdepthtoks@{}%
 \global\widthtoks@{}\global\Widthtoks@{}\global\heighttoks@{}%
 \global\Heighttoks@{}\global\depthtoks@{}\global\Depthtoks@{}%
}
\newcount\xoff@
\newcount\yoff@
\newcount\endcount@
\newcount\rcount@
\newdimen\firstx@
\newdimen\firsty@
\newdimen\secondx@
\newdimen\secondy@
\newdimen\tocenter@
\newdimen\charht@
\newdimen\charwd@
\def\outside@{\Err@{This arrow points outside the \string\CD}}
\newif\ifsvertex@
\newif\iftvertex@
\def\arrow@#1#2{\global\xoff@#1\relax\global\yoff@#2\relax
 \count@\rowcount@\advance\count@-\yoff@
 \ifnum\count@<\@ne\outside@\else\ifnum\count@>\Rowcount@\outside@\fi\fi
 \count@\colcount@\advance\count@\xoff@
 \ifnum\count@<\@ne\outside@\else\ifnum\count@>\Colcount@\outside@\fi\fi
 \tcolcount@\colcount@\advance\tcolcount@\xoff@
 \Width@\rowcount@\colcount@\divide\getdim@\tw@\tocenter@-\getdim@
 \ifdim\getdim@=\z@
  \firstx@\z@\firsty@\mathaxis@\svertex@true
 \else
  \svertex@false
  \ifHshort@
   \Colwidth@\colcount@\divide\getdim@\tw@
   \ifE@ \firstx@\getdim@ \else \firstx@-\getdim@ \fi
  \else
   \ifE@ \firstx@\getdim@ \else \firstx@-\getdim@ \fi
  \fi
  \ifE@
   \ifH@ \advance\firstx@\thr@@\p@ \else \advance\firstx@-\thr@@\p@ \fi  
  \else
   \ifH@ \advance\firstx@-\thr@@\p@ \else \advance\firstx@\thr@@\p@ \fi  
  \fi
  \ifN@
   \Height@\rowcount@\colcount@ \firsty@\getdim@                         
   \ifV@ \advance\firsty@\thr@@\p@ \fi                                   
  \else
   \ifV@
    \Depth@\rowcount@\colcount@ \firsty@-\getdim@                        
    \advance\firsty@-\thr@@\p@                                           
   \else
    \firsty@\z@                                                          
   \fi
  \fi
 \fi
 \ifV@
 \else
  \Colwidth@\colcount@\divide\getdim@\tw@
  \ifE@\secondx@\getdim@\else\secondx@-\getdim@\fi
  \ifE@\else\getcgap@\colcount@\advance\secondx@-\getdim@\fi
  \endcount@\colcount@\advance\endcount@\xoff@
  \count@\colcount@
  \ifE@
   \advance\count@\@ne
   \loop
    \ifnum\count@<\endcount@
    \Colwidth@\count@\advance\secondx@\getdim@
    \getcgap@\count@\advance\secondx@\getdim@
    \advance\count@\@ne
    \repeat
  \else
   \advance\count@\m@ne
   \loop
    \ifnum\count@>\endcount@
    \Colwidth@\count@\advance\secondx@-\getdim@
    \getcgap@\count@\advance\secondx@-\getdim@
    \advance\count@\m@ne
    \repeat
  \fi
  \Colwidth@\count@\divide\getdim@\tw@
  \ifHshort@
  \else
   \ifE@\advance\secondx@\getdim@\else\advance\secondx@-\getdim@\fi
  \fi
  \ifE@\getcgap@\count@\advance\secondx@\getdim@\fi
  \rcount@\rowcount@\advance\rcount@-\yoff@
  \Width@\rcount@\count@\divide\getdim@\tw@
  \tvertex@false
  \ifH@\ifdim\getdim@=\z@\tvertex@true\Hshort@false\fi\fi
  \ifHshort@
  \else
   \ifE@\advance\secondx@-\getdim@\else\advance\secondx@\getdim@\fi
  \fi
  \iftvertex@
   \advance\secondx@.4\p@
  \else
   \ifE@\advance\secondx@-\thr@@\p@\else\advance\secondx@\thr@@\p@\fi    
  \fi
 \fi
 \ifH@
 \else
  \ifN@
   \Rowheight@\rowcount@\secondy@\getdim@
  \else
   \Rowdepth@\rowcount@\secondy@-\getdim@
   \getrgap@\rowcount@\advance\secondy@-\getdim@
  \fi
  \endcount@\rowcount@\advance\endcount@-\yoff@
  \count@\rowcount@
  \ifN@
   \advance\count@\m@ne
   \loop
    \ifnum\count@>\endcount@
    \Rowheight@\count@\advance\secondy@\getdim@
    \Rowdepth@\count@\advance\secondy@\getdim@
    \getrgap@\count@\advance\secondy@\getdim@
    \advance\count@\m@ne
    \repeat
  \else
   \advance\count@\@ne
   \loop
    \ifnum\count@<\endcount@
    \Rowheight@\count@\advance\secondy@-\getdim@
    \Rowdepth@\count@\advance\secondy@-\getdim@
    \getrgap@\count@\advance\secondy@-\getdim@
    \advance\count@\@ne
    \repeat
  \fi
  \tvertex@false
  \ifV@\Width@\count@\colcount@\ifdim\getdim@=\z@\tvertex@true\fi\fi
  \ifN@
   \getrgap@\count@\advance\secondy@\getdim@
   \Rowdepth@\count@\advance\secondy@\getdim@
   \iftvertex@
    \advance\secondy@\mathaxis@
   \else
    \Depth@\count@\tcolcount@\advance\secondy@-\getdim@
    \advance\secondy@-\thr@@\p@                                          
   \fi
  \else
   \Rowheight@\count@\advance\secondy@-\getdim@
   \iftvertex@
    \advance\secondy@\mathaxis@
   \else
    \Height@\count@\tcolcount@\advance\secondy@\getdim@
    \advance\secondy@\thr@@\p@                                           
   \fi
  \fi
 \fi
 \ifV@\else\advance\firstx@\sxdimen@\fi
 \ifH@\else\advance\firsty@\sydimen@\fi
 \iftX@
  \advance\secondy@\tXdimen@ii
  \advance\secondx@\tXdimen@i
  \slope@
 \else
  \iftY@
   \advance\secondy@\tYdimen@ii
   \advance\secondx@\tYdimen@i
   \slope@
   \secondy@\secondx@\advance\secondy@-\firstx@
   \ifNESW@\else\multiply\secondy@\m@ne\fi
   \multiply\secondy@\tan@i\divide\secondy@\tan@ii\advance\secondy@\firsty@
  \else
   \ifa@
    \slope@
    \ifNESW@\global\advance\angcount@\exacount@\else
     \global\advance\angcount@-\exacount@\fi
    \ifnum\angcount@>23 \global\angcount@23 \fi
    \ifnum\angcount@<\@ne\global\angcount@\@ne\fi
    \slope@a\angcount@
    \ifY@
     \advance\secondy@\Ydimen@
    \else
     \ifX@
      \advance\secondx@\Xdimen@
      \dimen@\secondx@\advance\dimen@-\firstx@
      \ifNESW@\else\multiply\dimen@\m@ne\fi
      \multiply\dimen@\tan@i\divide\dimen@\tan@ii
      \advance\dimen@\firsty@\secondy@\dimen@
     \fi
    \fi
   \else
    \ifH@\else\ifV@\else\slope@\fi\fi
   \fi
  \fi
 \fi
 \ifH@\else\ifV@\else\ifsvertex@\else
  \dimen@6\p@\multiply\dimen@\tan@ii
  \count@\tan@i\advance\count@\tan@ii\divide\dimen@\count@
  \ifE@\advance\firstx@\dimen@\else\advance\firstx@-\dimen@\fi
  \multiply\dimen@\tan@i\divide\dimen@\tan@ii
  \ifN@\advance\firsty@\dimen@\else\advance\firsty@-\dimen@\fi
 \fi\fi\fi
 \ifp@
  \ifH@\else\ifV@\else
   \getcos@\pdimen@\advance\firsty@\dimen@\advance\secondy@\dimen@
   \ifNESW@\advance\firstx@-\dimen@ii\else\advance\firstx@\dimen@ii\fi
  \fi\fi
 \fi
 \ifH@\else\ifV@\else
  \ifnum\tan@i>\tan@ii
   \charht@\ten@\p@\charwd@\ten@\p@
   \multiply\charwd@\tan@ii\divide\charwd@\tan@i
  \else
   \charwd@\ten@\p@\charht@\ten@\p@
   \divide\charht@\tan@ii\multiply\charht@\tan@i
  \fi
  \ifnum\tcount@=\thr@@
   \ifN@\advance\secondy@-.3\charht@\else\advance\secondy@.3\charht@\fi
  \fi
  \ifnum\scount@=\tw@
   \ifE@\advance\firstx@.3\charht@\else\advance\firstx@-.3\charht@\fi
  \fi
  \ifnum\tcount@=12
   \ifN@\advance\secondy@-\charht@\else\advance\secondy@\charht@\fi
  \fi
  \iftY@
  \else
   \ifa@
    \ifX@
    \else
     \secondx@\secondy@\advance\secondx@-\firsty@
     \ifNESW@\else\multiply\secondx@\m@ne\fi
     \multiply\secondx@\tan@ii\divide\secondx@\tan@i
     \advance\secondx@\firstx@
    \fi
   \fi
  \fi
 \fi\fi
 \ifH@\harrow@\else\ifV@\varrow@\else\arrow@@\fi\fi}
\newdimen\mathaxis@
\mathaxis@90\p@\divide\mathaxis@36
\def\harrow@b{\ifE@\hskip\tocenter@\hskip\firstx@\fi}
\def\harrow@bb{\ifE@\hskip\xdimen@\else\hskip\Xdimen@\fi}
\def\harrow@e{\ifE@\else\hskip-\firstx@\hskip-\tocenter@\fi}
\def\harrow@ee{\ifE@\hskip-\Xdimen@\else\hskip-\xdimen@\fi}
\def\harrow@{\dimen@\secondx@\advance\dimen@-\firstx@
 \ifE@\let\next@\rlap\else\multiply\dimen@\m@ne\let\next@\llap\fi
 \next@{%
  \harrow@b
  \smash{\raise\pdimen@\hbox to\dimen@
   {\harrow@bb\arrow@ii
    \ifnum\arrcount@=\m@ne\else\ifnum\arrcount@=\thr@@\else
     \ifE@
      \ifnum\scount@=\m@ne
      \else
       \ifcase\scount@\or\or\char118 \or\char117 \or\or\or\char119 \or
       \char120 \or\char121 \or\char122 \or\or\or\arrow@i\char125 \or
       \char117 \hskip\thr@@\p@\char117 \hskip-\thr@@\p@\fi
      \fi
     \else
      \ifnum\tcount@=\m@ne
      \else
       \ifcase\tcount@\char117 \or\or\char117 \or\char118 \or\char119 \or
       \char120 \or\or\or\or\or\char121 \or\char122 \or\arrow@i\char125
       \or\char117 \hskip\thr@@\p@\char117 \hskip-\thr@@\p@\fi
      \fi
     \fi
    \fi\fi
    \dimen@\mathaxis@\advance\dimen@.2\p@
    \dimen@ii\mathaxis@\advance\dimen@ii-.2\p@
    \ifnum\arrcount@=\m@ne
     \let\leads@\null
    \else
     \ifcase\arrcount@
      \def\leads@{\hrule\height\dimen@\depth-\dimen@ii}\or
      \def\leads@{\hrule\height\dimen@\depth-\dimen@ii}\or
      \def\leads@{\hbox to\ten@\p@{%
       \leaders\hrule\height\dimen@\depth-\dimen@ii\hfil
       \hfil
      \leaders\hrule\height\dimen@\depth-\dimen@ii\hskip\z@ plus2fil\relax
       \hfil
       \leaders\hrule\height\dimen@\depth-\dimen@ii\hfil}}\or
     \def\leads@{\hbox{\hbox to\ten@\p@{\dimen@\mathaxis@\advance\dimen@1.2\p@
       \dimen@ii\dimen@\advance\dimen@ii-.4\p@
       \leaders\hrule\height\dimen@\depth-\dimen@ii\hfil}%
       \kern-\ten@\p@
       \hbox to\ten@\p@{\dimen@\mathaxis@\advance\dimen@-1.2\p@
       \dimen@ii\dimen@\advance\dimen@ii-.4\p@
       \leaders\hrule\height\dimen@\depth-\dimen@ii\hfil}}}\fi
    \fi
    \cleaders\leads@\hfil
    \ifnum\arrcount@=\m@ne\else\ifnum\arrcount@=\thr@@\else
     \arrow@i
     \ifE@
      \ifnum\tcount@=\m@ne
      \else
       \ifcase\tcount@\char119 \or\or\char119 \or\char120 \or\char121 \or
       \char122 \or \or\or\or\or\char123 \or\char124 \or
       \char125 \or\char119 \hskip-\thr@@\p@\char119 \hskip\thr@@\p@\fi
      \fi
     \else
      \ifcase\scount@\or\or\char120 \or\char119 \or\or\or\char121 \or\char122
      \or\char123 \or\char124 \or\or\or\char125 \or
      \char119 \hskip-\thr@@\p@\char119 \hskip\thr@@\p@\fi
     \fi
    \fi\fi
    \harrow@ee}}%
  \harrow@e}%
 \iflabel@i
  \dimen@ii\z@\setbox\ZER@\hbox{$\m@th\tsize@@\label@i$}%
  \ifnum\arrcount@=\m@ne
  \else
   \advance\dimen@ii\mathaxis@
   \advance\dimen@ii\dp\ZER@\advance\dimen@ii\tw@\p@
   \ifnum\arrcount@=\thr@@\advance\dimen@ii\tw@\p@\fi
  \fi
  \advance\dimen@ii\pdimen@
  \next@{\harrow@b\smash{\raise\dimen@ii\hbox to\dimen@
   {\harrow@bb\hskip\tw@\ldimen@i\hfil\box\ZER@\hfil\harrow@ee}}\harrow@e}%
 \fi
 \iflabel@ii
  \ifnum\arrcount@=\m@ne
  \else
   \setbox\ZER@\hbox{$\m@th\tsize@\label@ii$}%
   \dimen@ii-\ht\ZER@\advance\dimen@ii-\tw@\p@
   \ifnum\arrcount@=\thr@@\advance\dimen@ii-\tw@\p@\fi
   \advance\dimen@ii\mathaxis@\advance\dimen@ii\pdimen@
   \next@{\harrow@b\smash{\raise\dimen@ii\hbox to\dimen@
    {\harrow@bb\hskip\tw@\ldimen@ii\hfil\box\ZER@\hfil\harrow@ee}}\harrow@e}%
  \fi
 \fi}
\let\tsize@\tsize
\def\tsizeCDlabels{\let\tsize@\tsize}
\def\ssizeCDlabels{\let\tsize@\ssize}
\def\tsize@@{\ifnum\arrcount@=\m@ne\else\tsize@\fi}
\def\varrow@{\dimen@\secondy@\advance\dimen@-\firsty@
 \ifN@\else\multiply\dimen@\m@ne\fi
 \setbox\ZER@\vbox to\dimen@
  {\ifN@\vskip-\Ydimen@\else\vskip\ydimen@\fi
   \ifnum\arrcount@=\m@ne\else\ifnum\arrcount@=\thr@@\else
    \hbox{\arrow@iii
     \ifN@
      \ifnum\tcount@=\m@ne
      \else
       \ifcase\tcount@\char117 \or\or\char117 \or\char118 \or\char119 \or
       \char120 \or\or\or\or\or\char121 \or\char122 \or\char123 \or
       \vbox{\hbox{\char117}\nointerlineskip\vskip\thr@@\p@
       \hbox{\char117}\vskip-\thr@@\p@}\fi
      \fi
     \else
      \ifcase\scount@\or\or\char118 \or\char117 \or\or\or\char119 \or
      \char120 \or\char121 \or\char122 \or\or\or\char123 \or
      \vbox{\hbox{\char117}\nointerlineskip\vskip\thr@@\p@
      \hbox{\char117}\vskip-\thr@@\p@}\fi
     \fi}%
    \nointerlineskip
   \fi\fi
   \ifnum\arrcount@=\m@ne
    \let\leads@\null
   \else
    \ifcase\arrcount@\let\leads@\vrule\or\let\leads@\vrule\or
    \def\leads@{\vbox to\ten@\p@{%
     \hrule\height1.67\p@\depth\z@\width.4\p@
     \vfil
     \hrule\height3.33\p@\depth\z@\width.4\p@
     \vfil
     \hrule\height1.67\p@\depth\z@\width.4\p@}}\or
    \def\leads@{\hbox{\vrule\height\p@\hskip\tw@\p@\vrule}}\fi
   \fi
  \cleaders\leads@\vfill\nointerlineskip
   \ifnum\arrcount@=\m@ne\else\ifnum\arrcount@=\thr@@\else
    \hbox{\arrow@iv
     \ifN@
      \ifcase\scount@\or\or\char118 \or\char117 \or\or\or\char119 \or
      \char120 \or\char121 \or\char122 \or\or\or\arrow@iii\char123 \or
      \vbox{\hbox{\char117}\nointerlineskip\vskip-\thr@@\p@
      \hbox{\char117}\vskip\thr@@\p@}\fi
     \else
      \ifnum\tcount@=\m@ne
      \else
       \ifcase\tcount@\char117 \or\or\char117 \or\char118 \or\char119 \or
       \char120 \or\or\or\or\or\char121 \or\char122 \or\arrow@iii\char123 \or
       \vbox{\hbox{\char117}\nointerlineskip\vskip-\thr@@\p@
       \hbox{\char117}\vskip\thr@@\p@}\fi
      \fi
     \fi}%
   \fi\fi
   \ifN@\vskip\ydimen@\else\vskip-\Ydimen@\fi}%
 \ifN@
  \dimen@ii\firsty@
 \else
  \dimen@ii-\firsty@\advance\dimen@ii\ht\ZER@\multiply\dimen@ii\m@ne
 \fi
 \rlap{\smash{\hskip\tocenter@\hskip\pdimen@\raise\dimen@ii\box\ZER@}}%
 \iflabel@i
  \setbox\ZER@\vbox to\dimen@{\vfil
   \hbox{$\m@th\tsize@@\label@i$}\vskip\tw@\ldimen@i\vfil}%
  \rlap{\smash{\hskip\tocenter@\hskip\pdimen@
  \ifnum\arrcount@=\m@ne\let\next@\relax\else\let\next@\llap\fi
  \next@{\raise\dimen@ii\hbox{\ifnum\arrcount@=\m@ne\hskip-.5\wd\ZER@\fi
   \box\ZER@\ifnum\arrcount@=\m@ne\else\hskip\tw@\p@\fi}}}}%
 \fi
 \iflabel@ii
  \ifnum\arrcount@=\m@ne
  \else
   \setbox\ZER@\vbox to\dimen@{\vfil
    \hbox{$\m@th\tsize@\label@ii$}\vskip\tw@\ldimen@ii\vfil}%
   \rlap{\smash{\hskip\tocenter@\hskip\pdimen@
   \rlap{\raise\dimen@ii\hbox{\ifnum\arrcount@=\thr@@\hskip4.5\p@\else
    \hskip2.5\p@\fi\box\ZER@}}}}%
  \fi
 \fi
}
\newdimen\goal@
\newdimen\shifted@
\newcount\Tcount@
\newcount\Scount@
\newbox\shaft@
\newcount\slcount@
\def\getcos@#1{%
 \ifnum\tan@i<\tan@ii
  \dimen@#1%
  \ifnum\slcount@<8 \count@9 \else \ifnum\slcount@<12 \count@8 \else
   \count@7 \fi\fi
  \multiply\dimen@\count@\divide\dimen@\ten@
  \dimen@ii\dimen@\multiply\dimen@ii\tan@i\divide\dimen@ii\tan@ii
 \else
  \dimen@ii#1%
  \count@-\slcount@\advance\count@24
  \ifnum\count@<8 \count@9 \else \ifnum\count@<12 \count@8
   \else\count@7 \fi\fi
  \multiply\dimen@ii\count@\divide\dimen@ii\ten@
  \dimen@\dimen@ii\multiply\dimen@\tan@ii\divide\dimen@\tan@i
 \fi}
\newdimen\adjust@
\def\Nnext@{\ifN@\let\next@\raise\else\let\next@\lower\fi}
\def\arrow@@{\slcount@\angcount@
 \ifNESW@
  \ifnum\angcount@<\ten@
   \let\arrowfont@\arrow@i\global\advance\angcount@\m@ne
   \global\multiply\angcount@13
  \else
   \ifnum\angcount@<19
    \let\arrowfont@\arrow@ii\global\advance\angcount@-\ten@
    \global\multiply\angcount@13
   \else
    \let\arrowfont@\arrow@iii\global\advance\angcount@-19
    \global\multiply\angcount@13
  \fi\fi
  \Tcount@\angcount@
 \else
  \ifnum\angcount@<5
   \let\arrowfont@\arrow@iii\global\advance\angcount@\m@ne
   \global\multiply\angcount@13 \global\advance\angcount@65
  \else
   \ifnum\angcount@<14
    \let\arrowfont@\arrow@iv\global\advance\angcount@-5
    \global\multiply\angcount@13
   \else
    \ifnum\angcount@<23
     \let\arrowfont@\arrow@v\global\advance\angcount@-14
     \global\multiply\angcount@13
    \else
     \let\arrowfont@\arrow@i\global\angcount@117
  \fi\fi\fi
  \ifnum\angcount@=117 \Tcount@115 \else\Tcount@\angcount@\fi
 \fi
 \Scount@\Tcount@
 \ifE@
  \ifnum\tcount@=\z@\advance\Tcount@\tw@\else\ifnum\tcount@=13
   \advance\Tcount@\tw@\else\advance\Tcount@\tcount@\fi\fi
  \ifnum\scount@=\z@\else\ifnum\scount@=13 \advance\Scount@\thr@@\else
   \advance\Scount@\scount@\fi\fi
 \else
  \ifcase\tcount@\advance\Tcount@\thr@@\or\or\advance\Tcount@\thr@@\or
  \advance\Tcount@\tw@\or\advance\Tcount@6 \or\advance\Tcount@7
  \or\or\or\or\or\advance\Tcount@8 \or\advance\Tcount@9 \or
  \advance\Tcount@12 \or\advance\Tcount@\thr@@\fi
  \ifcase\scount@\or\or\advance\Scount@\thr@@\or\advance\Scount@\tw@\or
  \or\or\advance\Scount@4 \or\advance\Scount@5 \or\advance\Scount@\ten@
  \or\advance\Scount@11 \or\or\or\advance\Scount@12 \or\advance
  \Scount@\tw@\fi
 \fi
 \ifcase\arrcount@\or\or\global\advance\angcount@\@ne\else\fi
 \ifN@\shifted@\firsty@\else\shifted@-\firsty@\fi
 \ifE@\else\advance\shifted@\charht@\fi
 \goal@\secondy@\advance\goal@-\firsty@
 \ifN@\else\multiply\goal@\m@ne\fi
 \setbox\shaft@\hbox{\arrowfont@\char\angcount@}%
 \ifnum\arrcount@=\thr@@
  \getcos@{1.5\p@}%
  \setbox\shaft@\hbox to\wd\shaft@{\arrowfont@
   \rlap{\hskip\dimen@ii
    \smash{\ifNESW@\let\next@\lower\else\let\next@\raise\fi
     \next@\dimen@\hbox{\arrowfont@\char\angcount@}}}%
   \rlap{\hskip-\dimen@ii
    \smash{\ifNESW@\let\next@\raise\else\let\next@\lower\fi
      \next@\dimen@\hbox{\arrowfont@\char\angcount@}}}\hfil}%
 \fi
 \rlap{\smash{\hskip\tocenter@\hskip\firstx@
  \ifnum\arrcount@=\m@ne
  \else
   \ifnum\arrcount@=\thr@@
   \else
    \ifnum\scount@=\m@ne
    \else
     \ifnum\scount@=\z@
     \else
      \setbox\ZER@\hbox{\ifnum\angcount@=117 \arrow@v\else\arrowfont@\fi
       \char\Scount@}%
      \ifNESW@
       \ifnum\scount@=\tw@
        \dimen@\shifted@\advance\dimen@-\charht@
        \ifN@\hskip-\wd\ZER@\fi
        \Nnext@
        \next@\dimen@\copy\ZER@
        \ifN@\else\hskip-\wd\ZER@\fi
       \else
        \Nnext@
        \ifN@\else\hskip-\wd\ZER@\fi
        \next@\shifted@\copy\ZER@
        \ifN@\hskip-\wd\ZER@\fi
       \fi
       \ifnum\scount@=12
        \advance\shifted@\charht@\advance\goal@-\charht@
        \ifN@\hskip\wd\ZER@\else\hskip-\wd\ZER@\fi
       \fi
       \ifnum\scount@=13
        \getcos@{\thr@@\p@}%
        \ifN@\hskip\dimen@\else\hskip-\wd\ZER@\hskip-\dimen@\fi
        \adjust@\shifted@\advance\adjust@\dimen@ii
        \Nnext@
        \next@\adjust@\copy\ZER@
        \ifN@\hskip-\dimen@\hskip-\wd\ZER@\else\hskip\dimen@\fi
       \fi
      \else
       \ifN@\hskip-\wd\ZER@\fi
       \ifnum\scount@=\tw@
        \ifN@\hskip\wd\ZER@\else\hskip-\wd\ZER@\fi
        \dimen@\shifted@\advance\dimen@-\charht@
        \Nnext@
        \next@\dimen@\copy\ZER@
        \ifN@\hskip-\wd\ZER@\fi
       \else
        \Nnext@
        \next@\shifted@\copy\ZER@
        \ifN@\else\hskip-\wd\ZER@\fi
       \fi
       \ifnum\scount@=12
        \advance\shifted@\charht@\advance\goal@-\charht@
        \ifN@\hskip-\wd\ZER@\else\hskip\wd\ZER@\fi
       \fi
       \ifnum\scount@=13
        \getcos@{\thr@@\p@}%
        \ifN@\hskip-\wd\ZER@\hskip-\dimen@\else\hskip\dimen@\fi
        \adjust@\shifted@\advance\adjust@\dimen@ii
        \Nnext@
        \next@\adjust@\copy\ZER@
        \ifN@\hskip\dimen@\else\hskip-\dimen@\hskip-\wd\ZER@\fi
       \fi	
      \fi
  \fi\fi\fi\fi
  \ifnum\arrcount@=\m@ne
  \else
   \loop
    \ifdim\goal@>\charht@
    \ifE@\else\hskip-\charwd@\fi
    \Nnext@
    \next@\shifted@\copy\shaft@
    \ifE@\else\hskip-\charwd@\fi
    \advance\shifted@\charht@\advance\goal@-\charht@
    \repeat
   \ifdim\goal@>\z@
    \dimen@\charht@\advance\dimen@-\goal@
    \divide\dimen@\tan@i\multiply\dimen@\tan@ii
    \ifE@\hskip-\dimen@\else\hskip-\charwd@\hskip\dimen@\fi
    \adjust@\shifted@\advance\adjust@-\charht@\advance\adjust@\goal@
    \Nnext@
    \next@\adjust@\copy\shaft@
    \ifE@\else\hskip-\charwd@\fi
   \else
    \adjust@\shifted@\advance\adjust@-\charht@
   \fi
  \fi
  \ifnum\arrcount@=\m@ne
  \else
   \ifnum\arrcount@=\thr@@
   \else
    \ifnum\tcount@=\m@ne
    \else
     \setbox\ZER@
      \hbox{\ifnum\angcount@=117 \arrow@v\else\arrowfont@\fi\char\Tcount@}%
     \ifnum\tcount@=\thr@@
      \advance\adjust@\charht@
      \ifE@\else\ifN@\hskip-\charwd@\else\hskip-\wd\ZER@\fi\fi
     \else
      \ifnum\tcount@=12
       \advance\adjust@\charht@
       \ifE@\else\ifN@\hskip-\charwd@\else\hskip-\wd\ZER@\fi\fi
      \else
       \ifE@\hskip-\wd\ZER@\fi
     \fi\fi
     \Nnext@
     \next@\adjust@\copy\ZER@
     \ifnum\tcount@=13
      \hskip-\wd\ZER@
      \getcos@{\thr@@\p@}%
      \ifE@\hskip-\dimen@\else\hskip\dimen@\fi
      \advance\adjust@-\dimen@ii
      \Nnext@
      \next@\adjust@\box\ZER@
     \fi
  \fi\fi\fi}}%
 \iflabel@i
  \rlap{\hskip\tocenter@
  \dimen@\firstx@\advance\dimen@\secondx@\divide\dimen@\tw@
  \advance\dimen@\ldimen@i
  \dimen@ii\firsty@\advance\dimen@ii\secondy@\divide\dimen@ii\tw@
  \global\multiply\ldimen@i\tan@i\global\divide\ldimen@i\tan@ii
  \ifNESW@\advance\dimen@ii\ldimen@i\else\advance\dimen@ii-\ldimen@i\fi
  \setbox\ZER@\hbox{\ifNESW@\else\ifnum\arrcount@=\thr@@\hskip4\p@\else
   \hskip\tw@\p@\fi\fi
   $\m@th\tsize@@\label@i$\ifNESW@\ifnum\arrcount@=\thr@@\hskip4\p@\else
   \hskip\tw@\p@\fi\fi}%
  \ifnum\arrcount@=\m@ne
   \ifNESW@\advance\dimen@.5\wd\ZER@\advance\dimen@\p@\else
    \advance\dimen@-.5\wd\ZER@\advance\dimen@-\p@\fi
   \advance\dimen@ii-.5\ht\ZER@
  \else
   \advance\dimen@ii\dp\ZER@
   \ifnum\slcount@<6 \advance\dimen@ii\tw@\p@\fi
  \fi
  \hskip\dimen@
  \ifNESW@\let\next@\llap\else\let\next@\rlap\fi
  \next@{\smash{\raise\dimen@ii\box\ZER@}}}%
 \fi
 \iflabel@ii
  \ifnum\arrcount@=\m@ne
  \else
   \rlap{\hskip\tocenter@
   \dimen@\firstx@\advance\dimen@\secondx@\divide\dimen@\tw@
   \ifNESW@\advance\dimen@\ldimen@ii\else\advance\dimen@-\ldimen@ii\fi
   \dimen@ii\firsty@\advance\dimen@ii\secondy@\divide\dimen@ii\tw@
   \global\multiply\ldimen@ii\tan@i\global\divide\ldimen@ii\tan@ii
   \advance\dimen@ii\ldimen@ii
   \setbox\ZER@\hbox{\ifNESW@\ifnum\arrcount@=\thr@@\hskip4\p@\else
    \hskip\tw@\p@\fi\fi
    $\m@th\tsize@\label@ii$\ifNESW@\else\ifnum\arrcount@=\thr@@\hskip4\p@
    \else\hskip\tw@\p@\fi\fi}%
   \advance\dimen@ii-\ht\ZER@
   \ifnum\slcount@<9 \advance\dimen@ii-\thr@@\p@\fi
   \ifNESW@\let\next@\rlap\else\let\next@\llap\fi
   \hskip\dimen@\next@{\smash{\raise\dimen@ii\box\ZER@}}}%
  \fi
 \fi
}
\def\outCD@#1{\def#1{\Err@{\noexpand#1must not be used within \string\CD}}}
\newskip\preCDskip@
\newskip\postCDskip@
\preCDskip@\z@
\postCDskip@\z@
\def\preCDspace#1{\RIfMIfI@
 \onlydmatherr@\preCDspace\else\advance\preCDskip@#1\relax\fi\else
 \onlydmatherr@\preCDspace\fi}
\def\postCDspace#1{\RIfMIfI@
 \onlydmatherr@\postCDspace\else\advance\postCDskip@#1\relax\fi\else
 \onlydmatherr@\postCDspace\fi}
\def\predisplayspace#1{\RIfMIfI@
 \onlydmatherr@\predisplayspace\else
 \advance\abovedisplayskip#1\relax
 \advance\abovedisplayshortskip#1\relax\fi
 \else\onlydmatherr@\preCDspace\fi}
\def\postdisplayspace#1{\RIfMIfI@
 \onlydmatherr@\postdisplayspace\else
 \advance\belowdisplayskip#1\relax
 \advance\belowdisplayshortskip#1\relax\fi
 \else\onlydmatherr@\postdisplayspace\fi}
\def\PreCDSpace#1{\global\preCDskip@#1\relax}
\def\PostCDSpace#1{\global\postCDskip@#1\relax}
\def\CD#1\endCD{%
 \outCD@\cgaps\outCD@\rgaps\outCD@\Cgaps\outCD@\Rgaps
 \preCD@#1\endCD
 \advance\abovedisplayskip\preCDskip@
 \advance\abovedisplayshortskip\preCDskip@
 \advance\belowdisplayskip\postCDskip@
 \advance\belowdisplayshortskip\postCDskip@
 \vcenter{\offinterlineskip
  \vskip\preCDskip@\Let@\global\colcount@\@ne\global\rowcount@\z@
  \everycr{%
   \noalign{%
    \ifnum\rowcount@=\Rowcount@
    \else
     \getrgap@\rowcount@\vskip\getdim@
     \global\advance\rowcount@\@ne\global\colcount@\@ne
    \fi}}%
  \tabskip\z@
  \halign{&\global\xoff@\z@\global\yoff@\z@
   \getcgap@\colcount@\hskip\getdim@
   \hfil\vrule\height\ten@\p@\width\z@\depth\z@
   $\m@th\displaystyle{##}$\hfil
   \global\advance\colcount@\@ne\cr
   #1\crcr}\vskip\postCDskip@}%
 \preCDskip@\z@\postCDskip@\z@
 \def\getcgap@##1{\ifcase##1\or\getdim@\z@\else\getdim@\standardcgap\fi}%
 \def\getrgap@##1{\ifcase##1\getdim@\z@\else\getdim@\standardrgap\fi}%
 \let\Width@\relax\let\Height@\relax\let\Depth@\relax\let\Rowheight@\relax
 \let\Rowdepth@\relax\let\Colwidth@\relax
}
\let\enddocument\bye
\def\alloc@#1#2#3#4#5{\global\advance\count1#1by\@ne
  \ch@ck#1#4#2%
  \allocationnumber=\count1#1%
  \global#3#5=\allocationnumber
  \wlog{\string#5=\string#2\the\allocationnumber}}
\catcode`\@=\active

\font\blackboard=msbm10
\font\sblackboard=msbm7
\font\ssblackboard=msbm5

\font\sanss=cmss10
\font\ssanss=cmss8 scaled 900
\font\sssanss=cmss8 scaled 600
\font\smss=cmss8 scaled 800
\font\teneurm=eurm10
\newfam\bbfam
\textfont\bbfam=\blackboard
\scriptfont\bbfam=\sblackboard
\scriptscriptfont\bbfam=\ssblackboard
\def\bb#1{{\fam\bbfam\relax#1}}

\newfam\ssfam
\textfont\ssfam=\sanss
\scriptfont\ssfam=\ssanss
\scriptscriptfont\ssfam=\sssanss
\def\ss#1{{\fam\ssfam\relax#1}}

\def\sss#1{\hbox{\ssanss #1}}
\def\eu#1{\hbox{$\teneurm \relax #1$}}
\catcode`\"=12
\mathchardef\za="710B  
\mathchardef\zb="710C  
\mathchardef\zg="710D  
\mathchardef\zd="710E  
\mathchardef\zve="710F 
\mathchardef\zz="7110  
\mathchardef\zh="7111  
\mathchardef\zvy="7112 
\mathchardef\zi="7113  
\mathchardef\zk="7114  
\mathchardef\zl="7115  
\mathchardef\zm="7116  
\mathchardef\zn="7117  
\mathchardef\zx="7118  
\mathchardef\zp="7119  
\mathchardef\zr="711A  
\mathchardef\zs="711B  
\mathchardef\zt="711C  
\mathchardef\zu="711D  
\mathchardef\zvf="711E 
\mathchardef\zq="711F  
\mathchardef\zc="7120  
\mathchardef\zw="7121  
\mathchardef\ze="7122  
\mathchardef\zy="7123  
\mathchardef\zvp="7124  
\mathchardef\zvr="7125 
\mathchardef\zvs="7126 
\mathchardef\zf="7127  
\mathchardef\zG="7000  
\mathchardef\zD="7001  
\mathchardef\zY="7002  
\mathchardef\zL="7003  
\mathchardef\zX="7004  
\mathchardef\zP="7005  
\mathchardef\zS="7006  
\mathchardef\zU="7007  
\mathchardef\zF="7008  
\mathchardef\zC="7009  
\mathchardef\zW="700A  
\catcode`\"=\active
\loadmsam
\TagsOnRight

\input paper.st\relax
\hsize=33pc
\hoffset=40pt
\vsize=54pc
\voffset=6pt
\def\*{{\textstyle *}}
\def\s*{{\scriptstyle *}}
\def\proof{\demo{Proof}}
\def\endproof{\hfill \vrule height4pt width6pt depth2pt \enddemo}
\def\R{{\bb R}}

\def\r{\ss r}
\def\sT{\ss T}
\def\ssT{\sss T}
\def\sG{\ss G}
\def\sH{\ss H}
\def\smT{\hbox{\smss T}}
\def\ezt{\eu\zt}
\def\ezp{\eu\zp}
\def\ezk{\eu\zk}
\def\eza{\eu\za}

\def\xd{\operatorname{d}\!}
\def\sdt{\text{\sevenrm d}_{\smT}}
\def\dt{\xd_{\sss T}}
\def\vt{\text{\tenrm v}_{\sss T}}

\def\idt{\operatorname{id}}

\def\xi{\operatorname{i}}

\def\ll{\text{\it \$}}
\def\T{\Cal T}
\def\V{\Cal V}
\def\H{\Cal G}
\def\I{\Cal J}
\def\i{\Cal J^\*}
\def\G{\Cal H}
\def\lb{[\mkern -7mu [\mkern 7mu}
\def\rb{]\mkern -3mu ] \mkern 3mu}
\def\ad{\operatorname{ad}}

\newsymbol\leqs 1336
\newsymbol\geqs 133E
\overfullrule=0pt

                \title
                Tangent and cotangent lifts  and graded Lie algebras\\
associated with Lie algebroids
                \endtitle

    \author
                J. Grabowski$^1$ and
                P. Urba\'nski$^2$
        \endauthor
        \date{ \today}
        \affil
     $^1$ Institute of Mathematics, University of Warsaw\\
        Banacha 2, 02-097 Warsaw, Poland\\and\\
        Mathematical Institute, Polish Academy of Sciences\\
        \'Sniadeckich 8, P.~O. Box 137, 00-950 Warsaw, Poland\\
     {\sanss email:\ jagrab\@mimuw.edu.pl}\\
        \\
        $^2$    Division of Mathematical Methods in Physics,
                University of Warsaw \\
                Ho\D{z}a 74, 00-682 Warsaw, Poland\\
        {\sanss email:\ urbanski\@fuw.edu.pl}
        \endaffil
     \address{
     Institute of Mathematics, University of Warsaw, Banacha 2,
     02-097 Warsaw, Poland }

     \address{
     Division of Mathematical Methods in Physics, University of
     Warsaw, Ho\D{z}a 74, 00-682 Warsaw, Poland }

      \thanks{
    Supported by KBN, grant No 2 PO3A 074 10 }

   \abstract
        Generalized Schouten, Fr\"olicher-Nijenhuis and
Fr\"olicher-Richardson brackets are defined for an arbitrary
Lie algebroid.
        Tangent and cotangent lifts of Lie algebroids are introduced
and discussed and the behaviour of the related graded Lie
brackets under these lifts is studied.  In the case of the
canonical Lie algebroid on the tangent bundle, a new  common
generalization of the Fr\"olicher-Nijenhuis and the  symmetric
Schouten brackets, as well as  embeddings of the Nijenhuis-Richardson
and the Fr\"olicher-Nijenhuis bracket into the Schouten bracket, are
obtained.
   \endabstract
        \document
        \maketitle

        \bigskip
        { \eightrm Keywords:
     Lie algebroid, Schouten bracket, Fr\"olicher-Nijenhuis bracket,
     Poisson structure, tangent lift, cotangent lift}

    { \eightrm MSC 1991: 17B70, 53C15, 17B66, 58F05}

\vskip 1cm
        \hl1"0"{Introduction}
        There are three main ingredients of the presented paper: to
study natural graded Lie brackets of tensor fields associated
with Lie algebroids over  a manifold, procedures of lifting
of these tensor fields to tensor fields over derived manifolds (like
lifting tensor fields on $M$ to tensor fields on $\sT M$), and
the behaviour of the brackets under the lifting procedures. All
this is done in the setting of an arbitrary Lie algebroid. We
tried to make the presentation as systematic as possible, but
still most of the presented results seem to be new in the
literature (the most important can be found  in Sections
5-7). We work coordinate-free as long as it is
reasonable, but for the readers convenience we provide also
detailed calculations in local coordinates.

        Graded Lie algebras and  Lie super-algebras play an important role in
the contemporary geometry and mathematical physics. Fundamental
examples are, defined in a natural way, graded extensions of the
Lie bracket of vector fields: the Schouten (Schouten-Nijenhuis)
bracket and the Fr\"olicher-Nijenhuis bracket.
        Another example is  the Lie algebra of a Lie group.
A {\it Lie algebroid}  is a common generalization of the notion of tangent
bundle and the Lie algebra of a Lie group [20].  Lie
algebroids have the basic ingredients which make possible usual
constructions of differential geometry, like exterior derivative,
Lie derivative, contractions, Schouten bracket, etc. (cf. [22]).

        In the first section of the paper, we introduce natural
analogues of the Schouten, Fr\"olicher\-Ni\-jen\-huis and
Nijenhuis-Richardson brackets in the framework of an arbitrary
Lie algebroid.

        An important example of a  Lie algebroid is provided by the
cotangent bundle of a Poisson manifold. In Section 2, we present
basic facts concerning Lie brackets related to Poisson
structures, like Poisson bracket of functions, brackets of
1-forms, and their extensions.

        Definitions of tangent and cotangent lifts of a Lie
algebroid are given in Section 3, and in Section 4 we study
tangent lifts of associated tensor fields.  The results of Section
5 show the behaviour of the introduced graded brackets under the
vertical and tangent lifts and generalize results of [31] and
[10]. Similar questions for the cotangent lift, which seems to be
less functorial, are studied in Section 6.

        All above is applied in the last section to the case of the
canonical Lie algebroid -- the tangent bundle $\sT M$, in
relation with our previous work [10]. Among other lifts, there
is an embedding of the Fr\"olicher-Nijenhuis bracket on a
manifold $M$ into a graded Lie bracket of differential forms
introduced earlier by the first author [9], into the
Fr\"olicher-Nijenhuis bracket on  $\sT^\* M$ (discovered by
Dubois-Violette and Michor [7] and regarded as a part of a 'common
generalization' of the Fr\"olicher-Nijenhuis and the symmetric
Schouten brackets), and embeddings of the Nijenhuis-Richardson
and the Fr\"olicher-Nijenhuis bracket on $M$ into the Schouten
bracket on $\sT^\* M$.

        \hl1{Graded Lie brackets of a Lie algebroid}
        \newpre\tag{1.} \Reset\tag{1}
        Let $M$  be a manifold and let $\zt\colon E\rightarrow M$ be
a vector bundle. By $\zF(\zt)$ we denote the graded exterior
algebra generated by sections of $\zt$, $\zF(\zt) =\oplus_{k\in
\bb Z}\zF^k(\zt)$, where $\zF^k(\zt) = \zG(M,\bigwedge^k E)$ for
$k\geqs 0$ and  $\zF^k(\zt) =  \{0\}$ for $k<0$. The dual vector
bundle we denote by $\zp \colon E^{\textstyle *} \rightarrow M$.
For $\zt = \ezt_M \colon \sT M \rightarrow M$, we get the graded
algebra of multivector fields on $M$, and for $\zt =\ezp_M
\colon \sT^\* M \rightarrow M$- the graded algebra of
differential forms on $M$.

        Let  us consider a Lie algebroid structure  on $\zt\colon
E\rightarrow M$ with the Lie bracket $[\,,\,]$ and the anchor $\za_\zt
\colon E\rightarrow \sT M$.  The Lie bracket $[\,,\,]$ of sections
of $E$ induces the well-known generalization of the standard
calculus of differential forms and vector fields ({[20]}, [21],
{[22]}).

The exterior derivative $\xd_\zt \colon \zF^k(\zp) \rightarrow
\zF^{k+1}(\zp)$ is defined by the formula
                $$ \multline
        \xd_\zt\zm(X_1,\dots,X_{k+1}) = \sum_i (-1)^{i+1}
\za_\zt(X_i)(\zm(X_1,\dots,\widehat{X}_i,\dots ,X_{k+1})) +\\
 + \sum_{i<j} (-1)^{i+j}\zm ([X_i,X_j], X_1,\dots , \widehat{X}_i,
\dots ,\widehat{X}_j, \dots , X_{k+1}),
                        \endmultline            \tag \label{1.1}$$
        where $X_i \in \zF^1(\zt)$ and the hat over a symbol means
that this is to be omitted.   For $X \in \zF^k(\zt)$, the
contraction $\xi_X \colon  \zF^p(\zp) \rightarrow
\zF^{p-k}(\zp)$ is defined in the standard way and the Lie
differential operator
                $$ \ll_X\colon \zF^p(\zp)\rightarrow \zF^{p-k+1}(\zp)
                                $$
         is defined by the graded commutator
                $$ \ll_X = \xi_X \circ \xd_\zt - (-1)^k \xd_\zt \circ
\xi_X.
                                $$
        The following theorem contains a list of well-known
formulae.
        \claim{Theorem}{}
        Let  $\zm\in \zF^k(\zp),\ \zn\in \zF(\zp)$ and $X,Y \in
\zF^1(\zt)$. We have
        \list
        \item $\xd_\zt \circ \xd_\zt =0$,
        \item $\xd_\zt(\zm\wedge \zn) = \xd_\zt(\zm) \wedge \zn + (-1)^k\zm
\wedge \xd_\zt(\zn)$,
        \item $\xi_X(\zm\wedge \zn) = \xi_X(\zm) \wedge \zn + (-1)^k\zm
\wedge \xi_X(\zn)$,
        \item $\ll_X(\zm\wedge \zn) = \ll_X\zm\wedge \zn +
\zm\wedge \ll_X \zn$,
        \item $\ll_X\circ \ll_Y - \ll_Y \circ \ll_X =\ll_{[X,Y]}$,
        \item $\ll_X \circ \xi_Y - \xi_Y\circ \ll_X = \xi_{[X,Y]}$.
        \endlist
                 \endclaim
        The last formula can be generalized in the following way
(cf. [24], [16]).
        \claim{Theorem}{}
        For $X\in \zF^{k+1}(\zt) $ and $Y\in \zF^{l+1}(\zt)$ there is a
unique $Z \in \zF^{k+l+1}(\zt)$ (denoted by $[X,Y]$) such that
                $$  \ll_Y\circ \xi_X   - (-1)^{(k+1)l}\xi_X\circ \ll_Y =
- \xi_{Z}.
                                $$
        \endclaim
        In particular, we have for $X\in \zF^1(\zt)$ and $f\in
\zF^0(\zp) = C^\infty(M)$
                $$ [X,f] = \za_\zt(X)(f).
                                $$
        The bracket
                $$[\,,\,] \colon \zF^{k+1}(\zt)\times \zF^{l+1}(\zt)
\rightarrow \zF^{k+l+1}(\zt)
                                $$
is an extension of the Lie algebroid bracket on $\zF^1(\zt)$ and it is
called {\it the generalized Schouten bracket }of the Lie algebroid.
It defines a graded Lie algebra structure on $\zF(\zt)$ with the
shifted grading  $\zF(\zt) = \oplus_{k\in \bb Z}\zF^{(k)}
(\zt)$, where $\zF^{(k)}(\zt) = \zF^{k+1}(\zt)$.
        Moreover, for $X\in \zF^{(k)}(\zt)$, the adjoint action
$\ad_X = [X,\cdot ]$ is a graded derivation of rank $k$  of the
exterior algebra, i.~e.,
                $$ [X,Y\wedge Z] = [X,Y]\wedge Z  + (-1)^{kl}Y
\wedge [X,Z]
                        \tag \label{1.2}$$
          for $Y\in \zF^{l}(\zt)$.

        The fact, that the generalized Schouten bracket is a graded
        Lie bracket, means that it is
        \list
        \item graded skew-symmetric, i.~e., for $X\in \zF^{(k)}
(\zt)$, $Y\in \zF^{(l)}(\zt)$, we have
                $$ [X,Y] = - (-1)^{kl}[Y,X];
                \tag \label{1.3}$$

        \item satisfies the graded Jacobi identity, i.~e., for $X,Y$
as before and for  $Z\in \zF^{(m)}(\zt)$, we have
                $$ (-1)^{km}[[X,Y],Z] + (-1)^{lk} [[Y,Z],X] +
(-1)^{ml}[[Z,X]Y] = 0.
                        \tag \label{1.4}$$
        \endlist
        The Jacobi identity can be written in the form
                $$[X,[Y,Z]] - (-1)^{kl}[Y,[X,Z]] = [[X,Y],Z],
                \tag \label{1.5}$$
        which means that the graded commutator
                $$ [\ad_X, \ad_Y] = \ad_X\circ \ad_Y - (-1)^{kl}\ad_Y
\circ \ad_X
                                $$
        equals  $\ad_{[X,Y]}$.

        From \Ref{1.2}, it follows that for $X_1,\dots, X_k, Y_1,
\dots , Y_l \in \zF^1(\zt)$, we have
                $$ [X_1\wedge\dots \wedge X_k, Y_1\wedge \dots \wedge
Y_k] = \sum_{i,j}(-1)^{i+j} [X_i,Y_j] \wedge X_1\wedge
\dots\wedge \widehat{X}_i \wedge \dots \wedge X_k\wedge
Y_1\wedge \dots \wedge \widehat{Y}_j\wedge \dots \wedge Y_l
                \tag \label{1.6}$$
         (see [16], [24]).

        There is also a {\it symmetric Schouten bracket} which extends
the Schouten bracket on $ \zF^0(\zt)\oplus \zF^1(\zt)$ to
symmetric multisections. This bracket is defined by the following
formula, similar to \Ref{1.6} (see [7])
                $$ [X_1\vee\dots \vee X_k, Y_1\vee \dots \vee
Y_k] = \sum_{i,j}[X_i,Y_j] \vee X_1\vee
\dots\vee \widehat{X}_i \vee \dots \vee X_k\vee
Y_1\vee \dots \vee \widehat{Y}_j\vee \dots \vee Y_l.
                        \tag \label{1.7}$$
        The symmetric Schouten bracket provides the symmetric algebra
$\zF_s(\zt) = \oplus_{k\in \bb Z} \zF_s^k(\zt) $, $\zF_s^k(\zt)
=\zG(M,S^kE)$, with the abstract Poisson algebra structure
([2]).

        In the case of the canonical Lie algebroid structure in the
tangent bundle $\ezt_M\colon \sT M \rightarrow M$, we obtain the
classical skew-symmetric and symmetric Schouten brackets ([28], [25]).

        Before we pass to other brackets, let us introduce the space
$\zF^k_l(\zt) = \zG(M, \bigwedge^kE\otimes \bigwedge^lE^\*)$ of
tensor fields of mixed type. Of course, we can identify
$\zF^k_l(\zt)$ and $\zF^l_k(\zp)$. For   $K\in \zF^k_1(\zp)$, we
define in the natural way the contraction
                $$ \xi_K\colon \zF^n(\zp) \rightarrow \zF^{n+k-1}(\zp).
                                        $$
        For simple tensors  $K = \zm\otimes X$, where $\zm \in \zF^k
(\zp),\, X\in \zF^1(\zt) $, we put
                $$ \xi_K\zn = \zm\wedge \xi_X\zn .
                        \tag \label{1.8}$$
        The corresponding Lie differential is defined by the formula
                $$ \ll_K = \xi_K \circ \xd_\zt + (-1)^k\xd_\zt \circ \xi_K
                                        $$
        and, in particular,
                $$ \ll_{\zm\otimes X  } = \zm\wedge \ll_X + (-1)^k \xd_\zt
\zm \wedge \xi_X.
                        \tag \label{1.9}$$
        This definition is compatible with the previous one in the
case of $K\in \zF^0_1(\zp) = \zF^1(\zt)$.

        The contraction (insertion) $\xi_K$ can be extended to an operator
                $$ \xi_K\colon \zF^n_1(\zp) \rightarrow \zF_1^{n+k-1}(\zp)
                                        $$
        by the formula
                $$ \xi_K(\zm\otimes X) = \xi_K(\zm)\otimes X .
                                        $$
        \claim{Theorem}{}
        The bracket
                $$ [\,,\,]^{N-R} \colon \zF_1^{k+1}(\zp) \times
\zF_1^{l+1}(\zp) \rightarrow \zF_1^{k+l+1}(\zp),
                                        $$
given by the formula
                $$ [K,L]^{N-R} = \xi_KL- (-1)^{kl}\xi_LK,
                        \tag \label{1.10}$$
defines a graded Lie algebra structure on the graded space
$\zF_1(\zp) = \oplus_{k\in \bb Z} \zF_1^{(k)}(\zp)$,  where
$\zF_1^{(k)}(\zp) = \zF_1^{k+1}(\zp)$.

        For simple tensors $\zm\otimes X \in \zF_1^{k}(\zp)$ and
$\zn \otimes Y \in \zF_1^{l}(\zp)$, the following formula holds
                $$[\zm\otimes X, \zn\otimes Y]^{N-R} = \zm\wedge
\xi_X\zn \otimes  Y + (-1)^k \xi_Y \zm\wedge \zn \otimes X.
                        \tag \label{1.11}$$
        \endclaim
        \proof
        The bracket is obviously skew-symmetric. Now, let  $K\in
\zF_1^{(k)}(\zp)$, $L \in \zF_1^{(l)}(\zp)$ and $N\in
\zF_1^{(n)}(\zp)$. It is a simple and standard task to verify that
                $$ \xi_{[K,L]^{N-R}} = \xi_K\circ \xi_L - (-1)^{kl}\xi_L
\circ \xi_K.
                        \tag \label{1.12}$$
         From \Ref{1.12} we get
                $$\multline [K, [L,N]^{N-R}]^{N-R} = \xi_K (\xi_L N -
(-1)^{ln}\xi_N L) - (-1)^{k(l+n)}(\xi_L\circ \xi_N - (-1)^{ln}
\xi_N\circ \xi_L) K = \\
        = \xi_K\xi_LN -(-1)^{ln}\xi_K\xi_NL - (-1)^{k(l+n)}\xi_L\xi_N
K + (-1)^{kl+kn +ln}\xi_N \xi_L K
                        \endmultline    $$
and, after easy calculations,
                $$(-1)^{kn}[K,[L ,N ]^{N-R}]^{N-R} + (-1)^{lk}[L,[N ,
K]^{N-R}]^{N-R} + (-1)^{nl}[N,[K ,L ]^{N-R}]^{N-R} =0 .
                                $$
        \endproof
        The bracket $[\,,\,]^{N-R}$ is called  {\it the generalized
Nijenhuis-Richardson} bracket.

{\bf Remark.} The generalized  Nijenhuis-Richardson bracket is a
purely vector bundle bracket and does not depend on any additional Lie
algebroid structure. For $\zt =\zt_M$ we get the classical
Nijenhuis-Richardson bracket of vector-valued forms ([23], [24]
[26]).

        The  {\it generalized Fr\"olicher-Nijenhuis bracket} is
defined for simple tensors $\zm \otimes X \in \zF_1^{k}(\zp)$ and
$\zn \otimes Y \in \zF_1^{l}(\zp)$ by the formula
                $$\multline
        [\zm \otimes X, \zn \otimes Y]^{F-N} = \\
        =\zm\wedge \zn \otimes
[X,Y] + \zm\wedge \ll_X\zn\otimes Y -\ll_Y\zm \wedge \zn \otimes
X + (-1)^k(\xd_\zt\zm \wedge \xi_X\zn \otimes Y + \xi_Y\zm
\wedge \xd_\zt\zn \otimes X) = \\
        = (\ll_{\zm\otimes X}\zn) \otimes Y - (-1)^{kl}(\ll_{\zn
\otimes Y}\zm)\otimes X + \zm\wedge \zn \otimes [X,Y].
                \endmultline                    \tag \label{1.13}$$

        \claim{Theorem}{}
        The formula \Ref{1.13} defines a graded Lie bracket
$[\,,\,]^{F-N}$ on the graded space $\zF_1^{}(\zp) =\oplus_{k\in
\bb Z}\zF_1^{k}(\zp)$. Moreover,
                $$\ll_{[K,L]^{F-N}} = \ll_K\circ \ll_L -
(-1)^{kl}\ll_L\circ \ll_K
                        \tag \label{1.14}$$
         for $K\in \zF_1^{k}(\zp)$, $L\in \zF_1^{l}(\zp)$.
        \endclaim
        \proof
        The right-hand side in formula \Ref{1.13} is linear with
respect to $\zm,\ \zn,\ X,\ Y$, and for $f\in C^\infty(M)$ we have
                $$[f\zm\otimes X,\zn \otimes Y]^{F-N} = [\zm\otimes
fX,\zn\otimes Y]^{F-N}.
                                        $$
        It follows that this formula defines actually a graded skew
symmetric bracket
                $$[\,,\,]\colon \zF_1^{k}(\zp) \times \zF_1^{l}(\zp)
\rightarrow \zF_1^{k+l}(\zp) .
                                        $$
         Standard calculations show that \Ref{1.14} holds and we use
it to verify  the Jacobi identity. For $\zm\otimes X$,
$\zn\otimes Y$ as above and $\zvy \otimes Z \in\zF_1^{m}(\zp)$,
we have
                $$\multline
        [\zvy \otimes Z, [\zm\otimes X , \zn \otimes
Y]^{F-N}]^{F-N}= \\
        = \ll_{\zvy\otimes Z}\ll_{\zm\otimes X}\zn\otimes Y -
(-1)^{kl} \ll_{\zvy\otimes Z}\ll_{\zn\otimes Y}\zm\otimes X +
\ll_{\zvy \otimes Z}(\zm\wedge \zn) \otimes [X,Y]- \\
-(-1)^{(k+l)m}\left(\ll_{\zm\otimes X}\ll_{\zn\otimes Y} \zvy -
(-1)^{kl} \ll_{\zn\otimes Y}\ll_{\zm\otimes X} \zvy\right)
\otimes Z + \\
        + \zvy\wedge \ll_{\zm \otimes X}\zn \otimes [Z,Y] -
(-1)^{kl} \zvy\wedge \ll_{\zn\otimes Y}\zm\otimes [Z,Y] +
\zvy\wedge \zm\wedge \zn \otimes [Z,[X,Y]].
                        \endmultline    $$
        Now, with some effort and with the use of \Ref{1.14},
one checks that
                $$\multline
        (-1)^{ml}[\zvy\otimes Z, [\zm\otimes X, \zn\otimes
Y]^{F-N}]^{F-N} + (-1)^{km}[\zm\otimes , [\zn\otimes Y,
\zvy\otimes Z]^{F-N}]^{F-N} +\\ + (-1)^{kl}[\zn\otimes Y,
[\zvy\otimes Z, \zm\otimes X]^{F-N}]^{F-N} =0.
                        \endmultline    $$
        \endproof
        In the case of the canonical Lie algebroid  (the tangent
bundle), we obtain the classical Fr\"olicher -Nijenhuis bracket on
vector-valued forms ([8], [13], [7], [23]). Let us make
a remark that, since, in general, the correspondence $K\rightarrow
\ll_K$ is not injective,  the identity \Ref{1.14} cannot be used
as a definition of  the generalized Fr\"olicher-Nijenhuis bracket.

        All the introduced graded Lie algebra structures have their
purely algebraic versions. Lie algebroids correspond  (depending
on authors) to pseudo-Lie algebras, Lie-Rinehart algebras,
Lie-Cartan pairs etc. ([16], [29]).

        \hl1{Graded Lie brackets on Poisson manifolds}
        \newpre\tag{2.} \Reset\tag{1}
        Let us consider a Poisson structure on a manifold $M$, i.~e., a
bi-vector field $P\in \zF^2(\zt_M)$ such that $[P,P] = 0$. The tensor $P$
defines a Lie algebra structure on $ C^\infty(M)$ by the
Poisson bracket
                $$ \{f,g\}_P = \xi_P (\xd f\wedge \xd g).
                                        $$
          It is known ([3], [12], [18], [19])
          that $P$ induces also a Lie algebroid
structure on  the cotangent bundle $\ezp_M\colon \sT^\*M
\rightarrow M$.
The anchor map of this structure is the corresponding to $P$
morphism of vector bundles $\widetilde{P} \colon \sT^\*M
\rightarrow \sT M $. The Lie bracket  is given by the
formula
                $$ [\zm,\zn]_P = \ll_{\widetilde{P}(\zm)}\zn -
\ll_{\widetilde{P}(\zn)} \zm -  \xd(\xi_P(\zm\wedge \zn))
                                \tag \label{2.1}$$
         and the Lie algebroid exterior derivative $\xd_{\ezp_M}
\colon \zF(\ezt_M) \rightarrow  \zF(\ezt_M)$ is given by the
Schouten bracket
                $$ \xd_{\ezp_M}(X) = [P,X]
                                \tag \label{2.2}$$
([2], [16], [11]).

        The Lie algebroid structure on $\sT^\* M$ induces a
generalized Schouten bracket $[\,,\,]_P$ of differential
forms --  {\it the Koszul-Schouten bracket} and a generalized
Fr\"olicher-Nijenhuis bracket $[\,,\,]^{F-N}_P$ of multivector valued
$1$-forms, described in the previous section.

        The vector bundle morphism $\widetilde{P} \colon \sT^\* M
\rightarrow \sT M$ is a morphism of Lie algebroids and,
consequently, the induced mapping
                $$\ss \zL_P= \bigwedge\widetilde{P}\colon \zF(\ezp_M)
\rightarrow \zF(\ezt_M),\ \ \bigwedge ^k\widetilde{P} \colon
\bigwedge^k\sT^\* M\rightarrow \bigwedge ^k\sT M,
                                        $$
                $$ \ss \zL_P(\zm_1\wedge \dots\wedge \zm_k) =
\widetilde{P}(\zm_1)\wedge \dots \wedge \widetilde{P}(\zm_k)
                                                \tag \label{2.3}$$
         gives us a homomorphism of the corresponding Schouten
brackets ([17], [14], [15]).

        \claim{Theorem}{}
        \list
        \item"(a)" The mapping $\ss \zL_P$ is a homomorphism between the
Koszul-Schouten and the Schouten brackets
                $$ \ss \zL_P\colon \left(\zF(\ezp_M),[\,,\,]_P \right)
\rightarrow
\left(\zF(\ezt_M), [\,,\,]\right) .
                                        $$
        It is an isomorphism if and only if $\widetilde{P}$ is
invertible, i.\,e., if and only if $P$ defines a symplectic structure.
        \item"(b)"If $\widetilde{P}$ is invertible, then the formula
                $$ \ss F_P(\zm\otimes X) = \ss \zL_P(\zm)\otimes
\ss \zL_P^{-1}(X)
                                        $$
         defines an isomorphism of the Fr\"olicher-Nijenhuis brackets
                $$\ss F_P\colon \left( \zF_1(\ezp_M), [\,,\,]^{F-N}\right)
\rightarrow  \left( \zF_1(\ezt_M), [\,,\,]^{F-N}_P\right).
                                        $$
        \endlist
        \endclaim
        \proof
        Part (a) is essentially due to Koszul ([17]). Part (b) is
obvious, because $\widetilde{P}$ is an isomorphism of Lie
algebroids and, consequently, respects the Lie algebroid brackets,
exterior derivatives, contractions, Lie differentials and,
finally, the Fr\"olicher-Nijenhuis brackets.
        \endproof

        Beside the homomorphism $\ss \zL_P$ of graded commutative
algebras, we can consider a module valued graded derivation of
degree $-1$, introduced by Michor [23]:
                $$\ss R_P \colon \zF^k(\ezp_M) \rightarrow \zF^{k-1}_1
(\ezp_M),
                                        $$
        characterized by the properties
                $$ \ss R_P|_{\zF^0(\ezp_M)}=0,\ \ \ss R_P|_{\zF^1(\ezp_M)}=
\widetilde{P},
                                        $$
        and
                $$ \ss R_P(\zm\wedge \zn) = \ss R_P(\zm)\wedge \zn +
(-1)^k\zm\wedge \ss R_P(\zn),
                        \tag \label{2.4}$$
        for $\zm \in \zF^k(\ezp_M)$, so that for 1-forms
$\zm_1,\dots,\zm_k$ we have
                $$\ss R_P\colon (\zm_1\wedge \cdots\wedge \zm_k) = \sum_i
(-1)^{i+1}\zm_1\wedge \cdots \wedge \widehat{\zm_i}\wedge \cdots
\wedge \zm_k \otimes \widetilde{P}(\zm_i).
                        \tag \label{2.5}$$
        We can define $\ss R_P$ in a more intrinsic way. Let
$\widetilde{\zm} \colon \bigwedge^{k-1}\sT M\rightarrow \sT^\* M$
be the bundle morphism corresponding to $\zm \in \zF^k(\ezp_M)$.
The morphism $\widetilde{\ss R_P(\zm)} \colon \bigwedge^{k-1} \sT M
\rightarrow \sT M$ corresponding to $\ss R_P(\zm) \in
\zF^{k-1}_1(\ezp_M)$ is defined by the formula
                $$ \widetilde{\ss R_P(\zm)} = \widetilde{P} \circ
\widetilde{\zm}.
                                                        $$
           Let us define the hamiltonian map $\sH_P = \ss R_P\circ
\xd$ and  the total hamiltonian map $\sG_P =\ss
\zL_P\circ\xd\,$. We have
                $$ \sH_P\colon \zF^{k}(\ezp_M) \rightarrow \zF^{k}_1(\ezp_M)
                                                        $$
and
                $$ \sG_P\colon \zF^{k}(\ezp_M) \rightarrow \zF^{k+1}(\ezt_M).
                                                        $$

        For $f\in  C^\infty(M)$, $\sH_P(f) = \sG_P(f) =
\widetilde{P}(\xd f) = -[P,f]$ is the hamiltonian vector field
corresponding to $f$. For a fixed $P$, we shall write also
$\ss R_\zm, \ \sH_\zm$,  and $\sG_\zm$, instead of $\ss R_P(\zm), \
\sH_P(\zm)$ and $\sG_P(\zm)$. The Koszul-Schouten bracket may be
written now in the form
                $$ [\zm,\zn]_P = \xi_{\sH_\zm}\zn - (-1)^k\ll_{\ss R_\zm}\zn
                                        \tag\label{2.6}$$
for $\zm \in \zF^{k}(\ezp_M)$ ([1], [9]).

        There is another graded Lie bracket on $\zF(\ezp_M)$, which
extends directly the Poisson bracket $\{\,,\,\}_P$ of functions.

        \claim{Theorem}{([9])}
        On the graded space $\zF^{}(\ezp_M) = \oplus_{k\in \bb Z}
\zF^{k}(\ezp_M)$, the formula
                $$ \{\zm,\zn\}_P = \ll_{\sss H_\zm}\zn +\xd \ll_{\sss R_\zm}
\zn
                                        \tag \label{2.7}$$
        defines  a graded Lie algebra structure such that on
functions  $\{,\}_P$ is the Poisson bracket and $\{\xd\zm,
\zn\}_P = \xd \{\zm,\zn \}_P$. Moreover, for $g_0,\dots,g_k,
f_0,\dots, f_l \in  C^\infty (M)$, we have
                $$ \multline \{g_0 \xd g_1\wedge \cdots\wedge \xd g_k
, f_0\xd f_1\wedge \cdots \wedge \xd f_l \}_P = \\
        = \{g_0,f_0\}\xd g_1\wedge \cdots\wedge\xd g_k \wedge \xd f_1
\wedge \cdots\wedge\xd f_l \\
        - g_0\sum_{i>0,j}(-1)^{i+j}\xd\{g_i,f_j\} \wedge \xd g_1
\wedge \cdots\wedge \widehat{\xd g_i}\wedge \cdots\wedge \xd g_k
\wedge \xd f_0 \wedge \cdots\wedge \widehat{\xd f_j}\wedge
\cdots\wedge \xd f_l\\
        - (-1)^kf_0\sum_{j>0,i}(-1)^{i+j}\xd\{g_i,f_j\} \wedge \xd g_0
\wedge \cdots\wedge \widehat{\xd g_i}\wedge \cdots\wedge \xd g_k
\wedge \xd f_1 \wedge \cdots\wedge \widehat{\xd f_j}\wedge
\cdots\wedge \xd f_l\\
        - \sum_{i,j>0}(-1)^{i+j}\{g_i,f_j\}  \xd g_0
\wedge \cdots\wedge \widehat{\xd g_i}\wedge \cdots\wedge \xd g_k
\wedge \xd f_0 \wedge \cdots\wedge \widehat{\xd f_j}\wedge
\cdots\wedge \xd f_l.
                                \endmultline    \tag \label{2.8}$$
        \label{t6} \endclaim

        {\bf Remark.} The bracket \Ref{2.7} is similar to brackets
introduced by  Michor ([23]) and Cabras and Vinogradov ([6]).
Their brackets, however, are either not skew-symmetric or do not
satisfy the Jacobi identity, so that they are not graded Lie
brackets. These brackets  define a proper graded Lie bracket on
co-exact forms which coincides with \Ref{2.7} projected to
co-exact forms. The following theorem shows, how the Poisson
bracket \Ref{2.7} is related to other graded Lie brackets
associated with  the Poisson structure.
        \claim{Theorem}{({\rm [9]})}
        \list
        \item"(a)" The exterior derivative
                $$ \xd\, \colon \left(\zF^{}(\ezp_M),\{\,,\,\}_P\right)
\rightarrow  \left(\zF^{}(\ezp_M),[\,,\,]_P\right)
                                        $$
        is a homomorphism of the Koszul-Schouten bracket into the
bracket \Ref{2.7}, i.~e.,
                $$ [\xd \zm,\xd \zn]_P = \xd\, \{\zm,\zn\}_P.
                                        $$
        The bracket $\{\,,\,\}_P$ may be then understood as an
\lq integral' of the Koszul-Schouten bracket.

        \item"(b)" The hamiltonian map
                $$\sH_P \colon \left(\zF^{}(\ezp_M),\{\,,\,\}_P\right)
\rightarrow  \left(\zF^{}_1(\ezp_M), [\,,\,]^{F-N}\right)
                                                $$
is a homomorphism of $\{\,,\,\}_P$ into the Fr\"olicher-Nijenhuis bracket:
                $$ \sH_P(\{\zm,\zn\}_P) = [\sH_P(\zm),\sH_P(\zn)]^{F-N}.
                                                $$
        \item "(c)" The total hamiltonian map
        $$ \sG_P \colon \left(\zF^{}(\ezp_M),\{\,,\,\}_P\right)
\rightarrow  \left( \zF^{}(\ezt_M), [\,,\,] \right)
                                                $$
        is a homomorphism of $\{\,,\,\}_P$ into the Schouten bracket:
                $$ \sG_P (\{\zm,\zn\}_P) = [\sG_P(\zm), \sG_P(\zm)].
                                                $$
        \item"(d)" If $\widetilde{P}$ is invertible, i~.e., if $P$ defines a
symplectic form $\zw$, then the isomorphism
                $$\ss \zL_P \colon \zF^{}(\ezp_M) \rightarrow \zF^{}(\ezt_M)
                                                $$
transports $\{\,,\,\}_P$ to an \lq integral' $\{\,,\,\}_\zw$ of
the Schouten bracket.
        \endlist
        \endclaim

        \hl1{Tangent and cotangent lifts of Lie algebroids \label{3}}
        \newpre\tag{3.} \Reset\tag{1}
        Let  $\zt\colon E\rightarrow M$ be a vector bundle with a
Lie algebroid structure and let $\za_\zt \colon  E\rightarrow \sT M$ be
the anchor map. Let $P$ be the corresponding linear Poisson
structure on the dual bundle $\zp \colon E^\* \rightarrow M$.
Linearity of this Poisson structure means that the bracket of
functions which are linear on fibers is linear on fibers ([22],
[4], [5]). There is a one-to-one correspondence $X\rightarrow \zi(X)$
between sections of $E$ and linear functions on $E^\*$:
$\zi(\zm) = \langle X(\zp(\zm)),\zm\rangle $.
The relation between the Lie algebroid  and the Poisson
bracket is given by
                $$ \zi([X,Y]) = \{\zi(X), \zi(Y)\}_P.
                                        \tag \label{3.1}$$
         For each $X\in \zF^{1}(\zt)$, we define a vector field
$\H(X)$ on $E^\*$ putting
                $$ \H(X) = \widetilde{P}(\xd\zi(X)) = \sH_P(\zi(X)) =
-[P,\zi(X)].
                                        $$
        The vector field $\H(X)$ is projectable and
                $$ \za_\zt(X) = \zp_\* \H(X),\ \text{i.~e.,}\
\zp^\*(\za_\zt(X)(f)) = \{\zi(X), \zp^\*f\}_P.
                                        \tag \label{3.2}$$
The tangent space $\sT E^\*$ carries the structure of a double
vector bundle ([27], [20], [21])
        $$\CD \sT E^\* @()\L{\sT\zp} @(1,0) @() \L{\ezt_{E^\s*}}@(0,-1)  & \sT
M
@()\l{\ezt_M} @(0,-1)  \\
                        E^\* @() \L{\zp}@(1,0)   &    M
        \endCD \quad .                  \tag\label{3.3}$$
        There is a natural lift of the Poisson tensor $P$ on $E^\*$
to a Poisson tensor $\dt P$ on $\sT E^\*$ ([5], [10]) which is
linear with respect to both vector bundle structures on $\sT
E^\*$. It follows, that we have two bundles with Lie algebroid
structures - dual bundles with respect to the horizontal and
vertical bundle structures in \Ref{3.3}.

        The bundle dual to $\ezt_{E^\s*} \colon \sT E^\*
\rightarrow E^\*$ is $\ezp_{E^\s*}\colon \sT^\* E^\* \rightarrow
E^\*$ with the canonical pairing. The bundle dual to $\sT \zp
\colon \sT E^\* \rightarrow \sT M$ can be identified with the
bundle $\sT \zt\colon \sT E \rightarrow \sT M$. The pairing is
the tangent pairing obtained from the canonical pairing
$\langle ,\rangle \colon E\times _ME^\* \rightarrow \R$ by
applying the tangent functor ([30], [10]).

        The Lie algebroid structures on  $\ezp_{E^\s*}\colon \sT^\*
E^\* \rightarrow E^\*$ and $\sT \zt\colon \sT E \rightarrow \sT
M$ are called the {\it tangent} and {\it cotangent} lifts of the Lie
algebroid structure on $E$.

Let us see all this in local coordinates.
 Let $(x^a)$ be a local coordinate system on $M$ and  let
$e_1,\dots,e_m$ be basis of local sections of $E$. The dual basis
of local sections of $E^\*$ we denote by ${e^*}^1,\dots,{e^*}^m$.
The correspondig coordinates in $E^\*$ we denote by $(x^a,
\zx_i)$ and the coordinates in $E$ by $(x^a, y^i)$, i.~e.,
$\zi(e_i) = \zx_i$ and $\zi({e^*}^i) = y^i$.
        Since $P$ is a linear Poisson structure, we have
                $$ P= \frac{1}{2}c_{ij}^k(x) \zx_k\partial _{\zx_i}\wedge
\partial _{\zx_j}  + \zd_i^a(x) \partial _{\zx_i}\wedge \partial _{x^a},
                                        \tag \label{3.4}$$ where $c^k_{ij},
\zd_i^a$
depend on $x$ only and where we clearly use the summation
convention. The anchor map is given by
                $$\za_\zt(y^ie_i(x)) = y^i\zd_i^a(x)\partial _{x^a}
                                        \tag \label{3.5}$$
and the Lie bracket by
                $$ [e_i,e_j] =c_{ij}^ke_k,
                                        \tag \label{3.6}$$
so that
                $$ [f^ie_i,g^je_j] = \left( f^ig^j c^k_{ij} + f^i
\zd_i^a\frac{\partial g^k}{\partial x^a} -g^j\zd_j^a\frac{\partial
f^k}{\partial x^a} \right)e_k.          \tag \label{3.7}$$
        The exterior derivative  $\xd_\zt$ is determined by
                $$ \xd_\zt(f) =\zd_i^a \frac{\partial f}{\partial x^a}
{e^*}^i,\ \ \xd _\zt({e^*}^k) = \frac{1}{2}c^k_{ji} {e^*}^i\wedge
{e^*}^j.
                                        \tag \label{3.8}$$
        Let us introduce the adopted coordinate systems
        \list
        \item"" \phantom{and} $(x^a, {\dot x}^b)$ in $\sT M$,
        \item"" \phantom{and} $(x^a, y^i,{\dot x}^b, {\dot y}^j)$ in $\sT E$,
        \item""  and $(x^a, \zx_i,{\dot x}^b, {\dot \zx}_j)$ in $\sT E^\*$.
        \endlist
        In these coordinates, we have ([5], [10])
                $$ \multline
                \dt P = c_{ij}^k(x)\zx_k \partial _{\zx_i} \wedge
\partial _{{\dot \zx}_j} + \frac{1}{2} c_{ij}^k(x){\dot \zx}_k
\partial _{{\dot \zx}_i} \wedge \partial _{{\dot \zx}_j} + \zd_i^a
\partial _{\dot \zx_i} \wedge \partial _{x^a} + \\
        + \zd_i^a(x)\partial _{\zx_i}\wedge \partial _{{\dot x}^a} +
\frac{\partial \zd_i^a}{\partial x^b}(x){\dot x}^b \partial
_{\dot\zx_i} \wedge \partial _{{\dot x}^a} + \frac{1}{2}
\frac{\partial c_{ij}^k}{\partial x^b}(x) {\dot x}^b
\zx_k\partial _{\dot \zx_i}\wedge \partial _{{\dot \zx}_j} .
                                \endmultline            \tag \label{3.9}$$
        For sections $\xd x^a$ and $\xd \zx_i$ of the Lie algebroid
bundle $\ezp_{E^\s*} \colon  \sT ^\* E^\*\rightarrow E^\*$, we have
$\zi(\xd x^a) = {\dot x}^a$, $\zi (\xd \zx_i) = \dot \zx_i$, and
                $$ \split  [\xd x^a, \xd x^b] =& 0,\\
                [\xd \zx_i, \xd x^a] = &  \frac{\partial
\zd_i^a}{\partial x^b}(x) \xd x^b, \\
                [\xd \zx_i, \xd \zx_j] =& c_{ij}^k (x) \xd \zx_k +
\frac{\partial c^k_{ij}}{\partial x^b} \zx_k \xd x^b .
                        \endsplit               \tag \label{3.10}$$
         The anchor map is given by
                $$\split \za_{\ezp_{E^*}}(\xd x^a) &=
-\zd_i^a(x)\partial _{\zx_i},\\
                \za_{\ezp_{E^*}}(\xd \zx_i) &= \zd_i^a(x) \partial _{x^a}
+ c_{ij}^k(x) \zx_k \partial_{\zx_j},
                        \endsplit       \tag \label{3.11}$$

and the exterior derivative $\xd_{\ezp_{E^*}}$ on
$\zF(\ezt_{E^\s*})$ is determined by
                $$\split \xd_{\ezp_{E^*}}(f) & =  \left(\zd_i^a
\frac{\partial f}{\partial x^a} + c_{ij}^k \zx_k \frac{\partial
f}{\partial \zx_j}  \right) \partial _{\zx_i} - \zd_i^a
\frac{\partial f}{\partial \zx_i} \partial _{x^a}, \\
        \xd_{\ezp_{E^*}}(\partial _{x^a}) & = - \frac{\partial
\zd_i^b}{ \partial x^a} \partial_{\zx_i}\wedge \partial _{x^b} -
\frac{1}{2}   \frac{\partial c_{ij}^k}{\partial x^a} \zx_k
\partial _{\zx_i}\wedge \partial _{\zx_j},\\
        \xd_{\ezp_{E^*}}(\partial _{\zx_k}) &=
\frac{1}{2}c_{ji}^k\partial _{\zx_i} \wedge \partial _{\zx_j}.
                        \endsplit               \tag \label{3.12}$$
        Now, let us turn to the Lie algebroid bundle $\sT \zt \colon
\sT E \rightarrow \sT M$ and its dual $\sT \zp \colon  \sT E^\*
\rightarrow \sT M$.
We choose a basis $(\overline{e}_i, \dot e_j)$ of local
sections of $\sT \zt \colon  \sT E \rightarrow \sT M$, such that
$\zi (\overline{e}_i) = \zx_i$ and $\zi(\dot e_j) = \dot \zx _j$,
and  a similar basis $({\overline{e}}^{\* j}, {\dot e}^{\* j} )$ of
local sections of $\sT \zp \colon  \sT E^\* \rightarrow \sT M$,
such that $\zi ({{\overline{e}}^\*}^i) = y^i$, $\zi ({\dot
e}^{\* j}) = {\dot y}^j$. In the chosen bases, section
$\overline{e}_i$ is dual to ${\dot e}^{\* i}$ and $\dot e_j$ is
dual to ${\overline{e}^\*}^j$, with respect to the pairing between
$\sT E$ and $\sT E^\*$.

        The anchor map $\za_{\ssT \zt} \colon
\sT E \rightarrow \sT\sT M$ is given by
                $$ \split \za_{\ssT \zt}(\dot e_i) =& \zd_i^a(x)
\partial _{x^a} + \frac{\partial \zd_i^a}{\partial x^b}(x) {\dot
x}^b \partial _{{\dot x}^a}, \\
                        \za_{\ssT \zt}(\overline{e}_i) =& \zd_i^a(x)
\partial _{{\dot x}^a},
                        \endsplit       \tag \label{3.13}$$
        and the Lie bracket by
                $$ \split
                [\overline{e}_i,\overline{e}_j] = & 0 ,\\
        [\overline{e}_i, \dot e_j] = & c_{ij}^k(x) \overline{e}_k, \\
        [\dot e_i , \dot e_j] = & c_{ij}^k(x) \dot e_k +
\frac{\partial c_{ij}^k}{\partial x^a}(x) \dot x^a \overline{e}_k .
                        \endsplit       \tag \label{3.14}$$
        The exterior derivative is determined by

\newpage
                $$ \split
        \xd_{\ssT \zt} (f) =& \left( \zd_i^a\frac{\partial f}{
\partial x^a} + \frac{\partial \zd_i^a}{\partial x^b}{\dot x}^b
\frac{\partial f}{\partial {\dot x}^a}\right)
{\overline{e}^\*}^i + \zd_i^a \frac{\partial f}{{\dot x}^a}{\dot
e}_i^\* ,\\
        \xd_{\ssT \zt}({\overline{e}^\*}^k) =& \frac{1}{2}
c_{ij}^k(x) {\overline{e}^\*}^i \wedge {\overline{e}^\*}^j ,\\
        \xd_{\ssT \zt}({\dot e}^{\* k}) =& c_{ij}^k(x) {\dot
e}^{\* i} \wedge  {\overline{e}^\*}^j + \frac{1}{2}\frac{
\partial c_{ij}^k}{\partial x^a}(x) {\dot x}^a
{\overline{e}^\*}^i \wedge {\overline{e}^\*}^j.
                        \endsplit       \tag \label{3.15}$$

        To make our presentation more complete, we end this section
with a simple theorem on iterated lifts.

        \claim""{Proposition}{}
        Let $\zt\colon E\rightarrow M$ be a Lie algebroid. Then the Lie
algebroid structures provided by tangent-cotangent,
cotangent-tangent and cotangent-cotangent lifts are canonically isomorphic.
        \endclaim
        \proof
        The tangent and cotangent lifts of the algebroid correspond
to the tangent lift $\dt P$ of the Poisson structure $P$ on the dual bundle
$ \zp \colon E^\* \rightarrow M $. Thus the iterated lifts
correspond to the tangent of the tangent Poisson structure,
i.~e., the Poisson structure  $\dt\dt P$ on  $\sT \sT E^\*$.
Since $\dt P$ is linear with respect to both vector bundle
structures on $\sT E^\*$, the Poisson structure $\dt\dt P$ is
linear with respect to the three vector bundle structures on $\sT\sT
E^\*$.

         The Lie algebroid on the bundle  $\ezp_{\ssT E^*}\colon \sT^\*
\sT E^\* \rightarrow \sT E^\* $, dual to $\ezt_{\ssT E^*} \colon
\sT \sT E^\* \rightarrow \sT E^\*$, is, by  definition, the
cotangent lift of the tangent lift of the Lie algebroid and, at the
same time, the cotangent lift of   the cotangent lift of the
Lie algebroid.

        The bundle dual to $\sT\ezt_{ E^*} \colon \sT \sT E^\*
\rightarrow \sT E^\* $ is the bundle $\sT \ezp_{E^*} \colon
\sT \sT^\* E^\* \rightarrow \sT E^\*$ and the corresponding Lie
algebroid structure is the tangent lift of the tangent lift of
the algebroid $\zt \colon E\rightarrow M$.

        There is well known canonical isomorphism (see Section~7)
$\ezk_{E^*}$ between $\ezt_{\ssT E^*} \colon \sT \sT E^\* \rightarrow
\sT E^\*$ and $\sT\ezt_{ E^*} \colon \sT \sT E^\* \rightarrow
\sT E^\* $, and the Poisson structure $\dt\dt P$ is invariant
with respect to it (cf. [10]).  It follows, that the dual bundles
are isomorphic and this isomorphism is an isomorphism of Lie
algebroid structures.
        \endproof

        \hl1{Tangent lifts of sections of vector bundles associated
with a Lie algebroid}
        \newpre\tag{4.} \Reset\tag{1}
        In [31] (see also the references there) there are described
different types of lifts of tensor fields from a manifold to the
tangent and cotangent bundles. In [10], the tangent lifts were
studied from the categorial point of view and were applied to
Poisson geometry (Courant in [4], [5] did the first steps in that
direction). All these lifts are related to the canonical
Lie algebroid structure in the tangent bundle $\ezt_M \colon \sT M
\rightarrow M $. Here we  show, how the lifting procedures
can be generalized to the case of an arbitrary Lie algebroid.

        Let us recall first  lifting processes in the case of the
canonical Lie algebroid. For a function $f \in  C^\infty (M)$, we
define the vertical and complete lifts $\vt(f), \dt(f) \in
C^\infty (\sT M) $:
                $$ \vt(f) = \ezt_M^\* (f),\ \ \  \dt(f)(v) = \langle
v,\xd f\rangle ,
                                        \tag \label{4.1}$$
        i.~e., $\vt(f)$ is the pull-back and $\dt(f)$ is the
exterior derivative $\xd f$, regarded as a function on $\sT M$.
        The vertical and complete lifts $\vt(X), \dt(X)$ of a vector
field $X$ on $M$ can be defined as unique vector fields on $\sT
M$ such that
                $$ \split
        \vt(X)(\vt(f)) &=0,\\
        \vt(X)(\dt(f)) &= \dt(X)(\vt(f)) = \vt(X(f)), \\
        \dt(X)(\dt(f)) &= \dt(X(f))
                        \endsplit       \tag \label{4.2}$$
(see [31], [10]).

        In local coordinates ($(x^a) $ on $M$ and $(x^a, {\dot
x}^b)$ on $\sT M$), we have
                $$ \vt(f)(x,\dot x) = f(x),\ \ \ \dt(f)(x,\dot x) =
\frac{\partial f}{ \partial x^a}(x){\dot x}^a,
                                        \tag \label{4.3}$$

         and

                $$\vt(f^a\partial _{x^a}) = f^a\partial _{{\dot x}^a},
\ \  \dt(f^a\partial _{x^a}) = f^a \partial _{x^a} +
\frac{\partial f^a}{ \partial x^b}{\dot x}^b \partial _{{\dot
x}^a}.
                                        \tag \label{4.4}$$

        Let $(\zt_j\colon E_j \rightarrow M)$ be a family of vector
bundles. We are going to define the {\it complete} and {\it
vertical} lifts of sections of the bundle $\underset i \to
\otimes{}_ME_j \rightarrow M $ to sections  of the bundle $ \underset j
\to \otimes{} _{\ssT M} \sT E_j \rightarrow  \sT M$.

        In the case of a section $X$ of a vector bundle $\zt\colon
E\rightarrow M$, we can proceed in the following way. The section $X$ defines
a linear function $\zi(X)$ on $E^\*$. The functions $\dt(\zi(X))$ and
$\vt(\zi(X))$  on $\sT E^\*$ are linear with respect to
the vector bundle structure $\sT \zp\colon \sT E^\* \rightarrow
\sT M$. It follows that these functions define sections  of the
dual bundle $\sT\zt \colon \sT E \rightarrow \sT M$. We denote
by $\T (X)$ the section which corresponds to $\dt (\zi(X))$ and by
$\V (X)$ the section which corresponds to $\vt(\zi(X))$.

        We can extend this procedure to the case of sections of the bundle
$\underset j \to \otimes{}_ME_j \rightarrow M $.
        For every  such section $X$,  we define a multilinear
function $\widetilde{X} \colon  \underset j \to \times{}_M
E^\*_j \rightarrow \R$ by
                $$ \widetilde{X}(\zm_1,\dots,\zm_k) =
\xi_X(\zm_1\wedge \cdots \wedge \zm_k).
                                        $$
           In the case of $X\in C^\infty(M)$ understood as a section
of the 0-tensor product, we put $\widetilde{X} = X $.

        \claim{Theorem}{}
        For every section $X \in \zG(M ,\underset j \to
\otimes{}_ME_j)$, there are, uniquely determined, sections
$\V(X), \T(X) \in \zG(M, \underset j \to \otimes{}_{\sT M} \sT
E_j) $  such that
                $$\split
        \text{(a)}\qquad  \widetilde{\V(X)} & =\vt(\widetilde{X}),\\
        \text{(b)}\qquad  \widetilde{\T(X)} & =\dt(\widetilde{X})
                        \endsplit               \tag \label{4.5}$$
        (we identify $\sT (\underset j \to \times {}_M E_j)$ and
$\underset j \to \times {}_{\ssT M} \sT E _j$ ).

        In particular,
                $$\V(f) = \vt(f)\ \ \text{and}\ \  \T(f) = \dt f \ \
\text{for}\ f\in C^\infty(M),
                                        $$
        and
                $$\V(e_{1,i_1}\otimes \cdots \otimes e_{k,i_k}) =
\overline{e}_{1,i_1} \otimes \cdots \otimes
\overline{e}_{k,i_k}, \ \ \T(e_{1,i_1}\otimes \cdots \otimes
e_{k,i_k}) = \sum_j \overline{e}_{1,i_1} \otimes \cdots \otimes
{\dot e}_{j,i_j} \otimes \cdots \otimes \overline{e}_{k,i_k},
                                        $$
where $(e_{j,i})$ is a basis of local sections of $\zt_j$.
\endclaim
        \proof
        Let $(e_{j,i})$ be a basis of local sections of $\zt_j$ and
let $(x^a,y^i_j), (x^a, \zx_{j,i})$ be the corresponding
coordinate system in $E_j, E_j^\*, \ j=1,\dots,k$. These
coordinate systems induce a coordinate system $(x^a, \zx_{1,i_1},
\dots ,\zx_{k,i_k})$ on $\underset j \to \times {}_M E^\*_j$. In the
introduced coordinate systems and bases,
                $$ X = f^{i_1\dots i_k} e_{1,i_1}\otimes \cdots
\otimes e_{k,i_k}
                                        $$
and
                $$\widetilde{X} =  f^{i_1\dots i_k}\zx_{1,i_1}\cdots
\zx_{k, i_k}.
                                        \tag \label{4.6}$$

         In the adopted coordinates $(x^a, \zx_{1,i_1}, \dots
,\zx_{k,i_k},{\dot x}^a, {\dot \zx}_{1,i_1}, \dots ,{\dot
\zx}_{k,i_k})$ on $\underset j \to \times {}_{\ssT M} \sT E^\*_j $,
the vertical lift $\vt(\widetilde{X})$ looks formaly like
\Ref{4.6}. It follows that $\vt(\widetilde{X}) =
\widetilde{\V(X)} $, where
                $$ \V(X) =  f^{i_1\dots i_k}\overline{e}_{1,i_1}
\otimes \cdots \otimes \overline{e}_{k,i_k}.
                                        \tag \label{4.7}$$
         The complete lift $\dt(\widetilde{X})$ is given in local
coordinates by (cf. \Ref{4.3})
                $$
                \dt(\widetilde{X}) = \frac{\partial f^{i_1\dots i_k}}{
\partial x^a} {\dot x}^a \zx_{1,i_1} \cdots \zx_{k,i_k}
+ \sum_j f^{i_1\dots i_k} \zx_{1,i_1} \cdots {\dot \zx}_{j,i_j}
\cdots \zx_{k,i_k}
                                        $$
         and, consequently,
                $$ \T(X) = \frac{\partial f^{i_1\dots i_k}}{ \partial
x^a} {\dot x}^a \overline{e}_{1,i_1}\otimes \cdots \otimes
\overline{e}_{k,i_k} + \sum_j f^{i_1\dots i_k}
\overline{e}_{1,i_1}\otimes \cdots  \otimes {\dot e}_{j,i_j}
\dots \otimes \overline{e}_{k,i_k}.
                                        \tag \label{4.8}$$
        Let us notice that both lifts commute with permutations and,
consequently,  lifts of skew-sym\-met\-ric and
symmetric tensors are skew-symmetric and symmetric, respectively.
        \endproof

        There is another (equivalent) way to introduce the vertical
and complete lifts. Again, we begin with the case  of a section
$X \colon M\rightarrow E$ of $\zt $. Applying the tangent
functor to $X$, we obtain $\sT X\colon \sT M \rightarrow \sT E$.
Since $\zt \circ X = \idt_M$, we have $\sT \zt \circ \sT X =
\idt_{\ssT M}$, i.~e., $\sT M$ is a section of $\sT \zt \colon
\sT E \rightarrow \sT M$. It is easy to check  that $\sT X
= \T(X)$. In order to get a similar construction of $\T(X)$ for
$X\in \zG(M, \underset i \to \otimes{}_ME_i) $, we observe first
that applying the tangent functor to the tensor product map
$\otimes \colon \underset i \to \times {}_{ M}  E_i
\rightarrow  \underset i \to \otimes{}_ME_i$, we get a
bundle morphism with respect to the tangent bundle structures
                $$ \sT \otimes \colon \sT \underset i \to \times
{}_{ M}  E_i \simeq \underset i \to \times {}_{\ssT M}
\sT E_i \rightarrow \sT \underset i \to \otimes{}_ME_i.
                                        $$
         This morphism  is multilinear and, consequently,
defines a linear morphism
                $$ \otimes{}_{\ssT} \colon \underset i \to
\otimes{}_{\ssT M} \sT E_i \rightarrow \sT \underset i \to
\otimes{}_ME_i.
                                        $$
        We have also the dual morphism of the dual bundles
                $$ \otimes{}^{\ssT}\colon \sT \underset i \to
\otimes{}_ME^\*_i \rightarrow  \underset i \to \otimes{}_{\ssT
M} \sT E^\*_i .
                                        $$
        If we replace $E_i$  by $E^\*_i$, we obtain
                $$ \otimes^{\ssT}\colon \sT \underset i \to
\otimes{}_ME_i \rightarrow  \underset i \to \otimes{}_{\ssT
M} \sT E_i.
                                        $$
        The complete lift $\T(X)$ of a section $X \colon M
\rightarrow \underset i \to \otimes{}_ME_i$ can be now defined
by
                $$ \T(X) = \otimes^{\ssT}\circ \sT X .
                                        $$

        For the vertical lift $\V(X)$, we provide a different approach
which makes use of the structure of a double  vector bundle
on $\sT E$:
        $$\CD \sT E @()\L{\sT\zt} @(1,0) @() \L{\zt_{E}}@(0,-1)  & \sT M
@()\l{\ezt_M} @(0,-1)  \\
                        E @() \L{\zt}@(1,0)   &    M
        \endCD \quad .$$
        Since $\zt_E$ ($\sT \zt$) is a vector bundle morphism with
respect to the horizontal (vertical) vector bundle structures, we
have two kernel subbundles of these projections:
        \list
        \item"" $V_{\ssT}E = \ker \zt_E$ is a subbundle of $\sT \zt
\colon \sT E \rightarrow \sT M$,
        \item""  $V_{\zt}E = \ker \sT\zt $ is a subbundle of $ \zt_E
\colon \sT E \rightarrow  E$.
        \endlist
        It is known that there are canonical  vector bundle
isomorphisms
                $$ V_{\ssT }E \simeq \sT M \times _M E , \ \ \  V_\zt
\simeq E\times _ME,
                                        $$
         which can be extended to isomorphisms of tensor products
                $$  \split \sT M\times _M \underset i \to
\otimes{}_ME_i &\simeq \underset i \to \otimes{}_{\ssT M}
V_{\ssT}  E_i \subset \underset i \to \otimes{}_{\ssT M} \sT E_i,
\\
        E\times _M \otimes{} _M^kE & \simeq  \otimes{} _E^k V_\zt E
\subset \otimes{}_E^k \sT E .
                                \endsplit \tag \label{4.9}$$
        With these identifications, we can define vertical lifts
$\V_{\ssT}$ and $\V_\zt$ of sections of $\underset i \to
\otimes{}_ME_i$ and $\otimes{} _M^k E $ respectively. In local
coordinates, we have
                $$\split \V_{\ssT} (f^{i_1\dots i_k} e_{1,i_1}\otimes
\cdots \otimes e_{k,i_k}) &= f^{i_1\dots i_k}\overline{e}_{1,i_1}
\otimes \cdots \otimes \overline{e}_{k, i_k},\\
        \V_\zt (f^{j_1\dots j_k}e_{j_1}\otimes \cdots \otimes
e_{j_k}) & = f^{j_1\dots j_k}\partial _{y^{j_1}} \otimes \cdots
\otimes  \partial _{y^{j_k}}.
                        \endsplit       \tag \label{4.10}$$
         The vertical lift $\V_{\ssT}$ coincides with $\V$, the
vertical lift $\V_\zt$ will be used in the definition of
cotangent lifts.

        \claim{Theorem}{}
        Let $\zt_i \colon E_i\rightarrow M$ and $\zt'_j \colon E'_j
\rightarrow  M$ be vector bundles over $M$, and let $X \in
\zG(M,\underset i \to \otimes{}_ME_i)$  or $X=f\in C^\infty
(M)$, $Y\in \zG(M, \underset j \to \otimes{}_ME'_j)$. Then,

                $$\split
        \text{\tenrm (a)}\qquad \V(X\otimes _MY) &= \V(X)\otimes _{\ssT
M}\V(Y), \\
        \text{\tenrm (b)} \qquad \T(X\otimes _M Y) &= \T(X) \otimes _{\ssT
M}\V(Y) + \V(X) \otimes _{\ssT M}\T(Y).
                        \endsplit       \tag \label{4.11}$$
        \endclaim
        \proof
        Using the canonical projections
        $$\left(\underset i \to \times {}_{ M}  E_i^\*\right) \times
_M \left(\underset j \to \times {}_{ M}  {E'}^\*_j\right) \rightarrow
\underset i \to \times {}_{ M}  E_i^\*
                                        $$
        and
        $$\left(\underset i \to \times {}_{ M}  E_i^\* \right)\times
_M \left( \underset j \to \times {}_{ M}  {E'}^\*_j \right)
\rightarrow \underset j \to \times {}_{ M}  {E'}^\*_j,
                                        $$
        we can consider $\widetilde{X}$ and $\widetilde{Y}$ as
functions on $\left( \underset i \to \times {}_{ M}  E_i^\* \right)
\times _M \left( \underset j \to \times {}_{ M}  {E'}^\*_j \right)$.
We have then
                $$\widetilde{X\otimes _MY} = \widetilde{X}\cdot
\widetilde{Y}.
                                        $$
          The statements of the theorem are now direct consequences
of the following properties of the vertical and complete lifts of functions
(cf. [10]):
                $$\split \vt(\widetilde{X}\cdot \widetilde{Y}) &=
\vt(\widetilde{X}) \cdot \vt(\widetilde{Y}),\\
        \dt(\widetilde{X} \cdot \widetilde{Y} ) &= \dt(
\widetilde{X}) \cdot  \vt(\widetilde{Y}) + \vt(\widetilde{X})
\cdot \dt(\widetilde{Y}).
                        \endsplit       $$
        \endproof
        \claim""{Corollary}{} Let $\zt \colon E\rightarrow M$ be a
vector bundle.
        \list
        \item"(a)" The vertical lift
                $$\V \colon \zF^{}(\zt) \rightarrow \zF^{}(\sT
\zt)
                                        $$
          is a homomorphism of graded commutative algebras:
                $$ \V(X\wedge Y) = \V(X)\wedge \V(Y).
                                        $$
        \item"(b)" The complete lift
                $$ \T \colon \zF^{}(\zt) \rightarrow \zF^{}(\sT
\zt)
                                        $$
        is a graded $\V$-derivation of degree 0:
                $$ \T(X\wedge Y) = \T(X)\wedge \V(Y) + \V(X)\wedge
\T(Y).
                                        $$
        \endlist
        \endclaim

        \hl1{Tangent lifts and graded Lie brackets}
                \newpre\tag{5.} \Reset\tag{1}
        Throughout this section $\zt \colon E\rightarrow M$ is a
vector bundle of a Lie algebroid  with  the anchor $\za_\zt
\colon E\rightarrow \sT M$. We shall prove theorems describing
the behaviour of contractions, Lie differentials, graded Lie
brackets, etc., with respect to the tangent lifts.
        \claim{Theorem}{}
        Let   $X\in \zF^{1}(\zt)$. Then,
        \list
        \item"(a)"   $\za_{\ssT \zt}(\V(X)) = \vt (\za_\zt(X))$,
        \item"(b)"   $\za_{\ssT \zt}(\T(X)) = \dt (\za_\zt(X))$.
        \endlist
        \endclaim
        \proof
        (a)\ \  With the notation of Section 3, we have, in view of
\Ref{3.13},
                $$ \za_{\ssT \zt}(\V(f^ie_i)) = \za_{\ssT
\zt}(f^i\overline{e}_i) = f^i\za_{\ssT \zt}(\overline{e}_i) =
f^i\zd^a_i\partial _{{\dot x}^a}.
                                        $$
        On the other hand, it follows from \Ref{4.4} that
                $$ \vt(\za_{\ssT \zt}(f^ie_i)) = \vt( f_i\zd_i^a
\partial _{x^a}) = f^i\zd_i^a\partial _{{\dot x}^a}.
                                        $$

        (b)\ \ From \Ref{3.13} and \Ref{4.8}, we get
                $$\multline \za_{\ssT \zt}(\T(f^ie_i)) = \za_{\ssT
\zt}\left( \frac{\partial f^i}{\partial x^a}{\dot x}^a
\overline{e}_i  + f^i\dot e_i \right) = \frac{\partial
f^i}{\partial x^a} {\dot x}^a \za_{\ssT \zt} (\overline{e}_i) +
f^i \za_{\ssT \zt}(\dot e_i) =\\
        = \frac{\partial f^i}{\partial x^a}{\dot x}^a \zd_i^b \partial
_{{\dot x}^b} + f^i\left( \zd_i^a\partial _{x^a} + \frac{\partial
\zd_i^b}{\partial x^a} {\dot x}^a \partial _{{\dot x}^b}\right)
        = f^i \zd_i^a\partial _{x^a}  + \frac{\partial }{\partial x^a}
\left( f^i\zd_i^b\right) {\dot x}^a \partial _{{\dot x}^b}.
                \endmultline            $$
        On the other hand, from \Ref{4.4}, we have
                $$ \dt(\za_{\ssT \zt}(f^i e_i)) = \dt(f^i
\zd_i^b \partial _{x^b})
        = f^i\zd^b_i \partial _{x^b} + \frac{\partial }{\partial x^a}
\left( f^i \zd_i^b \right) {\dot x}^a \partial _{{\dot x}^b}
                                        $$
and the theorem follows.
        \endproof
        The behaviour of the Schouten bracket with respect to
vertical and complete lifts is described by the following
theorem.
        \claim{Theorem}{}
        Let  $X,Y \in \zF^{}(\zt)$. Then,
        \list
        \item"(a)" $[\V(X), \V(Y)] = 0$,
        \item"(b)" $[\V(X), \T(Y)] = [\T(X), \V(Y)] = \V([X,Y])$,
        \item"(c)" $[\T(X) ,\T(Y)] = \T([X,Y])$,
i.~e., the complete lift is a homomorphism of the Schouten
brackets.
        \endlist
        \endclaim
        \proof
        (a)\ \ \ We have for the tangent Lie algebroid (\Ref{3.13},
\Ref{3.14}) $\za_{\ssT \zt}(\overline{e_i}) =  \zd_i^a \partial
_{{\dot x}^a}$ and $[\overline{e}_i, \overline{e}_j] =0$. Since
$\V(f^ie_i) = f^i\overline{e}_i$, we get $[\V(X), \V(Y)] = 0$.

        (b)\ \ \ The proof goes by  induction with respect to the
degree of $X$ and $Y$. Let us take first $X\in \zF^{1}(\zt)$ and
$Y=f\in C^\infty(M)$. It follows from Theorem 10 and \Ref{4.2}
that
                $$[\V(X),\T(f)] = \za_{\ssT \zt}(\V(X))(\dt f) =
\vt(\za_\zt(X))(\dt f) = \vt(\za_\zt(X)(f)) = \V([X,f]).
                                        $$
Similarly, $[\T(X),\V(f)] = \V([X,f])$.

Assume now that $ X,Y  \in \zF^{1}(\zt)$. The tangent lift
$\{\,,\,\}_{\sdt P}$ of the Poisson bracket $\{\,,\,\}_P$ has the
following property (cf. [4], [5], [10])
                $$\vt(\{f,g\}_P) = \{\dt(f), \vt(g)\}_{\sdt P}=
\{\vt(f), \dt(g)\}_{\sdt P}.
                                        $$
        Since $\zi([X,Y]) = \{\zi(X), \zi(Y)\}_P$, we get by
definitions of the lifts $\V,\ \T$ (\Ref{4.5}) that
                $$\multline
        \zi([\V(X), \T(Y)]) = \{\zi(\V(X)), \zi(\T(Y))\}_{\sdt P} =
\{\vt(\zi(X)), \dt(\zi(Y))\}_{\sdt P}= \\ = \vt(\{\zi(X), \zi(Y)\}_P)=
\vt( \zi([X, Y])) = \zi(\V([X,Y])).
                                \endmultline    $$
        It follows that  $[\V(X), \T(Y)] = \V([X,Y])$ and that
                $$[\T(X), \V(Y)] = -[\V(Y),\T(X)] = -\V([X,Y]) =
\V([Y,X]).
                                        $$
        Now, it remains to show that if  (b) holds for $X\in
\zF^{k}(\zt), \ Y\in \zF^{l}(\zt)$, then it holds also for $X$
and $Y\wedge Z$, $Z\in \zF^{1}(\zt)$. We get from (a), from
Corollary to Theorem~9, and from \Ref{1.2}  that
                $$\multline
        [\V(X), \T(Y\wedge Z)] = [\V(X), \T(Y)\wedge \V(Z) + \V(Y)
\wedge \T(Z)] = \\ = [\V(X), \T(Y)] \wedge \V(Z) + (-1)^{(k-1)l}\T(Y)
\wedge [\V(X), \V(Z)] + [\V(X), \V(Y)] \wedge \T(Z)
+ \\ + (-1)^{(k-1)l}\V(Y) \wedge [\V(X), \T(Z)]
        = \V([X,Y]) \wedge \V(Z) + (-1)^{(k-1)l}\V(Y) \wedge \V([X
, Z]) = \\ = \V\left([X,Y]\wedge Z + (-1)^{(k-1)l} Y\wedge [X,Z]
\right) = \V([X,Y\wedge Z]).
                                \endmultline    $$
        The equality $[\T(X), \V(Y\wedge Z)] = \V([X,Y\wedge Z])$
follows in an analogous way.

        (c)\ \ \ As in (b), we check  the equality for $X \in
\zF^{1}(\zt)$ and $Y =f \in C^\infty(M)$. We have in this case
                $$ \multline
        [\T(X),\T(f)] = \za_{\ssT \zt}(\T(X))(\dt f) = \dt(
\za_\zt(X))(\dt f) = \\
        =\dt(\za_\zt(X)(f)) = \dt([X,f]) = \T([X,f]).
                                \endmultline    $$
        If $X,Y \in \zF^{1}(\zt)$, then
                $$ \multline
        \zi([\T(X), \T(Y)]) = \{\zi(\T(X)), \zi(\T(Y)) \}_{
\sdt P} = \{\dt\zi(X), \dt\zi(Y)\}_{\sdt P} =\\
        = \dt\{\zi(X), \zi(Y)\}_P = \dt(\zi([X,Y])) =
\zi(\T([X,Y])).
                                \endmultline    $$
        Now, the proof proceeds like in the part (b).
        \endproof

The following theorem contains a collection of formulae
describing the behaviour of inserting operators, exterior
derivatives, and Lie derivatives under the tangent lifts.
        \claim{Theorem}{}
        Let  $\zm \in \zF^{}(\zp)$ and $X \in \zF^{1}(\zt)$. Then,
         $$\spreadmatrixlines{5pt} \matrix \format \r\quad & \c & \quad \l \\
\text{\tenrm (1)}& \text{\tenrm a)}  & \xi_{\V(X)} \V(\zm) = 0,\\
                &       \text{\tenrm b)}  & \xi_{\V(X)} \T(\zm)  = \xi_{\T(X)}
\V(\zm) = \V(\xi_X\zm),\\
                &       \text{\tenrm c)}  & \xi_{\T(X)} \T(\zm) =
\T(\xi_X\zm),\\
        \text{\tenrm (2)}&
                        \text{\tenrm a)}  &\xd_{\ssT \zt}\V(\zm) =
\V(\xd_\zt\zm), \\
                &       \text{\tenrm b)} &\xd_{\ssT \zt}\T(\zm) =
\T(\xd_\zt\zm),\\
        \text{\tenrm (3)}&
                        \text{\tenrm a)} &\ll_{\V(X)} \V(\zm) = 0,\\
        &       \text{\tenrm b)}  &\ll_{\V(X)} \T(\zm) = \ll_{\T(X)}
\V(\zm) = \V(\ll_X\zm),\\
        &       \text{\tenrm c)} & \ll_{\T(X)} \T(\zm) = \T(\ll_X\zm).
                                        \endmatrix $$
        \endclaim
        \proof
        (1)\ \ \  As in the proof of Theorem~8, we choose a basis
$(e_i)$ of local sections  of $\zt \colon E\rightarrow M$. We
have then $X = f^je_j$, $\zm = \sum_{i_1<\cdots <i_k}
\zm_{i_1\dots i_k} e^{* i_1} \wedge \cdots e^{* i_k}$, $\V(X) =
f^j \overline{e}_j$, and
                $$\V(\zm) = \sum_{i_1<\cdots <i_k} \zm_{i_1\dots i_k}
\overline{e}^{* i_1} \wedge \cdots \overline{e}^{* i_k}.
                                                $$
Since the basis $({\dot e}^*_i, \overline{e}_j^*)$ is dual (with
respect to the tangent pairing) to the basis $(\overline{e}_i,
{\dot e}_j)$, we have  $\xi_{\V(X)} \V(\zm) = 0$.
        With
                $$\T(X)= \frac{\partial f^j}{\partial x^a}{\dot x}^a
\overline{e}_j + f^j {\dot e}_j,
                                                $$
we get
                $$\multline
                \xi_{\T(X)}\V(\zm) =  \sum_{i_1<\cdots <i_k} (-1)^{j+1}
\zm_{i_1\dots i_k} f^j \overline{e}^{* i_1} \wedge \cdots \wedge
\widehat{\overline{e}^{*i_j}} \wedge \cdots \wedge
\overline{e}^{* i_k}= \\
        = \V\left( \sum_{i_1<\cdots <i_k} (-1)^{j+1} \zm_{i_1\dots
i_k} f^j e^{* i_1} \wedge \cdots \wedge \widehat{e^{*i_j}}
\wedge \cdots \wedge  e^{* i_k}\right) = \V(\xi_X\zm).
                                \endmultline    $$
Similarly, we obtain for
                $$ \T(\zm) = \sum_{i_1<\cdots <i_k,j}
\zm_{i_1\dots i_k} \overline{e}^{* i_1} \wedge \cdots \wedge
{\dot e}^{*i_j} \wedge \cdots \wedge \overline{e}^{* i_k} +
        \sum_{i_1<\cdots <i_k,a} \frac{\partial \zm_{i_1\dots i_k}}{
\partial x^a } {\dot x}^a  \overline{e}^{* i_1} \wedge \cdots
\wedge \overline{e}^{* i_k}
                                                $$
that $\xi_{\V(X)}\T(\zm) = \V(\xi_X \zm)$.

        Finally,
                $$\multline
        \xi_{\T(X)}\T(\zm) = \sum_{i_1<\cdots< i_k,a,j } (-1)^{j+1}
\frac{\partial \zm_{i_1\dots i_k}}{ \partial x^a } {\dot x}^a
f^{i_j} \overline{e}^{*i_1} \wedge \cdots \wedge
\widehat{\overline{e}^{*i_j}} \wedge \cdots \wedge
\overline{e}^{*i_k}  + \\
+
        \sum_{i_1<\cdots< i_k, s\neq j}(-1)^{j+1}\zm_{i_1\dots i_k}
f^{i_j} \overline{e}^{*i_1} \wedge \cdots \wedge \widehat{
\overline{e}^{*i_j}} \wedge \cdots \wedge {\dot e}^{*i_s}\wedge
\dots \wedge \overline{e}^{*i_k}
+ \\ +
        \sum_{i_1<\cdots< i_k,a,j} \zm_{i_1\dots i_k} \frac{\partial f^{i_j}}{
\partial x^a } {\dot x}^a \overline{e}^{*i_1} \wedge \cdots \wedge
\widehat{\overline{e}^{*i_j}} \wedge \cdots \wedge
\overline{e}^{*i_k}  =\\
        =       \sum_{i_1<\cdots< i_k,a,j}(-1)^{j+1} \frac{\partial
}{\partial x^a}\left(\zm_{i_1\dots i_k} f^{i_j}\right) {\dot
x}^a \overline{e}^{*i_1} \wedge \cdots \wedge \widehat{
\overline{e}^{*i_j}} \wedge \cdots \wedge \overline{e}^{*i_k} +\\
+       \sum_{i_1<\cdots< i_k,s\neq j}(-1)^{j+1} \zm_{i_1\dots i_k}
f^{i_j} \overline{e}^{*i_1} \wedge \cdots \wedge \widehat{
\overline{e}^{*i_j}} \wedge \cdots \wedge {\dot e}^{*i_s}\wedge
\dots \wedge \overline{e}^{*i_k} = \\
        =\T \left(\sum_{i_1<\cdots< i_k,j}(-1)^{j+1}\zm_{i_1\dots i_k} f^{i_j}
\overline{e}^{*i_1} \wedge \cdots \wedge \widehat{
\overline{e}^{*i_j}} \wedge \cdots \wedge \overline{e}^{*i_k}\right)
        = \T(\xi_X\zm).
                        \endmultline            $$

        (2)\ \ \  For $\zm  \in \zF^{k}(\zp)$, $\zn \in \zF^{}
(\zp)$, we have
                $$ \V(\xd_\zt(\zm \wedge \zn)) = \V(\xd_\zt( \zm)
\wedge \zn + (-1)^k \zm \wedge \xd_\zt (\zn)) =  \V(\xd_\zt( \zm))
\wedge \V(\zn) + (-1)^k \V (\zm) \wedge \V(\xd_\zt (\zn)).
                                        $$
        On the other hand,
                $$\xd_{\ssT \zt}(\V(\zm \wedge \zn)) = \xd_{\ssT
\zt}(\V(\zm) \wedge \V (\zn)) = \xd_{\ssT \zt}(\V(X))\wedge
\V(\zn) + (-1)^k \V(\zm) \wedge \xd_{\ssT \zt}(\V(\zn)).
                                        $$
        It follows, that in order to prove (a), it is enough to
consider the case  of $\zm $ being a function and $\zm = e^{*k}$.
For $\zm = f\in C^\infty(M)$ and $X\in \zF^{1}(\zt)$, we have,
using (1) and Theorem~10,
                $$ \multline
                \langle \V(\xd_\zt f), \T(X)\rangle = \xi_{\T(X)}
\V(\xd_\zt f) = \V(\xi_X \xd_\zt(f)) = \V(\za_\zt(X) (f)) = \dt
(\za_\zt (X)) (\vt (f)) = \\
        = \za_{\ssT \zt} (\T(X))(\V(f)) = \langle \xd_{\ssT \zt}
\V(f), \T(X)\rangle .
                                \endmultline    $$
         We have also
                $$\langle \V(\xd_\zt f), \V(X)\rangle = 0 =
\vt(\za_\zt(X)) (\vt(f)) = \za_{\ssT \zt}(\V(X)) (\V(f))=
\langle \xd_{\ssT \zt} \V(f), \V(X) \rangle  .
                                        $$
        Since $\V(X),\T(X)$ generate sections of the
$\sT\zt$-bundle, it follows that $\xd_{\ssT \zt} \V(f) =
\V(\xd_\zt f)$.
        Similarly, if $\zm = e^{*k}$, then
                $$\V(\xd_\zt(\zm)) = \V(\frac{1}{2}c^k_{ij} e^{*i}
\wedge e^{*j} ) = \frac{1}{2} c_{ij}^k \overline{e}^{*i} \wedge
\overline{e}^{*j} = \xd_{\ssT \zt}(\overline{e}^{*k}) =
\xd_{\ssT \zt} (\V(\zm)).
                                        $$
        To prove (b), we show first, using (a), that
                $$\multline
        \T(\xd_\zt(\zm\wedge \zn))= \T\left( \xd_\zt(\zm) \wedge
\zn + (-1)^k \zm \wedge \xd_\zt (\zn) \right)= \\
        = \T(\xd_\zt(\zm)) \wedge \V(\zn) + \V(\xd_\zt(\zm)) \wedge
\T(\zn) + (-1)^k \T(\zm) \wedge \V(\xd_\zt(\zn)) + (-1)^k \V(\zm)
\wedge \T(\xd_\zt(\zn)) = \\
        =  \T(\xd_\zt(\zm)) \wedge \V(\zn) + \xd_{\ssT \zt}
(\V(\zm)) \wedge \T(\zn) + (-1)^k \T(\zm) \wedge \xd_{\ssT \zt}
(\V(\zn)) + (-1)^k \V(\zm) \wedge \T(\xd_\zt(\zn)).
                                \endmultline    $$
         On the other hand,
                $$\multline
                \xd_{\ssT \zt}(\T(\zm \wedge \zn)) = \xd_{\ssT \zt}
\left( \T(\zm) \wedge \V(\zn) + \V(\zm) \wedge \T(\zn) \right) =
\\ = \xd_{\ssT \zt}(\T(\zm))\wedge \V(\zn) + (-1)^k \T(\zm)
\wedge \xd_{\ssT \zt}\V(\zn) + \xd_{\ssT \zt}(\V(\zm)) \wedge
\T(\zn) + (-1)^k\V(\zm) \wedge \xd_{\ssT \zt}(\T(\zn)).
                                \endmultline    $$
        This shows that it is enough to prove (b) for $\zm = f$ and
$\zm = e^{*k}$.
        Let us take $\zm = f\in C^\infty(M)$ and $X\in \zF^{1}(\zt)$. We have
                $$\multline
                \langle \T(\xd_{\zt} f), \T(X) \rangle  = \xi_{\T(X)}
\T(\xd_\zt f) = \dt (\xi_X\xd_\zt f) = \dt(\za_\zt (X)(f)) = \dt
(\za_\zt(X))(\dt f) = \\
        = \za_{\ssT \zt}(\T(X))(\dt f) = \langle \xd_{\ssT
\zt}(\T(f)), \T(X)\rangle
                                \endmultline    $$
and
                $$ \langle \T(\xd_{\zt} f), \V(X) \rangle = \xi_{\V(X)}
\T(\xd_\zt f) = \vt (\xi_X\xd_\zt f) = \za_{\ssT \zt}(\V(X))(\dt
f) =  \langle \xd_{\ssT \zt}(\T(f)), \V(X)\rangle.
                                                $$
        It follows that $ \T(\xd_\zt f) = \xd_{\ssT \zt} (\T(f))$.

        For $\zm = e^{*k}$, we have
                $$\T(\xd_\zt(e^{*k})) = \T(\frac{1}{2}c_{ji}^k e^{*i}
\wedge e^{*j} ) = c_{ji} ^k {\dot e}^{*i} \wedge
\overline{e}^{*j} + \frac{1}{2} \frac{\partial c_{ji}^k
}{\partial x^a} {\dot x}^a \overline{e}^{*i} \wedge
\overline{e}^{*j} = \xd_{\ssT \zt}({\dot e}^{*k}) = \xd_{\ssT
\zt}(\T(e^{*k})).
                                        $$

        (3)\ \ \  Since $\ll_X = \xi_X \circ \xd_\zt + \xd_\zt
\circ \xi_X$, we get, using (1) and (2), that
                $$ \multline
        \ll_{\V(X)}\T(\zm) =  \xi_{\V(X)} \xd_{\ss T\zt} \T(\zm) +
\xd_{\ssT \zt} \xi_{\V(X)} \T(\zm) = \xi_{\V(X)}\T(\xd_\zt\zm) +
\xd_{\ssT \zt} \V(\xi_X\zm) = \\ = \V(\xi_X  \xd_\zt \zm +
\xd_\zt \xi_X\zm) = \V(\ll_X \zm).
                \endmultline            $$
The rest of the proof is analogous.
        \endproof
        {\bf Remark.}   Sections $\T(X)$ span the $\sT \zt$-bundle
over each point except for the zero section of $\sT M$. It
follows that  in the proof above we can use the continuity
argument, instead of  considering  the case of sections $\V(X)$.

        \claim{Theorem}{}
        Let  $K,L \in \zF^{}_1(\zp)$. Then,
        \list
        \item"(a)" $[\V(K), \V(L)]^{N-R} = 0 $,
        \item"(b)" $[\V(K), \T(L)]^{N-R} = [\T(K), \V(L)]^{N-R} =
\V\left( [K, L]^{N-R} \right) $,
        \item"(c)" $[\T(K), \T(L)]^{N-R} = \T\left( [K, L]^{N-R}
\right) $, i.~e., the complete lift is a homomorphism between
Nijenhuis-Richardson brackets on $\zF^{}_1(\zp)$ and
$\zF^{}_1(\sT\zp)$.
        \endlist
        \endclaim
        \proof
        The theorem follows immediately from \Ref{1.11} and
Theorem~12~(1).
        Indeed, for $K = \zm\otimes X \in \zF^{k}_1(\zp)$, $L
=\zn\otimes Y$, we have
                $$ \multline
        [\V(K), \V(L)]^{N-R} = [\V(\zm)\otimes \V(X),
\V(\zn) \otimes \V(Y)]^{N-R} = \\ = \V(\zm)\wedge \xi_{\V(X)}\V(\zn)
\otimes \V(Y)  +(-1)^k \xi_{\V(Y)}\V(\zm) \wedge \V(\zn) \otimes
\V(X) =0 .
                                \endmultline            $$
Similarly,
                $$\multline
        [\V(K), \T(L)]^{N-R} = [\V(\zm)\otimes \V(X),
\T(\zn) \otimes \V(Y) + \V(\zn) \otimes \T(Y)]^{N-R} \\
        = \V(\zm)\wedge \xi_{\V(X)}\T(\zn) \otimes \V(Y) +
\V(\zm)\wedge \xi_{\V(X)}\V(\zn) \otimes \T(Y) + (-1)^k
\xi_{\V(Y)}\V(\zm) \wedge \T(\zn) \otimes \V(X)  \\ +(-1)^k
\xi_{\T(Y)}\V(\zm) \wedge \V(\zn) \otimes \V(X)
        = \V(\zm)\wedge \xi_{\V(X)}\T(\zn) \otimes \V(Y) + (-1)^k
\xi_{\T(Y)}\V(\zm) \wedge \V(\zn) \otimes \V(X)  \\ = \V\left(
\zm\wedge \xi_X \zn \otimes Y + (-1)^k \xi_Y \zm \wedge \zn
\otimes X\right) = \V([\zm\otimes X, \zn \otimes Y]^{N-R}).
                \endmultline                    $$
(c) can be proved by analogous calculations.
        \endproof

        We have also similar formulae for the Fr\"olicher-Nijenhuis
brackets.
        \claim{Theorem}{}
        Let  $K,L \in \zF^{}_1(\zp)$. Then,
        \list
        \item"(a)" $[\V(K), \V(L)]^{F-N} = 0 $,
        \item"(b)" $[\V(K), \T(L)]^{F-N} = [\T(K), \V(L)]^{F-N} =
\V\left( [K, L]^{F-N} \right) $,
        \item"(c)" $[\T(K), \T(L)]^{F-N} = \T\left( [K, L]^{F-N}
\right) $, i.~e., the complete lift $\T$ is a homomorphism between
Fr\"olicher-Nijenhuis brackets on $\zF^{}_1(\zp)$ and
$\zF^{}_1(\sT\zp)$.
        \endlist
        \endclaim
        \proof
        We make use of formulae of Theorem~11 and Theorem~12.
First, let us recall  that for $K =\zm\otimes X \in
\zF^{k}_1(\zp)$ and $L=\zn\otimes Y$, we have
                $$\multline
        [K,L]^{F-N} =\\ = \zm \wedge \zn \otimes [X,Y] + \zm
\wedge \ll_X\zn \otimes Y - \ll_Y\zm \wedge \zn \otimes X +
(-1)^k (\xd_\zt \zm \wedge \xi_X\zn \otimes Y + \xi_Y\zm \wedge
\xd_\zt \zn \otimes X),
                                \endmultline            $$
         so that
                $$\multline
                [\V(K), \V(L)]^{F-N} = \\ = \V(\zm) \wedge \V(\zn)
\otimes [\V(X), \V(Y)] + \V(\zm) \wedge \ll_{\V(X)}\V(\zn)
\otimes \V(Y) - \ll_{\V(Y)}\V(\zm) \wedge \V(\zn) \otimes \V(X) +\\
+ (-1)^k (\xd_{\ssT \zt} \V(\zm) \wedge \xi_{\V(X)} \V(\zn)
\otimes \V(Y) + \xi_{\V(Y)} \V(\zm) \wedge \xd_{\ssT \zt}
\V(\zn) \otimes \V(X)) =0.
                                \endmultline            $$
        Similarly,
                $$\multline
                [\V(K), \T(L)]^{F-N} = [\V(\zm) \otimes \V(X),
\T(\zn)\otimes \V(Y) + \V(\zn) \otimes \T(Y)]^{F-N} = \\
        = \V(\zm) \wedge \T(\zn) \otimes [\V(X),\V(Y)] + \V(\zm)
\wedge \ll_{\V(X)} \T(\zn)\otimes \V(Y) - \ll_{\V(Y)} \V(\zm)
\wedge \T(\zn)\otimes \V(X) +  \\ + (-1)^k\left( \xd_{\ssT\zt}\V(\zm)
\wedge \xi_{\V(X)} \T(\zn) \otimes \V(Y) + \xi_{\V(Y)} \V(\zm)
\wedge \xd_{\ssT\zt} \T(\zn) \otimes \V(X) \right) +\\
        + \V(\zm) \wedge \V(\zn) \otimes [\V(X),\T(Y)] + \V(\zm)
\wedge \ll_{\V(X)} \V(\zn)\otimes \T(Y) - \ll_{\T(Y)} \V(\zm)
\wedge \V(\zn)\otimes \V(X) + \\
        + (-1)^k\left( \xd_{\ssT\zt}V(\zm)
\wedge \xi_{\V(X)} \V(\zn) \otimes \T(Y)  + \xi_{\T(Y)} \V(\zm)
\wedge \xd_{\ssT\zt} \V(\zn) \otimes \V(X) \right) = \\
        = \V(\zm) \wedge \V(\ll_X \zn) \otimes \V(Y) + (-1)^k \V
(\xd_\zt \zm) \wedge \V(\xi_X\zn) \otimes \V(Y) + \V(\zm) \wedge
\V(\zn) \otimes \V([X,Y]) -\\
        - \V(\ll_Y\zm ) \wedge \V(\zn) \otimes \V(X) + (-1)^k
\V(\xi_Y\zm) \wedge \V(\xd_{\ssT \zt}\zn) \otimes \V(X) = \\
        = \V\left(\zm \wedge \zn \otimes [X,Y] + \zm \wedge \ll_X
\zn \otimes Y - \ll_Y\zm \wedge \zn \otimes X + \right. \\
        \left. +(-1)^k (\xd_\zt \zm \wedge \xi_X \zn \otimes Y +
\xi_Y \zm \wedge \zn \otimes X) \right) = \V([K,L]^{F-N}).
                                \endmultline    $$
        The remaining part of the proof consists of analogous easy
calculations and we skip them.
\       \endproof

        \hl1{Cotangent lifts and graded Lie brackets}
        \newpre\tag{6.} \Reset\tag{1}
        We shall study cotangent lifts of tensor fields associated
with a Lie algebroid. They are, however, defined only for
certain types of tensor fields and they do not have so nice
functorial properties as the tangent lifts. On the other hand,
they can be very useful, as we shall see later.
        Let $\zt \colon E \rightarrow M$ be a vector bundle with a
Lie algebroid structure, let $\za_\zt \colon E\rightarrow \sT M$
be the anchor map of this structure and let  $P$ be the
corresponding linear Poisson structure on the dual bundle $\zp
\colon E^\* \rightarrow M$. We define the (cotangent) {\it
complete lift}
                $$\H \colon \zF^{1}(\zt) \rightarrow \zF^{1}(\ezt_{E^\s*})
                                                $$
by the formula
                $$ \H(X) = \widetilde{P}( \xd \zi(X)) =\sH_P(\zi(X)) =
-[P, \zi(X)].
                                                $$
        In  the local coordinates introduced in Section~4, we have
                $$ \H(f^ie_i) = \left( f^i c^k_{ij} \zx_k -
\frac{\partial f^i}{ \partial x^a} \zd^a_j \zx_i \right)\partial
_{\zx_j} + f^i \zd^a_i \partial _{x^a}.
                                                \tag \label{6.1}$$

        For $\zm\in \zF^{k}(\zp)$, we have  the  {\it vertical lift}
$\V_\zp(\zm) \in \zF^{k}(\ezt_{E^\s*})$ (compare with \Ref{4.7}). In local
coordinates,
                $$ \V_\zp
(\sum_{i_1<\cdots<i_k}\zm_{i_1\dots i_k}e^{*i_1} \wedge \cdots
\wedge e^{*i_k})  =   \sum_{i_1<\cdots<i_k}\zm_{i_1\dots i_k}
\partial _{\zx_{i_1}} \wedge \cdots \wedge \partial _{\zx_{i_k}}.
                                        \tag\label{6.2}$$

        \claim{Theorem}{}
        Let  $\zm, \zn \in \zF^{}(\zp)$ and let   $Y, X \in
\zF^{1}(\zt)$. Then,
        \list
        \item"(a)" $\V_\zp(\zm \wedge \zn) = \V_\zp(\zm) \wedge
\V_\zp(\zn)$;
        \item"(b)" $[\V_\zp(\zm), \V_\zp(\zn)] =0$;
        \item"(c)" $[\zi(X), \V_\zp(\zm)] = -\V_\zp(\xi_X\zm)$;
        \item"(d)" $[P, \V_\zp(\zm)] = \V_\zp(\xd_\zt\zm)$;
        \item"(e)" $[\H(X), \V_\zp(\zm)] = \V_\zp(\ll_X\zm)$;
        \item"(f)" $[\H(X), \H(Y)] = \H([X,Y]) $;
        \item"(g)" $[\H(X), \zi(Y)] = \zi([X,Y])$,
        \endlist
where brackets on the left-hand sides are classical Schouten
brackets and brackets on the right-hand sides are Lie algebroid
brackets.
        \endclaim
        \proof
(a) and (b) follow immediately from \Ref{6.2}.

        (c)\ \ \  Since $\ad_{\zi(X)}$ and $\xi_X$ are graded
derivations of degree  -1 and $\V_\zp$ is, in view of (a), a
homomorphism, both sides of (c) are graded derivations of degree
-1 with respect to $\zm$. It follows by  induction that it is
enough to consider the case of $\zm = e^{*k}$. We have for $\zm
= e^{*k}$, $X= f^ie_i$,
                $$[\zi(X), \V_\zp(e^{*k})] = [f^i\zx_i, \partial
_{\zx_k}] = -f^k = -\V_\zp(\xi_X e^{*k}).
                                        $$

        (d)\ \ \ Similarly as in (c), we observe first that both
sides of (d) are graded derivations of degree 1 with respect to
$\zm$. It follows then that it is enough to prove (d) for $\zm
\in C^\infty(M)$ and $\zm = e^{*k}$. Let  $X\in \zF^{1}(\zt)$.
The formulae \Ref{3.2} and (c) imply that
                $$\multline
                \xi_{\xd\zi(Y)} [P,\V_\zp(f)] = [[P, \zp^\*f], \zi(Y)]
= -[[P,\zi(Y)], \zp^\*f] = [\H(Y), \zp^\*f] =\zp^\*
(\za_\zt(Y)(f)) =\\
        = \V_\zp(\xi_Y \xd_\zt f) = - [ \zi(Y), \V_\zp (\xd_\zt f)] =
\xi_{\xd\zi(Y)} \V_\zp(\xd_\zt f).
                                \endmultline    $$
        We conclude that  $[P ,\V_\zp(f)] = \V_\zp(\xd_\zt f)$.

        For $\zm = e^{*k}$, we have in local coordinates (with $P$ as in
\Ref{3.4})
                $$\multline
                [P, \V_\zp(e^{*k})] = [P, \partial _{\zx_k}] = [
\frac{1}{2} c^k_{ij} \zx_k \partial _{\zx_i} \wedge \partial
_{\zx_j} + \zd^a_i \partial _{\zx_i} \wedge \partial _{x^a},
\partial _{\zx_k}] = \frac{1}{2} c_{ji}^k \partial _{\zx_i}
\wedge \partial _{\zx_j} .
                                \endmultline    $$
        On the other hand,
                $$\V_\zp(\xd_\zt(e^{*k})) = \V_\zp (\frac{1}{2}
c_{ji}^k e^{*i} \wedge e^{*j} ) = \frac{1}{2} c_{ji}^k \partial
_{\zx_i} \wedge \partial _{\zx_j}
                                                $$
        and (d) follows.

        (e)\ \ \ According to (c) and (d), we have
                $$\multline
                \V_\zp(\ll_X\zm) =  \V_\zp (\xd_\zt \xi _X \zm) +
\V_\zp(\xi_X \xd_\zt \zm) = [P, \V_\zp(\xi_X \zm)] - [\zi(X),
\V_\zp (\xd_\zt \zm)] = \\
        = - [P,[\zi(X), \V_\zp(\zm)]] - [\zi(X),[P, \V_\zp (\zm)]] =
-[[P,\zi(X)], \V_\zp(\zm)] = [\H(X), \V_\zp(\zm)].
                                \endmultline    $$

        (f)\ \ \  Using the graded Jacobi identity and the formula
$ [P,P] =0$, we get
                $$ \multline
        [\H(X), \H(Y)] = [[P,\zi(X)], [P, \zi(Y)]] = [P,
[[P,\zi(X)], \zi(Y)]] = [P, \{ \zi(Y), \zi(X) \}_P] = \\
        =- [P, \zi([X,Y])] = \H([X,Y]).
                                \endmultline    $$

        (g)\ \ \  We have
                $$ [\H(X), \zi(Y)] = [-[P, \zi(X)], \zi(Y)] = \{
\zi(X), \zi(Y)\}_P = \zi([X,Y]).
                                        $$
        \endproof

        \claim{Theorem}{}
The formulae
                $$ \I(\zm \otimes X) = - \zi(X)\V_\zp(\zm)
                                        \tag \label{6.3}$$
        and
                $$ \H(\zm \otimes X) = [P,\I(\zm \otimes X) ] = \H(X)
\wedge \V_\zp(\zm) - \zi(X) \V_\zp (\xd_\zt \zm),
                                        \tag \label{6.4}$$
        for simple tensors $\zm \otimes X \in \zF^{}_1(\zp)$, define
linear mappings
                $$\I,\, \H \colon \zF^{}_1(\zp) \rightarrow
\zF^{}(\ezt_{E^*})
                                        $$
        which coincide on $\zF^{0}_1(\zp)$ with already defined
$-\zi$ and $\H$ respectively.
        \endclaim
        \proof
        The right-hand side of \Ref{6.3} is linear with respect to
$X$ and $\zm$ and for $f\in C^\infty (M)$ we have
                $$\I((f\zm) \otimes X) = -\zi (X) \V_\zp(f\zm) = -\zp^\*
(f) \zi(X) \V_\zp(\zm) = - \zi(fX) \V_\zp(\zm) = \I(\zm \otimes (fX)).
                                        $$
        Consequently, $\I$ is well defined. Since (6.4) is the
composition of $\I$ and the Schouten bracket with $P$,  also $\H$
is well defined. Of course, we have, for $X\in \zF^{1}(\zt)$,
                $$ \I(X) = -\zi(X), \ \ \ \ \H(X) = -[P,\zi(X)].
                                        $$
         \endproof

        We call  $\H(K)$ the {\it dual complete lift} of $K\in
\zF^{}_1(\zp)$. We shall show that $\I$ and $\H$ are
homomorphisms of gradad Lie algebras.

        \claim{Theorem}{}
        The mapping
                $$\I \colon  \zF^{}_1(\zp) \rightarrow  \zF^{}(\ezt_{E^*})
                                        $$
        is an injective homomorphism of the Nijenhuis-Richardson
bracket into the Schouten bra\-cket:
                $$   \I([K, L]^{N-R}) = [\I(K), \I(L)].
                                        $$
In particular, in the case of the canonical Lie algebroid, we get an
        embedding of the Nijenhuis-Richardson bracket on $M$ into the
        Schouten bracket on $\sT^\* M$.

        \endclaim

        \proof
        Let  $K= \zm \otimes X \in \zF^{k}_{1}(\zp)$ and $L
=\zn\otimes Y \in \zF^{}_{1}(\zp)$. It follows from Theorem~15
((a), (b),(c)) and from \Ref{1.5} that
                $$\multline
                [\I(K), \I(L)] = [\zi(X) \V_\zp(\zm), \zi(Y)
\V_\zp(\zn)] = [\zi(X) \V_\zp(\zm), \zi(Y)] \wedge \V_\zp(\zn) +\\
        + \zi(Y) [\zi(X) \V_\zp(\zm), \V_\zp(\zn)] =
         (-1)^{k+1} \zi(X)\V_\zp(\xi_Y \zm) \wedge \V_\zp (\zn) -
\zi(Y) \V_\zp(\zm) \wedge \V_\zp(\xi_X \zn) = \\
        - \left(\zi(X)\V_\zp (\zm \wedge \xi_X \zn) + (-1)^k \zi(X)
\V_\zp(\xi_Y \zm \wedge \zn) \right) = \\
        = \I\left(\zm \wedge \xi_X \zn \otimes Y + (-1)^k \xi_Y \zm
\wedge \zn \otimes X\right) = \I([K,L]^{N-R}).
                \endmultline            $$

        The injectivity follows from the local form of $\I$ on
$\zF^{k}_{1}(\zp)$:
                $$\I\left( \sum_{i_1<\cdots <i_k,j} \zm_{i_1\dots
i_k}^j(x) e^{*i_1} \wedge  \cdots \wedge e^{*i_k} \otimes e_j
\right) = - \sum_{i_1<\cdots <i_k, j} \zm_{i_1\dots i_k}^j(x)
\zx_j \partial _{\zx_{i_1}} \wedge \cdots \wedge  \partial
_{\zx_{i_k}}.
                                        $$
         It is evident that $\I(K) =0$ if and only if $K=0$.
        \endproof

        \claim{Theorem}{}
        The dual complete lift
                $$ \H \colon \zF^{}_{1}(\zp) \rightarrow
\zF^{}_{}(\ezt_{E^*})
                                        $$
         is a homomorphism of the Fr\"olicher--Nijenhuis bracket
into the Schouten bracket:
                $$ \H([K,L]^{F-N}) = [ \H(K), \H(L)] .
                                        $$
        \endclaim
        \proof
        Let  $K= \zm\otimes X \in \zF^{k}_{1}(\zp)$ and $L = \zn
\otimes Y \in \zF^{}_{1}(\zp) $. Using the graded Jacobi identity
and Theorem~15, we get three identities:
        $$\multline
        [\H(X)\wedge \V_\zp(\zm), \H(Y)\wedge \V_\zp(\zn)]
=[\H(X), \H(Y)] \wedge \V_\zp(\zm) \wedge \V_\zp(\zn) +\\
        +\H(X) \wedge [\V_\zp(\zm), \H(Y)] \wedge \V_\zp(\zn) +
\H(Y) \wedge \V_\zp(\zm) \wedge [\H(X), \V_\zp(\zn)] = \\
        = \H([X,Y]) \wedge \V_\zp(\zm \wedge \zn) - \H(X) \wedge
\V_\zp (\ll_Y\zm \wedge \zn) + \H(Y) \wedge \V_\zp(\zm \wedge
\ll_X \zn),
                        \endmultline \tag "(1)" $$

        $$\multline
        [\H(X)\wedge \V_\zp(\zm), \zi(Y) \V_\zp(\xd_\zt \zn)] =\\
        = (-1)^k [\H(X), \zi(Y)] \wedge \V_\zp(\zm) \wedge \V_\zp
(\xd_\zt\zn) +  \H(X) \wedge [\V_\zp(\zm), \zi(Y)] \wedge \V_\zp
(\xd_\zt \zn) + \\
        +(-1)^k \zi(Y) \V_\zp (\zm) \wedge [\H(X), \V_\zp(
\xd_\zt\zn)] = \\
        = (-1)^k \zi([X,Y]) \V_\zp (\zm\wedge \xd_\zt\zn) - (-1)^k
\H(X) \wedge \V_\zp(\xi_Y\zm \wedge \xd_\zt\zn) + (-1)^k \zi(Y)
\V_\zp (\zm \wedge \ll_X \xd_\zt\zn),
                \endmultline   \tag "(2)"               $$

        $$\multline
        [\zi(X)\V_\zp(\xd_\zt \zm), \zi(Y) \V_\zp(\xd_\zt \zn)] =\\
        =\zi(X) [\V_\zp(\xd_\zt \zm), \zi(Y)] \wedge
\V_\zp(\xd_\zt \zn) + \zi(Y) \V_\zp(\xd_\zt \zm) \wedge [\zi(X),
\V_\zp(\xd_\zt \zn)] = \\
        = (-1)^k \zi(X) \V_\zp(\xi_Y \xd_\zt \zm \wedge \xd_\zt
\zn) - \zi(Y) \V_\zp( \xd_\zt \zm \wedge \xi_X\xd_\zt \zn).
                        \endmultline    \tag "(3)"      $$

        Hence,
        $$\multline
        [\H(K), \H(L)] = [\H(X) \wedge \V_\zp(\zm) - \zi(X) \V_\zp
(\xd_\zt\zm), \H(Y) \wedge \V_\zp(\zn) - \zi(Y) \V_\zp
(\xd_\zt\zn)] = \\
        = \H([X,Y]) \wedge \V_\zp(\zm \wedge \zn) - \zi([X,Y])
\V_\zp(\xd_\zt\zm \wedge \zn + (-1)^k \zm \wedge \xd_\zt\zn) +
\H(Y) \wedge \V_\zp(\zm \wedge \ll_X \zn) - \\
        - \zi(Y) \V_\zp(\xd_\zt \zm \wedge \xi_X \xd_\zt \zn +
(-1)^k \zm \wedge \ll_X \xd_\zt\zn) - \H(X) \wedge \V_\zp (\ll_Y
\zm \wedge \zn) + \\
        +\zi(X) \V_\zp(\ll_Y \xd_\zt\zm \wedge \zn
+(-1)^k \xi_Y \xd_\zt \zm \wedge \xd_\zt \zn)
        + (-1)^k \H(Y)\wedge \V_\zp(\xd_\zt \zm \wedge \xi_X\zn) + \\
        +(-1)^k \H(X) \wedge \V_\zp(\xi_Y\zm \wedge \xd_\zt \zn) =
\H(\zm \wedge \zn \otimes [X,Y]) + \H(\zm \wedge \ll_X \zn
\otimes Y) + \\
        + \zi(Y) \V_\zp(\xd_\zt \zm \wedge \xd_\zt \xi_X \zn) -
\H(\ll_Y\zm \wedge \zn \otimes X)  - (-1)^k \zi(X)
\V_\zp(\xd_\zt \zi_Y\zm \wedge \xd_\zt \zn) +\\
        + (-1)^k \H(Y)\wedge \V_\zp(\xd_\zt \zm \wedge \xi_X\zn) +
(-1)^k \H(X) \wedge \V_\zp(\xi_Y\zm \wedge \xd_\zt \zn) = \\
        = \H(\zm \wedge \zn \otimes [X,Y]) + \H(\zm \wedge \ll_X
\zn \otimes Y) - \H(\ll_Y\zm \wedge \zn \otimes X) +(-1)^k
\H(\xd_\zt \zm \wedge \xi_X \zn \otimes Y) +\\
        +(-1)^k \H(\xi_Y \zm \wedge \xd_\zt \zn \otimes X)  =
\H([K,L]^{F-N}).
                        \endmultline$$

        \endproof
                \hl1{The  case of the canonical Lie algebroid}
        \newpre\tag{7.} \Reset\tag1
        It is obvious that the canonical Lie algebroid, including the
standard differential calculus on  a manifold, is of special
interest. It is exactly the case, when the corresponding linear
Poisson on the dual bundle is symplectic.
        Let us suppose that the linear Poisson structure $P$ on
$E^\*$ is nondegenerate, i.~e., the mapping $\widetilde{P}\colon \sT^\*
E^\* \rightarrow \sT E^\*$ is an isomorphism of vector bundles.
The inverse morphism $\widetilde{\zw} = \widetilde{P}^{-1}$
defines a symplectic form $\zw$ on $E^\*$:
                $$ \langle P, \zm \wedge \zn \rangle = \langle
\widetilde{P}(\zm) \wedge \widetilde{P}(\zn), \zw \rangle.
                                        $$
         We have in local coordinates (as in \Ref{3.4})
                $$ P= \frac{1}{2}c_{ij}^k(x) \zx_k\partial _{\zx_i}\wedge
\partial _{\zx_j}  + \zd_i^a(x) \partial _{\zx_i}\wedge \partial _{x_a}.
                                        $$
        Since $P$ is nonegenerate, the anchor map   $\za_\zt
\colon E \rightarrow \sT M$, given by
                $$\za_\zt(y^ie_i(x)) = y^i\zd_i^a(x)\partial _{x^a},$$
is an isomorphism of vector bundles. The anchor map is a
homomorphism of Lie brackets and it follows that  $\za_\zt$ defines
an isomorphism between the Lie algebroid structure on $E$ and the
canonical Lie algebroid structure on $\sT M$. The dual (Poisson)
bundle $(E^\*, P)$ can be identified with $(\sT^\*M, P_M)$ --
the cotangent bundle with the canonical Poisson structure. In the
following, we consider the case of canonical structures only.

        Let $(x^a)$ be a local coordinate system on $M$. We choose
$(e_a = \partial _{x^a})$ and $(e^{*a} = \xd x^a)$ as local
bases of sections of $\ezt_M \colon \sT M\rightarrow M$ and
$\ezp_M \colon \sT ^\* M\rightarrow M$, respectively.   We
denote by $(x^a, {\dot x}^b)$ and $(x^a, p_b)$, instead of
$(x^a, y^i)$  and $(x^a, \zx_i)$, the adopted coordinate systems.
In these coordinate systems, the canonical Poisson bivector
field is given by
                $$  P= \partial _{p_a}\wedge \partial _{x^a},
                                        \tag \label{7.1}$$
the anchor map $\za_{\ezt_M} \colon \sT M \rightarrow \sT M$ is
the identity, and $\xd_{\ezt_M} = \xd \,$ is the usual exterior
derivative. The formulae \Ref{3.6} and \Ref{3.8} take now the
form
                $$ [\partial _{x^a}, \partial _{x^b}] =0,\ \ \ \xd\,(f) =
\frac{\partial f}{\partial x^a}\xd x^a, \ \ \ \xd\,(\xd x^a) = 0.
                                        \tag \label{7.2}$$

        As in the previous sections, we can consider the tangent and
cotangent lifts of the canonical Lie algebroid structure on $\ezt_M
\colon \sT M \rightarrow M$ to Lie algebroid structures of bundles
$\sT \ezt_M \colon  \sT\sT M\rightarrow \sT M$ and $
\ezp_{\ssT^\s*} \sT^\* \sT^\* M \rightarrow \sT^\* M$. On the
other hand, there are canonical Lie algebroid structures of bundles
$\ezt_{\ssT M} \colon \sT\sT M\rightarrow \sT M$ and
$\ezt_{\ssT^\s* M} \colon \sT \sT^\* M\rightarrow \sT^\* M$.
In the following, we discuss relations between these Lie algebroids.

        There is the well-known isomorphism $\ezk_M$ of vector
bundles (cf. [30], [10]),
                $$ \CD \sT\sT M @() \L{\sT \ezt_M} @(0,-1) & \sT\sT M
@() \L{ \ezk_M} @(-1,0) @() \l{\ezt_{\ssT M}}@(0,-1 )\\
        \sT M  @()\a= @(1,0) & \sT M.
                                        \endCD          \tag \label{7.3}$$
As a mapping of manifolds, $\ezk_M$ is an involution: $\ezk_M \circ
\ezk_M = \idt$. As in Section~\Ref{3}, we introduce the basis
$(\overline{\partial _{x^a}},\dot{\partial _{x^b}} )$  of local
sections of $\sT\ezt_M$. The morphism $\ezk_M$ is
characterized by
                $$ \ezk_M \circ \partial _{{\dot x}^a} =
\overline{\partial _{x^a}}, \ \  \ezk_M\circ \partial _{x^b} =
\dot{\partial _{x^b}}
                                        \tag \label{7.4}$$

        and  it can be extended to  isomorphisms of tensor bundles
                $$ \ezk^r_M \colon  \otimes ^r_{\ezt_{\smT M}} \sT\sT M
\longrightarrow \otimes ^r_{\ssT \ezt_M} \sT\sT M
                                        $$
by
                $$\ezk^r_M (X_1\otimes_{\ezt_{\smT M}}\cdots \otimes
_{\ezt_{\smT M}} X_r ) = \ezk_M(X_1) \otimes _{\ssT \ezt_M}  \cdots
\otimes _{\ssT \ezt_M} \ezk_M(X_r).
                                        $$
        It follows that $\ezk^r_M$ can be restricted to skew-symmetric
and symmetric tensors:
                $$\ezk^r_M (\bigwedge^r_{\ezt_{\smT M}} \sT \sT M) =
\bigwedge^r _{\ssT \ezt_M} \sT\sT M , \ \ \  \ezk^r_M (\ss
S^r_{\ezt_{\smT M}}  \sT \sT M) = \ss S^r_{\ssT \ezt_M} \sT\sT M.
                                        $$
          For dual vector bundles $\sT \ezp_M \colon \sT \sT^\* M
\rightarrow \sT M$ and $\ezp_{\ssT M}\colon  \sT^\* \sT M
\rightarrow \sT M$, we have the dual isomorphism $\eza_M$
                $$ \CD
        \sT \sT^\* M @() \L{\sT\ezp_M } @(0,-1) @() \L{\eza_M}
@(1,0) & \sT^\* \sT M @() \l{\ezp_{\smT M}} @(0,-1) \\
        \sT M @()\a= @(1,0) & \sT M
                        \endCD          \tag \label{7.5}$$
and, consequently, isomorphisms of tensor bundles
                $$ \eza_M^r \colon \otimes ^r_{\ssT \ezp_M}  \sT\sT^\*
M \rightarrow \otimes ^r _{\ezp_{\smT M}} \sT^\* \sT M ,\ \ \
\eza^r_M \colon \bigwedge ^r_{\ssT \ezp_M} \sT\sT^\* M
\rightarrow \bigwedge ^r_{\ezp_{\smT M}} \sT^\* \sT M.
                                        \tag \label{7.6}$$
The induced mapping $\eza_M \colon \zF^{1}_{}(\sT\ezp_M)
\rightarrow \zF^{1}_{}(\ezp_{\ssT M})$ is given by
                $$  \eza_M(\overline{\xd p}_a) = \xd x^a,\ \ \  \eza_M
(\dot{\xd p}_b) = \xd {\dot x}^b.
                                        \tag \label{7.7}$$

        There are canonical lifts $\vt$ (vertical) and $\dt$
(complete) of multivector fields and forms on $M$ to multivector
fields and forms on $\sT M$ (cf. [31], [10] and \Ref{4.1},
\Ref{4.2}), which are strictly related to $\V$ and $\T$.
                \claim{Theorem}{}
        \list
                \item For $X\in \zG(M, \otimes ^r_{\ezt_M})$, we have
                $$ \ezk^r_M(\vt(X)) = \V(X), \ \ \ \ezk^r_M(\dt(X)) = \T(X).
                                        \tag \label{7.8}$$
                \item For $\zm \in \zG(M, \otimes ^r_{\ezp_M})$, we have
                $$ \vt(\zm) = \eza^r_M(\V(\zm)), \ \ \ \dt(\zm) =
\eza^r_M (\T(\zm)).
                                        \tag \label{7.9}$$
        \endlist
There are analogous formulae for skew-symmetric and symmetric
tensors, as well as for tensors of type $(r,s)$.
        \endclaim
        \proof
         Since
                $$ \vt( X \otimes Y) = \vt(X) \otimes \vt(Y),\ \
\dt(X \otimes Y) = \dt(X) \otimes \vt(Y) + \vt(X) \otimes \dt(Y)
                                                $$
(cf. [30], [10]) and
                $$ \V( X \otimes Y) = \V(X) \otimes \V(Y),\ \
\T(X \otimes Y) = \T(X) \otimes \V(Y) + \V(X) \otimes \T(Y)
                                                $$
(Theorem~9), it is sufficient to consider the case $r=1$.

        Let  $X = f^a \partial _{x^a}$. Then (\Ref{4.6} and
\Ref{4.7}), $\V(X) = f^a \overline{ \partial _{x^a}}$ and $\T(X)
= f^a \dot{\partial _{x^a}} + \dfrac{\partial f^a}{ \partial x^b}
{\dot x}^b \overline{\partial _{x^a}}$. It follows from \Ref{7.4}
that
                $$\split
                \T(X) = & \ezk_M(f^a\partial _{x_a} + \frac{\partial
f^a}{ \partial x^b} {\dot x}^b \partial _{{\dot x}^a}) =
\ezk_M(\dt(X)),\\
                \V(X) = & \ezk_M(f^a \partial _{{\dot x}^a}) = \vt(X).
                        \endsplit               $$
        Similar arguments prove \Ref{7.9}.
        \endproof

        Let us remark that the above theorem can be proved in a
coordinate-free way. First, we make use of the following
identities (cf. [10]):
                $$ \split
        \dt\widetilde{X} = & \widetilde{\dt X}\circ \eza_M, \\
        \vt\widetilde{X} = & \widetilde{\vt X} \circ \eza_M, \\
        \widetilde{\dt \zm} = & \dt \widetilde{\zm}\circ \ezk_M, \\
        \widetilde{\vt \zm} = & \vt \widetilde{\zm}\circ \ezk_M.
                        \endsplit               \tag \label{7.10}$$
        Then, using the pairings
                $$\langle \,,\,\rangle \colon \sT \sT M \times _{\ssT
M} \sT^\* \sT M \rightarrow \R,\ \  \langle \,,\,\rangle' \colon
\sT \sT M \times _{\ssT M} \sT \sT^\* M \rightarrow \R
                                                $$
        and the equality
                $$\langle v,\eza_M(w)\rangle = \langle \ezk_M(v),w
\rangle ' ,
                                \tag\label{7.11}        $$
        we get
                $$      \multline
        \langle \V(X), w \rangle '= \widetilde{\V(X)}(w) =
\vt(\widetilde{X})(w) = \widetilde{\vt(X)}( \eza_M(w))=\\
        = \langle \vt(X), \eza_M(w)\rangle = \langle \ezk_M(\vt(X)),
w\rangle' .
                        \endmultline            $$
        Hence, $\ezk_M(\V(X)) = \vt(X)  $, etc.

        \claim{Theorem}{}
        The isomorphism $\ezk_M$ of vector bundles is an isomorphism
of the canonical Lie algebroid structure on $\ezt_{\ssT M}
\colon \sT\sT M \rightarrow \sT M  $ and the tangent Lie algebroid
structure on $\sT \ezt_M \colon  \sT \sT M\rightarrow \sT M$

        The dual isomorphism $\eza_M$ of vector bundles is an
isomorphism of linear Poisson structures $(\sT \sT^\* M, \dt P_M)$
and $(\sT^\* \sT M, P_{\ssT M})$. Moreover,
                $$ \eza_M \circ  \xd_{\ssT \ezt_M} = \xd \circ \eza_M.
                                                \tag \label{7.13}$$
        \endclaim
        \proof
        For a section $X\colon M\rightarrow E$, we have $ \T(X) =
\sT X$. It follows that
                $$ \T(\za_\zt\circ X) =\sT (\za_\zt\circ X)= \sT\za_\zt
\circ \T(X).
                                                $$
        We have then, in view of Theorem~19 and Theorem~10, that
                $$ \sT\za_\zt\circ \T(X) = \T(\za_\zt\circ X) = \ezk_M
(\dt (\za_\zt\circ X )) = \ezk_M \circ \za_{\ssT \zt}.
                                                $$
        Since sections $\T(X)$ span the $\sT\zt$-bundle (except for the
zero sections), we obtain, by the continuity argument, that, for a
Lie algebroid $\zt \colon E\rightarrow M$ with the anchor $\za_\zt$,
we have  (cf. [22])
                $$ \ezk_M \circ \za_{\ssT \zt} = \sT \za_\zt.
                                                \tag \label{7.12}$$

        Since the anchor of the canonical Lie algebroid structure is the
identity on $\sT M$, we get from \Ref{7.12} that the anchor
$\za_{\ssT \ezt_M}$ of the tangent Lie algebroid structure is
$\ezk_M^{-1}$. As in the begining of this section, we obtain that
$\ezk_M$ is an isomorphism of Lie algebroids.

        The remaining part of the theorem follows directly from
\Ref{7.11}.
        \endproof
        \claim{Corollary}{}
        Let  $X,Y\in \zF^{}_{}(\ezt_M)$ or $X,Y \in
\zF^{}_{1}(\ezp_M)$ and let $\lb\,,\,\rb$ stays for the Schouten or
Nijenhuijs-Richardson or Fr\"olicher-Nijenhuis bracket. Then,
                $$ \split
                \lb \vt(X), \vt(Y) \rb &= 0;\\
                \lb \vt(X), \dt(Y) \rb &= \lb \dt(X), \vt(Y) \rb =
\vt\, \lb X,Y\rb;\\
                \lb \dt(X), \dt(Y) \rb &= \, \dt \lb X,Y\rb .
                        \endsplit               \tag \label{7.14}$$
        We have also all the formulae of Theorem~12 with $\V$
replaced by $\vt$, $\T$ replaced by $\dt$, and $\xd_{\ssT \zt}$,
$\xd_\zt$ replaced by $\xd\,$.
        \endclaim
        \proof
        All the formulae make use of the Lie algebroid constructions
only. We have already proved analogous formulae for lifts $\V,\
\T$, and the bracket of the tangent Lie algebroid  (Theorems
11, 13, 14). The formulae \Ref{7.14} are then simple consequences
of Theorem 20 and Theorem 19.
        \endproof

        Since $P_M$ is invertible, it follows from Theorem 5 that
the mapping
                $$\ss \zL_{P_M} \colon \zF^{}_{}(\ezp_{\ssT^\s* M})
\rightarrow \zF^{}_{}(\ezt_{\ssT^\s* M})
                                                $$
         is an isomorphism between the Koszul-Schouten and the
Schouten brackets. Since $P_M$ is skew-symmetric, for the  dual
map ${\ss \zL}^\*_{P_M}$, we have ${\ss \zL}^\*_{P_M}
(\zm) = (-1)^k \ss \zL_{P_M}(\zm)$ with $\zm \in
\zF^{k}_{}(\ezp_{\ssT^\s* M})$. Therefore it defines an
antihomomorphism of these brackets.
        We have also the well-known formula
                $$ {\ss \zL}^\*_{P_M} (\xd\zm) = [P, {\ss \zL}^\*_{P_M}
(\zm)].
                                                \tag \label{7.15}$$
        \claim{Theorem}{}
        For $\zm\in \zF^{k}_{}(\ezp_M)$, $X\in \zF^{1}_{}(\ezt_M)$,
we have
        \list
        \item ${\ss \zL}^\*_{P_M} (\ezp_M^\*(\zm)) = \V_{\ezp_M}(\zm)$,
        \item ${\ss \zL}^\*_{P_M} (\zi(X)\ezp_M^\*(\zm))= -\I(\zm\otimes X)$,
        \item ${\ss \zL}^\*_{P_M} (\xd(\zi(X)\ezp_M^\*(\zm)))= -\H(\zm\otimes
X)$
        \endlist
        \endclaim
        \proof
        \list
        \item Let  $\zm = \zm_1\wedge \dots \wedge \zm_k$, where $\zm_i
\in \zF^{1}_{}(\ezp_M)$. Since
                $${\ss \zL}^\*_{P_M} (\ezp_M^\*(\zm)) = (-1)^{k}
\widetilde{P_M}(\ezp_M^\*(\zm_1)) \wedge \cdots \wedge
\widetilde{P_M}(\ezp_M^\*(\zm_k))
                                                $$
        and
                $$
\V_{\ezp_M}(\zm) = \V_{\ezp_M}(\zm_1) \wedge \cdots \wedge
\V_{\ezp_M}(\zm_k),
                                                $$
        it is sufficient to prove (1) in the case $k=1$. For $\zm =
\zm_a \xd x^a$, we have
                $$ \widetilde{P_M}(\ezp_M^\*(\zm_a \xd x^a)) =
\widetilde{P_M}(\zm_a \xd x^a) = -\zm_a \partial _{p_a}
=-\V_{\ezp_M}(\zm_a \xd x^a)
                                                $$
(cf. \Ref{4.1} and \Ref{7.1}).
        \item It follows from (1) and \Ref{7.1} that
                $$ {\ss \zL}^\*_{P_M} (\zi(X)\ezp_M^\*(\zm)) =
\zi(X){\ss \zL}^\*_{P_M} (\ezp_M^\*(\zm)) = \zi(X) \V_{\ezp_M} =
-\I(\zm\otimes X).
                                $$
        \item
                $${\ss \zL}^\*_{P_M} (\xd(\zi(X)\ezp_M^\*(\zm))) = [P,
{\ss \zL}^\*_{P_M} (\zi(X)\ezp_M^\*(\zm))] = -[P, \I(\zm \otimes X)] =
-\H(\zm\otimes X ).
                                        $$
        \endlist
        \endproof

        The next theorems provide us with  homomorphisms of graded
Lie brackets related to certain cotangent lifts.
        \claim{Theorem}{}
        The mapping
                $$ \i \colon \zF^{}_{1}(\ezp_M) \ni \zm \otimes X
\mapsto \i(\zm \otimes X) = \zi(X) \ezp_M^\* (\zm) \in
\zF^{}_{}(\ezp_{\ssT ^\s* M})
                                        $$
        defines an injective homomorphism of the Fr\"olicher
-Nijenhuis bracket of vector-valued forms on $M$ into the extended Poisson
bracket $\{\,,\,\}$ of differential forms on $\sT ^\* M$ defined by the
Poisson structure $P_M$ (cf. \Ref{2.7}).
        \endclaim
        \proof
        Let  $\zm \otimes X ,\ \zn \otimes Y\ \in
\zF^{}_{1}(\ezp_M)$. We can always take $\zm,\zn$ such that
$\xd\zm = \xd \zn =0$, since such vector-valued forms span
$\zF^{}_{1}(\ezp_M)$. In such cases,
                $$ [\zm \otimes X, \zn \otimes Y]^{F-N} = \zm \wedge \zn
\otimes [X,Y] + \zm \wedge \ll_X \zn \otimes Y - \ll_Y\zm \wedge
\zn \otimes X
                                        \tag \label{7.16}$$
        and, locally,
                $$\zm = \xd \zm_1 \wedge \cdots \wedge \xd \zm_k, \
\zn = \xd \zn_1 \wedge \cdots \wedge \xd \zn_l,\ \ \zm_i, \zn_j
\in C^\infty(M).
                                        $$
        Let us denote $g_i = \ezp^\*_M(\zm_i)$, $f_j =
\ezp^\*_M(\zn_j)$ and $g_0= \zi(X),\ f_0 = \zi(Y) $.  Then
                $$ \i(\zm\otimes X) = g_0\xd g_1 \wedge \cdots \wedge
\xd g_k, \ \ \i(\zn\otimes Y) = f_0\xd f_1 \wedge \cdots \wedge
\xd f_l.
                                        $$
        Since $\{f_j,g_i\} = 0$ for $i,j>0$, we get, using \Ref{2.8},
                $$\multline
        \{\i(\zm\otimes X), \i( \zn\otimes Y)\} = \{\zi(X),
\zi(Y)\} \xd g_1 \wedge \cdots \wedge \xd g_k \wedge \xd f_1
\wedge \cdots \wedge \xd f_l  + \\
        + \zi(X)\sum_{i>0}(-1)^i \xd\{\zi(Y), g_i\} \wedge \xd g_1
\wedge \cdots  \wedge \widehat{\xd g_i} \wedge  \cdots \wedge
\xd g_k \wedge \xd f_1  \wedge \cdots \wedge \xd f_l  -\\
        -(-1)^k \zi(Y)\sum_{j>0}(-1)^j \xd\{\zi(X), f_j\} \wedge
\xd g_1 \wedge \cdots \wedge \xd g_k \wedge \xd f_1 \cdots
\wedge \widehat{\xd f_j} \wedge \wedge \cdots \wedge \xd f_l =
\\      = \zi([X,Y]) \ezp_M^\*(\zm \wedge \zn) - \sum_i \zi(X) \xd
g_1 \wedge \cdots \wedge \xd \{\zi(Y), g_i\} \wedge \cdots
\wedge \xd g_k \wedge \ezp_M^\* (\zn) +\\
        +\sum_j \zi(Y) \ezp_M^\*(\zm)\wedge \xd f_1 \wedge \dots
\wedge \xd \{\zi(X), f_j\} \wedge \cdots \wedge \xd f_l.
                \endmultline    \tag\label{7.17}        $$
        Since
                $$ \xd \{\zi(Y), g_i\} = \xd\{\zi(Y), \ezp_M^\*
(\zm_i) \} = \xd (\ezp^\*_M \ll_Y\zm_i) = \ezp_M^\* \ll_Y(\xd
\zm_i)\ \ \text{etc.},
                                        $$
        we finally get from \Ref{7.17}
                $$ \multline
        \{\i(\zm\otimes X), \i( \zn\otimes Y)\} =
\zi([X,Y])\ezp_M^\* (\zm \wedge \zn) - \zi(X) \ezp_M^\*(
\ll_Y\zm \wedge \zn) + \zi(Y)\ezp_M^\*(\zm \wedge \ll_X \zn) =\\
        = \i ([\zm\otimes X, \zn\otimes Y]^{F-N}).
                                \endmultline    $$
Thus  $\i$ defines a homomorphism of the brackets.
The injectivity follows directly from the form of $\i$ in local
coordinates.
        \endproof

        {\bf Remark.} It is well known that the mapping
                $$\gather  \zi \colon \ss S(\ezt_M) \rightarrow
\zF^{}_{}(\ezp_{\smT^\s*} M),\\
        \zi(X_1\vee\dots \vee X_k) = \zi(X_1)\cdots \zi(X_k)
\text{\ for\ } X_i \in \zF^1(\ezt_M),
                        \endgather      $$
        is an injective homomorphism of the symmetric Schouten
bracket on $M$ into the Poisson bracket on $\sT^\s* M$. Thus we
can consider $(\zF^{}_{}(\ezp_{\smT^\s* M}), \{\,,\,\})$ -- the
extended Poisson bracket of differential forms on $\sT^\* M$ --
as a common generalization of the Fr\"olicher-Nijenhuis and the
symmetric Schouten bracket (cf. [7]).
        \claim{Theorem}{}
        \list
        \item The mapping
                $$ \G\colon \zF^{}_{1}(\ezp_M) \rightarrow
\zF^{}_{1}(\ezp_{\smT^\s* M}) \colon \zm \otimes X
\mapsto \sH_{P_M}(\i(\zm \otimes X))
                                        $$
        is an injective homomorphism of the Fr\"olicher-Nijenhuis
bracket on $M$ and the Fr\"olicher-Nijenhuis bracket on $\sT^\*
M$.
        \item The mapping
                $$\H \colon \zF^{}_{1}(\ezp_M) \rightarrow
\zF^{}_{}(\ezt_{\smT^\s* M})
                                        $$
        is an injective homomorphism of the Fr\"olicher-Nijenhuis
bracket on $M$ into the Schouten bracket on $\sT^\* M$.
        \endlist
        \endclaim
        \proof
        (1)\quad  By Theorem 7~b and Theorem~23,  $\G$ is a
homomorphism of the Fr\"olicher-Nijen\-huis brackets. The local form of
$\G$ shows that it is injective:
                $$      \multline
        \G\left(\sum_{a_1<\dots <a_k,a} f^a_{a_1\dots a_k}\xd
x^{a_1} \wedge  \cdots \wedge \xd x^{a_k}\otimes \partial _{x^a}
\right) =\\
        = \sum_{a_1<\dots <a_k,a} f^a_{a_1\dots a_k} \left(\xd
x^{a_1} \wedge  \cdots \wedge \xd x^{a_k}\otimes \partial _{x^a}
- \sum_i(-1)^i \xd p_a \wedge \xd x^{a_1} \wedge  \cdots \wedge
\widehat{\xd x^{ a_i}} \wedge \cdots \right. \\
        \left. \cdots \wedge  \xd x^{a_k} \otimes \partial
_{p_{a_i}} \phantom{\sum} \right)
        - \sum_{a_1<\dots <a_k,a,b} \frac{\partial
f^a_{a_1\dots a_k}}{\partial x^b} p_a \left( \phantom{\sum}\xd
x^{a_1} \wedge
\cdots \wedge \xd x^{a_k}\otimes \partial _{p_b} + \right. \\
        \left. \sum_i(-1)^i \xd x^b \wedge \xd x^{a_1} \cdots
\wedge \widehat{\xd x^{ a_i}} \wedge \xd x^{a_k} \otimes
\partial _{p_{a_i}} \right) .
                                \endmultline
                \tag\label{7.18}        $$

        (2) \quad We have already proved that $\H$ is a homomorphism
(Theorem~18, but it follows also from Theorem~7c and Theorem 23).
        Locally, we have

        $$ \multline
        \H\left( \sum_{a_1<\dots <a_k,a} f^a_{a_1\dots a_k}\xd
x^{a_1} \wedge  \cdots \wedge \xd x^{a_k}\otimes \partial
_{x^a}\right) = \\
        = -\sum_{a_1<\dots <a_k,a} f^a_{a_1\dots a_k}\partial
_{x^a} \wedge \partial _{p_{a_1}} \wedge  \cdots \wedge \partial
_{ p_{a_k}} +  \sum_{a_1<\dots <a_k,a,b} \frac{\partial
f^a_{a_1\dots a_k}}{ \partial x^b} p_a\partial _{p_b} \wedge
\partial _{p_{a_1}} \wedge  \cdots \wedge \partial _{ p_{a_k}}
                                \endmultline  $$

        and the injectivity follows.
        \endproof

        {\bf Remarks.}\ The mapping $\G$, extended by the mapping
                $$ \G \colon  \ss  S(\ezt_M) \rightarrow
\zF^{}_{1}(\ezp_{\smT^\s* M}),
                                        $$
        is the ''common generalization'' of the
Fr\"olicher-Nijenhuis and Schouten brackets as  defined by
Dubois-Violette and Michor ([7]). The second part of the
theorem shows that we can define the classical
Fr\"olicher-Nijenhuis  bracket on $M$ by means of the classical
Schouten bracket on $\sT^\* M$: $[K,L]^{F-N}$ is the only
element of $\zF^{}_{1}(\ezp_M)$ such that  $\H([K,L]^{F-N})$ is
the Schouten bracket $[\H(K), \H(L)]$. In particular, for $N\in
\zF^{1}_{1}(\ezp_M)$ we have, that $[N,N]^{F-N} =0$ if and only
$[\H(N), \H(N)] =0$, i.e., $N$ is a Nijenhuis tensor on $M$ if
and only if $\H(N)$ is a Poisson tensor on $\sT^\* M$. It is easy
to verify that the Poisson structure defined by $\H(N)$ is linear
and, consequently, corresponds to an algebroid structure on the
tangent bundle $\ezt_M \colon \sT M \rightarrow M$. This sugests
a new approach to the theory of Poisson-Nijenhuis manifolds
which will be developed in a separate publication.

        \hl1""{Conclusions} We defined natural graded Lie brackets
associated with Lie algebroids, their lifts, and we studied
relations between the brackets and the lifted brackets. All the
constructions seem to be canonical, generalizing some known
results, but we also  obtained interesting results on the
classical level, as for example embeddings of the
Nijenhuis-Richardson and the
Fr\"olicher-Nijenhuis bracket over $M$ into the Schouten bracket
over $\sT^\* M$. However, the analysis presented in this paper
is far from being complete. We did not discuss the canonical
homomorphism between tangent Lie algebroid on  $\sT \zt \colon
\sT E \rightarrow \sT M $ and the canonical algebroid  on
$\ezt_E \colon \sT E \rightarrow E$.  We skipped also the
problem of iterated lifts. In Section 3, we mentioned that some
iterated lifts of an algebroid are canonically isomorphic.   The
existence of these isomorphisms provokes  the question, whether
iterated lifts of sections can be identified. For example, one
could expect that $\H \circ \T (K) $ can be identified with
$\T\circ \H (K)$.  In a forthcoming publication we shall discuss
identities of this kind, as well as other observations and
applications of the presented results.

\makebib

\bib \key{\bf 1 } \by J.~V.~Beltr\'an and J.~Monterde   \pp 383--402
        \paper Graded Poisson structures on the algebra of differential forms
        \yr 1995  \vol 70
        \jour Comment. Math. Helv.  \endbib

\bib \key{\bf 2} \by K. H. Bhaskara and K. Viswanath   \pp 68--72
        \paper Calculus on Poisson manifolds
        \yr  1988 \vol 20
        \jour Bull. London. Math. Soc.  \endbib

\bib \key{\bf 3} \by A.~Coste, P.~Dazord and A.~Weinstein
\pp 1--65 \paper Groupo\"{\char'020}des symplectiques
        \yr 1987        \jour Publications du D\'epartement de
Math\'ematiques de l'Universit\'e Claude-Bernard de Lyon  \endbib

\bib \key{\bf 4} \by T. J. Courant  \pp 631--661
        \paper Dirac manifolds
        \yr 1990  \vol 319
        \jour Trans. Amer. Math. Soc.  \endbib

\bib \key{\bf 5} \by T. J. Courant  \pp 5153--5160
        \paper Tangent Dirac Structures
        \yr  1990 \vol 23
        \jour  J. Phys. A \endbib
\bib \key{\bf 6} \by A. Cabras and A. M. Vinogradov   \pp 75--100
        \paper Extensions of the Poisson bracket to differential
forms and multi-vec\-tor fields
        \yr  1992 \vol 9
        \jour  J. Geom. Phys. \endbib

\bib \key{\bf 7} \by M. Dubois-Violette and P. Michor  \pp 51--66
        \paper A common generalization of the Fr\"olicher-Nijenhuis
bracket and the Schouten bracket for symmetric multivector fields
        \yr  1995 \vol 6
        \jour Indag. Math. N. S.  \endbib

\bib \key{\bf 8} \by A. Fr\"olicher and A. Nijenhuis   \pp 338--359
        \paper Theory of vector-valued differential forms, Part 1
        \yr 1956  \vol 18
        \jour  Indag. Math. \endbib

\bib \key{\bf 9 } \by J. Grabowski
        \paper  Z-graded extensions of Poisson brackets
        \jour \toappear \endbib

\bib \key{\bf 10} \by J. Grabowski and P. Urba\'nski   \pp 6743--6777
        \paper Tangent lifts of Poisson and related structures
        \yr  1995 \vol 28
        \jour J. Phys. A  \endbib

\bib \key{\bf 11} \by J. Huebschmann   \pp 57--113
        \paper Poisson cohomology and quantization
        \yr  1990 \vol 408
        \jour J. Reine Angew. Math.  \endbib

\bib \key{\bf 12} \by M. V. Karasev   \pp 497--527
        \paper Analogues of the objects of Lie group theory for
nonlinear Poisson brackets
        \yr 1987  \vol 28
        \jour Mat. USSR - Izv.  \endbib

\bib \key{\bf 13 } \by I. Kol\'a\v r, P. Michor and J. Slovak
        \book Natural Operations in Differential Geometry
        \yr 1991 \publ Springer-Ver\-lag
         \endbib
\bib \key{\bf 14} \by Y. Kosmann-Schwarzbach   \pp 153--165
        \paper Exact Gerstenhaber algebras and Lie bialgebroids
        \yr  1995 \vol 41
        \jour Acta Applicandae Math.  \endbib

\bib \key{\bf 15} \by Y. Kosmann-Schwarzbach   \paperinfo preprint
        \paper From Poisson algebras to Gerstenhaber algebras
        \yr 1995  \endbib

\bib \key{\bf 16} \by Y. Kosmann-Schwarzbach and F. Magri \pp
35--81
        \paper Poisson-Nijenhuis structures
        \yr 1990  \vol 53
        \jour Ann. Inst. Henri Poin\-ca\-r\'e, Phys. Th\'eor.  \endbib
\bib \key{\bf 17} \by  J.-L. Koszul
        \paper Crochet de Schouten\-Nijenhuis et cohomologie\
        \yr 1985  \inbook Elie Cartan et les math\'ematiques d'aujour\-d'hui
        \publ Ast\'erisque hors s\'erie   \endbib

\bib \key{\bf 18} \by  I. Krasilshchik \pp 79--100
        \paper Schouten bracket and canonical algebra, II
        \yr1988   \vol 1334 \inbook Global Analysis
        \publ Lect. Notes Math  \endbib

\bib \key{\bf 19} \by  I. Krasilshchik  \pp 70--76
        \paper Supercanonical algebras and Schouten brackets
        \yr  1991 \vol 49
        \jour Mathematical Notes  \endbib

\bib \key{\bf 20} \by K. C. H.  Mackenzie
        \book Lie Groupoids and Lie Algebroids in Differential
Geometry \bookinfo London Math. Soc. Lecture Notes Ser. 124
        \yr1987   \publ Cambridge University Press
           \endbib

\bib \key{\bf 21} \by K. C. H.  Mackenzie  \pp 180-239
        \paper Double Lie algebroids and second-order geometry I
        \yr 1992  \vol 94
        \jour Adv. Math  \endbib

\bib \key{\bf 22} \by K. C. H.  Mackenzie and P. Xu  \pp 415--452
        \paper Lie bialgebroids and Poisson groupoids
        \yr 1994  \vol 73
        \jour  Duke Math. J. \endbib

\bib \key{\bf 23} \by P. Michor  \pp 67--82
        \paper A generalization of Hamiltonian mechanics
        \yr 1985  \vol 2
        \jour J. Geom. Phys.  \endbib

\bib \key{\bf 24} \by  P. Michor  \pp 207--215
        \paper Remarks on the Schouten-Nijenhuis bracket
        \yr 1987  \vol 16
        \jour Suppl. Rend. Circ. Mat. di Palermo, Serie II  \endbib

\bib \key{\bf 25} \by A. Nijenhuis  \pp 390--403
        \paper Jacobi-type identities for bilinear differential
concomitants of certain tensor fields I
        \yr  1955 \vol 17
        \jour Indag. Math.   \endbib

\bib \key{\bf 26} \by A. Nijenhuis and R. Richardson  \pp 89--106
        \paper Deformations of Lie algebra structures
        \yr  1967 \vol 171
        \jour J. Math. Mech.  \endbib

\bib \key{\bf 27} \by  J. Pradines
        \paper Fibr\'es vectoriels doubles et calcul des jets non
holonomes
        \yr 1974  \paperinfo notes polycopi\'ees, Amiens
           \endbib

\bib \key{\bf 28} \by  J, A. Schouten \pp 449--452
        \paper \"Uber Differentialkonkomitanten zweier
kontravarianter Gr\"o\char'031 en
        \yr  1940 \vol 2
        \jour Indag. Math.  \endbib

\bib \key{\bf 29} \by S. M. Skryabin   \pp 1229--1336
        \paper An algebraic approach to the Lie algebra of Cartan type
        \yr 1993  \vol 21
        \jour Comm. Algebra  \endbib

\bib \key{\bf 30} \by W. M. Tulczyjew
        \book Geometric Formulations of Physical Theories
        \publ Bibliopolis, Napoli \yr 1989
         \endbib

\bib \key{\bf 31} \by K. Yano and S. Ishihara
        \book Tangent and Cotangent Bundles
          \publ Marcel Dekker, Inc., New York \yr 1973
          \endbib

        \enddocument